\documentclass[12pt,reqno]{amsart}
\usepackage[margin=1in]{geometry}
\usepackage{amsmath,amssymb,amsthm,graphicx,amsxtra, setspace}
\usepackage[utf8]{inputenc}
\usepackage{mathrsfs}
\usepackage{hyperref}
\usepackage{upgreek}
\usepackage{mathtools}
\usepackage{xcolor}
\usepackage{stmaryrd}
\usepackage[noadjust]{cite}

\usepackage{multirow}
\usepackage{booktabs}
\usepackage{caption}
\usepackage{subcaption}

\usepackage[mathcal]{euscript}
\allowdisplaybreaks

\usepackage{float}
\usepackage{multirow}

\usepackage{algorithm}
\usepackage{algorithmic}

\usepackage[pagewise]{lineno}

\DeclareMathAlphabet{\mathpzc}{OT1}{pzc}{m}{it}

\usepackage[cyr]{aeguill}

\colorlet{darkblue}{blue!50!black}

\hypersetup{
	colorlinks,%
	citecolor=blue,%
	filecolor=red,%
	linkcolor=red,%
	urlcolor=blue,%
	pdfnewwindow=true,%
	pdfstartview={FitH}
}

\newtheorem{theorem}{Theorem}[section]
\newtheorem{lemma}[theorem]{Lemma}
\newtheorem{proposition}[theorem]{Proposition}

\newtheorem{problem}[theorem]{Problem}
\newtheorem{property}[theorem]{Property}
\newtheorem{example}[theorem]{Example}
\newtheorem{remark}[theorem]{Remark}

\renewcommand{\arraystretch}{1.2}

\allowdisplaybreaks

\let\originalleft\left
\let\originalright\right
\renewcommand{\left}{\mathopen{}\mathclose\bgroup\originalleft}
\renewcommand{\right}{\aftergroup\egroup\originalright}

\let\emptyset\varnothing

\renewcommand{\d}{\/\mathrm{d}\/}

\def\S{\mathrm{S}}

\def\L{\mathbb{L}}

\def\f{\boldsymbol{f}}

\def\B{\mathcal{B}}

\def\bx{\boldsymbol{x}}

\def\bz{\boldsymbol{z}}
\def\bv{\boldsymbol{v}}
\def\V{\mathbb{v}}
\def\bw{\boldsymbol{w}}
\def\f{\boldsymbol{f}}

\def\N{\mathbb{N}}

\def\V{\mathbb{V}}

\def\P{\mathrm{P}}
\def\bu{\boldsymbol{y}}
\def\H{\mathbb{H}}
\def\bn{\boldsymbol{n}}

\newcommand{\R}{\mathbb{R}}

\renewcommand{\d}{\/\mathrm{d}\/}

\newcommand{\tnorm}[1]{\lvert\!\lvert\!\lvert #1 \rvert\!\rvert\!\rvert}

\newcommand{\Addresses}{{
		\footnote{
			
			\noindent \textsuperscript{1,3}Department of Mathematics, Indian Institute of Technology Roorkee-IIT Roorkee,
			Haridwar Highway, Roorkee, Uttarakhand 247667, INDIA.\par\nopagebreak
			\textsuperscript{2} Department of Mathematics, University of Iowa, Iowa City, IA 52242, USA. 
 			\noindent 
			
			\textit{e-mail:} \texttt{Wasim Akram: wakram2k11@gmail.com.}
			
			 \textit{e-mail:} \texttt{Weimin Han: weimin-han@uiowa.edu.}
			
			 \textit{e-mail:} \texttt{Manil T. Mohan: maniltmohan@ma.iitr.ac.in, maniltmohan@gmail.com.}

			\noindent \textsuperscript{*}Corresponding author.
			
			\textit{Key words:} Stationary convective Brinkman--Forchheimer Extended Darcy equations, hemivariational inequality, pseudomonotonicity, mixed finite element method, error estimates. 
			
			Mathematics Subject Classification (2020): 65N30, 35Q35, 47J20, 65J15, 65N15.
			
}}}

\begin{document}
	
	\title[Discontinuous Galerkin Methods for 2D and 3D CBFeD equations]{Discontinuous Galerkin Methods for Stationary  Convective Brinkman--Forchheimer Extended Darcy Equations with a Nonlinear Slip Boundary Condition of Friction Type
		\Addresses}
	\author[W. Akram, W. Han and M. T. Mohan]
	{Wasim Akram\textsuperscript{1}, Weimin Han\textsuperscript{2*}
		 and Manil T. Mohan\textsuperscript{3}}
	
	\maketitle

	\begin{abstract}
	This work develops and analyzes several discontinuous Galerkin methods for a variational inequality associated with the stationary convective Brinkman--Forchheimer extended Darcy (CBFeD) system. The CBFeD model generalizes the incompressible Navier--Stokes equations by incorporating additional damping and pumping terms. The resulting variational inequality characterizes the motion of a viscous, incompressible fluid in a saturated porous medium under a  friction-type slip boundary condition. For the CBFeD variational inequality, we first  explore the existence of a unique weak solution. Then, we establish the existence, uniqueness, and stability of the numerical solutions obtained using discontinuous Galerkin methods. Optimal error estimates are derived, measured in the broken $H^1$-norm for the velocity and in the $L^2$-norm for the pressure, when linear velocity elements and piecewise constant pressure approximations are employed. Numerical simulations are carried out to confirm the theoretically predicted convergence orders and to highlight key strengths of the method, including its ability to capture discontinuities, effectively resolve shear layers, handle advection-dominated regimes, and operate robustly on general polygonal meshes.
	\end{abstract}

	\section{Introduction}\label{sec1}\setcounter{equation}{0}
\subsection{The model}\label{subsec:model} Let $\mathcal{O}$ be a bounded domain in $\mathbb{R}^d$, where $d\in\{2,3\}$, and assume that its boundary $\Gamma=\partial\mathcal{O}$ is Lipschitz continuous. The stationary incompressible convective Brinkman--Forchheimer extended Darcy (CBFeD) problem in $\mathcal{O}$ is governed by the following momentum and continuity equations:
	\begin{equation} \label{eqn-CBF}
		\left\{
		\begin{aligned}
			-\mu\Delta\bu + (\bu \cdot \nabla)\bu + \alpha\bu + \beta|\bu|^{r-1}\bu - \kappa |\bu|^{q-1}\bu + \nabla p &= \f, &&\text{in } \mathcal{O}, \\
			\mathrm{div}\,\bu &= 0, &&\text{in } \mathcal{O}.
		\end{aligned}
		\right. 
	\end{equation}
	The velocity field is denoted by $\bu : \mathcal{O} \to \mathbb{R}^d$, the pressure by $p : \mathcal{O} \to \mathbb{R}$, and the external force by $\f : \mathcal{O} \to \mathbb{R}^d$.  The coefficient $\mu > 0$ represents the Brinkman term, which models the effective viscosity of the medium. The non-negative constant $\alpha$ corresponds to the Darcy and constant $\beta\ge 0$ corresponds to the Forchheimer damping terms, accounting for flow resistance induced by permeability and porosity, respectively.
	To incorporate possible pumping effects, an additional nonlinear term of the form $\kappa |\bu|^{q-1}\bu$ is introduced. When  $\kappa \ge 0$, this term counteracts the damping contribution, which is the case considered in this work. The exponent $r \in [1,\infty)$, referred to as the \emph{absorption exponent}, determines the intensity of the Forchheimer damping, with $r=3$ identified as the \emph{critical value} (\cite[Proposition 1.1]{KWH}). The parameter $q \in [1,r)$ regulates the strength of the pumping mechanism.
	Special cases of the model include the classical Navier--Stokes equations, obtained when $\alpha = \beta = \kappa = 0$, and the damped Navier--Stokes system, recovered for $\alpha,\beta > 0$ with $\kappa = 0$. The CBFeD formulation originates from an extended Darcy-Forchheimer law, as discussed in \cite{SGKKMTM, MTT}. Finally, the divergence-free constraint in the continuity equation ensures the incompressibility of the flow.

	Associated with the velocity field $\bu$, the rate-of-deformation tensor is defined by
		$$
		\boldsymbol{\varepsilon}(\bu) = \frac{1}{2}\big(\nabla \bu + (\nabla \bu)^\top\big),
		$$
		which represents the symmetric part of the velocity gradient and describes the fluid strain rate. Note that 
		$$\mathrm{div}(	\boldsymbol{\varepsilon}(\bu)) = \frac{1}{2}\left(\Delta \bu + \nabla(\mathrm{div\ }\bu)\right),$$ and the incompressibility condition provides 
		$2\mathrm{div}(	\boldsymbol{\varepsilon}(\bu) )= \Delta \bu .$ The corresponding Cauchy stress tensor is defined by
		\[
		\mathbb{T}(\bu,p)=2\mu\varepsilon(\bu)-p\mathrm{I},
		\]
		where \(\mathrm{I}\) denotes the identity tensor. The associated traction (or stress) vector acting on the boundary is
		\[
		\boldsymbol{t}(\bu,p)=\mathbb{T}(\bu,p)\boldsymbol{n},
		\]
		where \(\boldsymbol{n}=(n_1,\ldots,n_d)^{\top}\) is the unit outward normal vector to \(\Gamma\).  
		The tangential component of the traction vector is defined by $	\boldsymbol{t}_{\boldsymbol{\tau}}	=	\boldsymbol{t}	-	(\boldsymbol{t}\cdot \boldsymbol{n})\boldsymbol{n}.	$

	To complete the formulation of the CBF system \eqref{eqn-CBF}, suitable boundary conditions are imposed.
We decompose the boundary $\Gamma$ into two disjoint, non-empty components $\Gamma_{\mathrm{D}}$ and $\Gamma_{\mathrm{S}}$ such that $\Gamma = \overline{\Gamma}_{\mathrm{D}} \cup \overline{\Gamma}_{\mathrm{S}}$, $\Gamma_{\mathrm{D}} \cap \Gamma_{\mathrm{S}} = \emptyset$ and the Lebesgue measure of $\Gamma_{\mathrm{D}}$ denoted by $|\Gamma_{\mathrm{D}}|>0$. On the portion $\Gamma_{\mathrm{D}}$, we consider the homogeneous Dirichlet boundary condition given by:
\begin{align}
	\bu=\boldsymbol{0}\ \text{ on }\ \Gamma_{\mathrm{D}}. 
\end{align}
On the boundary segment $\Gamma_{\mathrm{S}}$, we distinguish between the normal and tangential directions. Let $\boldsymbol{n}=(n_1,\ldots,n_d)^{\top}$ be the unit outward normal vector along $\Gamma_{\mathrm{S}}$.
For any vector field $\boldsymbol{v}$ defined on the boundary, we decompose it into its normal and tangential components as
$
v_{\boldsymbol{n}} = \boldsymbol{v}\cdot \boldsymbol{n}, \
\boldsymbol{v}_{\boldsymbol{\tau}} =\boldsymbol{v}-(\boldsymbol{v}\cdot\boldsymbol{n})\boldsymbol{n}.
$
Assuming that $\bu$ is smooth up to the boundary, following the convention adopted in the mathematical literature on the Tresca slip boundary condition (\cite{FJWHWY}), we define the tangential stress by
\begin{align}\label{eqn-boundary-0}
\boldsymbol{\sigma}_{\boldsymbol{\tau}}(\bu) = \mu \frac{\partial \bu_{\boldsymbol{\tau}}}{\partial \boldsymbol{n}}=\mu(\boldsymbol{n}\cdot\nabla)\bu_{\boldsymbol{\tau}},
\end{align}
where $\mu$ is the Brinkman coefficient.
On $\Gamma_{\mathrm{S}}$, we impose a \emph{friction-type slip and non-leak boundary condition} of the type
\begin{align}\label{eqn-boundary}
	y_{\boldsymbol{n}}= 0,\  |\boldsymbol{\sigma}_{\boldsymbol{\tau}}| \leq g,\ \boldsymbol{\sigma}_{\boldsymbol{\tau}}\cdot \bu_{\boldsymbol{\tau}} + g|\bu_{\boldsymbol{\tau}}| =0\ \text{ on }\ \Gamma_{\mathrm{S}}.
\end{align}
The function $g : \Gamma_{\mathrm{S}} \to [0,\infty)$ is referred to as the \emph{threshold slip} or \emph{barrier function}. When $g \equiv 0$, condition \eqref{eqn-boundary} simplifies to the classical slip boundary condition, namely
$
y_{\boldsymbol{n}} = 0, \  \boldsymbol{\sigma}_{\boldsymbol{\tau}} = \boldsymbol{0}.
$
Moreover, the second and third relations in \eqref{eqn-boundary}  can be equivalently expressed in the form of the following logical implications:
\begin{align}\label{eqn-boundary-1}
\begin{cases}
	|\boldsymbol{\sigma}_{\boldsymbol{\tau}}| < g \;\; \Rightarrow \;\; \bu_{\boldsymbol{\tau}} = \boldsymbol{0}, \\[0.5em]
	\bu_{\boldsymbol{\tau}} \neq \boldsymbol{0} \;\; \Rightarrow \;\; \boldsymbol{\sigma}_{\boldsymbol{\tau}} = -g \dfrac{\bu_{\boldsymbol{\tau}}}{|\bu_{\boldsymbol{\tau}}|}.
\end{cases}
\end{align}
	This type of frictional boundary condition was first introduced by Fujita \cite{HFu} in the context of modeling blood flow in veins affected by arterial sclerosis, as well as fluid flow through canals whose bottoms are covered with layers of mud and pebbles.
	
	The numerical solution of problem \eqref{eqn-CBF}-\eqref{eqn-boundary} is challenging due to the presence of strong nonlinearities, the coupling between velocity and pressure, and the slip boundary condition formulated as an inequality.
	
	\subsection{Literature review}
	Since the early 1990s, boundary value problems with nonsmooth slip or leak friction-type boundary conditions have been widely studied following the pioneering work of Fujita \cite{HFu, HFu1}.  Since then, the mathematical analysis of the Stokes and Navier--Stokes equations with nonlinear slip boundary conditions has been extensively studied; see, for example, \cite{HFu1, HFuHK, HFu2, HFu3, YLKTL, NSa}. Numerical approximation of these problems has also been investigated using finite element methods and Uzawa-type iterative algorithms; see \cite{JKDJK, TKa, TKa1, YLRA, YLRA1, YLRA2, YLKL}.
	
	Hemivariational inequalities (HVIs) were introduced by Panagiotopoulos in the early 1980s \cite{PDP} to model engineering problems involving nonsmooth and nonmonotone phenomena. Since then, the theory of hemivariational inequalities and variational--hemivariational inequalities has undergone significant development, both in terms of mathematical analysis and numerical approximation; see, for example, \cite{WAMTM1, VGPR, VGBR, WHKCFJ, WHFJYY, WHSM1, WHSM, WHHQLM, WHCW, FJWHGZ, MLWHSZ, MLFW, ZNPDP, PDP1, WHMS, WHMS1, WHMS2}. In the context of fluid mechanics, HVIs provide a natural framework for modelling incompressible flows subject to frictional slip boundary conditions. In particular, hemivariational inequality formulations for the stationary Stokes and Navier--Stokes equations have been extensively studied, with existence, uniqueness, stability, and finite element approximation results established in \cite{WHKCFJ}. On the numerical side, discontinuous Galerkin methods for the stationary Navier--Stokes equations with nonlinear friction-type slip boundary conditions were developed and analyzed in \cite{FJWHWY}, while analogous methods for the Stokes equations were investigated in \cite{FWMLWH}, where stability and optimal-order error estimates were obtained.

	Recently, the incompressible convective Brinkman--Forchheimer (CBF) equations \cite{SGMTM, JJSGMTM2, KKMTM-2025, KKMTM-RMP, MTM9} (or the Navier--Stokes equations with damping \cite{WHHQLM,WWXLCWH,SMSDYZ}) and the CBFeD equations \cite{WAMTM1, WAMTM2, SGKKMTM, JJSGMTM1,  MTM7, MTM8} have attracted considerable attention owing to their wide range of applications in fluid flow through porous media. The papers \cite{WHHQLM,WWXLCWH} established the well-posedness of Stokes and stationary Navier--Stokes hemivariational inequalities with damping and frictional slip boundary conditions, developed mixed finite element methods with optimal error estimates, and proposed efficient iterative solution algorithms. The paper \cite{SMSDYZ} developed a penalty method for quasi-variational and quasi-variational--hemivariational inequalities with implicit constraints and applied the theory to the incompressible Navier--Stokes equations with damping and multivalued nonmonotone slip boundary conditions. In the hemivariational inequality setting, \cite{MTM7} established the well-posedness of the stationary CBFeD model with nonmonotone frictional slip boundary conditions, together with equivalent minimization formulations and an iterative solution algorithm. Subsequently, \cite{WAMTM1} developed and analyzed a mixed finite element method for the stationary CBFeD hemivariational inequality, proving optimal orders of convergence, while \cite{WAMTM2} investigated an optimal control problem for the same model, establishing the existence of optimal controls and the convergence of finite element approximations. The non-stationary CBFeD hemivariational inequality was studied in \cite{JJSGMTM1}, where the existence, uniqueness, and energy equality of weak solutions were established in both two and three spatial dimensions. Furthermore, \cite{JJMTM} considered a quasi-variational--hemivariational inequality for the CBFeD equations describing Bingham fluids and proved the existence of weak solutions. Motivated by these developments, the present work investigates the discontinuous Galerkin approximation of the stationary CBFeD hemivariational inequality, providing a new numerical framework for this class of nonsmooth fluid flow problems.
	
	\subsection{Our contributions and novelties}
	This paper develops and analyzes a class of discontinuous Galerkin methods for the (hemi)-variational inequality associated with the stationary CBFeD equations.
	The additional Darcy damping and nonlinear Forchheimer resistance terms significantly increase the mathematical complexity of both the continuous and discrete analyses.  These extensions render the CBFeD model particularly effective for modeling complex flows in porous media, including groundwater transport, filtration phenomena, and biological fluid flows. The model is supplemented with slip boundary conditions described by a nonsmooth, possibly nonconvex friction law, permitting tangential slip of the fluid along the boundary. The resulting nonsmooth boundary interaction is naturally formulated as a hemivariational inequality, where the friction law is represented by a locally Lipschitz, generally nonconvex superpotential defined on the slip boundary.
	
	The main contributions and novelties of this paper are summarized as follows:
	\begin{itemize}
\item 	We formulate a variational inequality associated with the stationary CBFeD equations under a friction-type slip boundary condition, see Problems~\ref{problem-1}--\ref{problem-3}. The formulation naturally incorporates the incompressibility constraint together with the nonlinear boundary condition arising from frictional slip. We establish the well-posedness of the continuous variational inequality by proving the existence and uniqueness of weak solutions under suitable assumptions in Theorems \ref{thm-exists-unique} and \ref{thm-equivalent}. This analytical result provides the foundation for the subsequent numerical analysis.
	
\item 	We develop several discontinuous Galerkin formulations (see Section~\ref{sec4}) for the CBFeD variational inequality. The proposed schemes extend existing DG methods for the Navier--Stokes variational inequalities to the considerably more challenging CBFeD model containing both Darcy and Forchheimer nonlinearities. In Theorem \ref{thm-dgstab}, we prove the existence, uniqueness, and stability of the discontinuous Galerkin approximations. The analysis requires new estimates to control the additional nonlinear operators introduced by the CBFeD equations while preserving the stability of the discrete problem. Under suitable regularity assumptions, we derive optimal a priori error estimates in the broken $H^1$-norm for the velocity and the $L^2$-norm for the pressure when piecewise linear velocity elements and piecewise constant pressure elements are employed, see Theorem~\ref{thm-dgerr}.
	
\item 	Finally, numerical simulations are presented in Section~\ref{sec7} to verify the theoretical convergence orders and to assess the performance of the proposed methods in capturing discontinuities caused by the friction-type boundary condition, resolving shear layers, handling advection-dominated flows, and performing robustly on general polygonal meshes.
\end{itemize}
	
	\subsection{Organization of the paper} The rest of the paper is organized as follows. In the next section, we introduce the mathematical and functional framework required for the analysis of the model \eqref{eqn-CBF} presented in Subsection~\ref{subsec:model}. After introducing the required preliminaries, we formulate the variational inequality associated with the system governed by equations \eqref{eqn-CBF}--\eqref{eqn-boundary-1}, thereby obtaining Problems \ref{problem-1} (in mixed formulation), \ref{problem-2} (variational inequality with pressure), and \ref{problem-3} (without pressure). We conclude the section by showing that Problems \ref{problem-1} and \ref{problem-2} are equivalent. In Subsection \ref{subsec:regprob}, we construct a Yosida regularization of the nondifferentiable boundary functional to obtain a smooth approximation of the original problem. Subsection \ref{subsec:existweak} establishes the existence of weak solutions to the regularized problem by means of a Galerkin approximation and compactness techniques. In Subsection \ref{subsec:wellp}, the limit as the regularization parameter tends to zero is analyzed, leading to the existence and uniqueness of weak solutions for Problem \ref{problem-3}. Finally, in Subsection \ref{subsec:comWell}, the Ladyzhenskaya--Babuška--Brezzi (LBB) condition is invoked to reconstruct the pressure and prove the existence and uniqueness of solutions to Problem \ref{problem-2}.

Section \ref{sec4} is devoted to the discontinuous Galerkin formulation of the variational inequality. We develop the DG discretization and formulate the corresponding finite-dimensional variational inequality. In Section \ref{sec5}, we investigate the stability and solvability of the discontinuous Galerkin formulation of the variational inequality arising from the CBFeD model. Section \ref{sec6} focuses on the a priori error analysis of the discontinuous Galerkin approximation. We derive optimal error estimates for the velocity in the broken $H^1$--norm and the pressure in the $L^2$--norm.
Finally, Section \ref{sec7} concludes the article by describing the proposed algorithm and presenting several numerical examples to validate the theoretical convergence results and illustrate the performance of the proposed method.
	
		\section{Variational Inequality}\label{sec2}\setcounter{equation}{0} In this section, we present the variational inequality formulation of the problem \eqref{eqn-CBF}-\eqref{eqn-boundary}. 
		\subsection{Functional setting} 
		We begin by introducing some notation. For an integer $m$, we denote by $H^m(\mathcal{O})$, the standard Sobolev space which is defined by 
		\begin{align*}
			H^m(\mathcal{O}):=\left\{v\in L^2(\mathcal{O})\big|\partial^{\boldsymbol{k}}v\in L^2(\mathcal{O}),\ \text{ for all }\ \boldsymbol{k}\ \text{ such that }\  |\boldsymbol{k}|\leq m\right\},
		\end{align*}
		where $\boldsymbol{k}=(k_1,\ldots,k_d)$, $k_1,\ldots k_d$ are non-negative integers, $|\boldsymbol{k}| = k_1 + \cdots+ k_d$, and
		\begin{align*}
			\partial^{\boldsymbol{k}}v=\frac{\partial^{|\boldsymbol{k}|}v}{\partial x_1^{k_1}\cdots\partial x_d^{k_d}}. 
		\end{align*}
		It is a Hilbert space  with the inner product
		\begin{align*}
			\left( u, v\right)_{H^m(\mathcal{O})}:= \sum_{|\boldsymbol{k}|\le m} \left( \partial^{\boldsymbol{k}} u, \partial^{\boldsymbol{k}} v \right).
		\end{align*}
		The norm and the associated seminorm on it are given by 
		\begin{align*}
			\|v\|_{m,\mathcal{O}}:=\left(\sum_{0\leq k\leq m}\int_{\mathcal{O}}|\partial^{\boldsymbol{k}}v(\boldsymbol{x})|^2d\boldsymbol{x}\right)^{1/2},\ 	|v|_{m,\mathcal{O}}:=\left(\sum_{|k|= m}\int_{\mathcal{O}}|\partial^{\boldsymbol{k}}v(\boldsymbol{x})|^2d\boldsymbol{x}\right)^{1/2}. 
		\end{align*}
		For functions satisfying a homogeneous Dirichlet condition on $\Gamma=\partial\mathcal{O},$
		\begin{align*}
			H_0^1(\mathcal{O}):=\left\{v\in H^1(\mathcal{O}):v\big|_{\Gamma}=0\right\},
		\end{align*}
	where 	$v\big|_{\Gamma}$ is understood in the sense of trace. We shall also consider the following space of functions whose average over the domain is zero:
	\begin{align*}
		L^2_0(\mathcal{O}):=\left\{q\in L^2(\mathcal{O}):\int_{\mathcal{O}}q(\boldsymbol{x})d\boldsymbol{x}=0\right\}. 
	\end{align*}
	Let us define $\mathbb{L}^p(\mathcal{O}):=L^p(\mathcal{O};\mathbb{R}^d),$ for $1\leq p\leq \infty$,  $\mathbb{H}^1(\mathcal{O})=H^1(\mathcal{O};\mathbb{R}^d)$ and 
	\begin{align*}
		\mathbb{V}:=\left\{\boldsymbol{v}\in \mathbb{H}^1(\mathcal{O}): \boldsymbol{v}\big|_{\Gamma_{\mathrm{D}}}=\boldsymbol{0},\ v_{\boldsymbol{n}}\big|_{\Gamma_{\mathrm{S}}}=0 \right\},\ Q=L^2_0(\mathcal{O}), \ \mathbb{H}=L^2(\mathcal{O};\mathbb{R}^d).
	\end{align*}
Since $|\Gamma_{\mathrm{D}}| > 0$, Poincar\'e's inequality applies, and thus $\mathbb{V}$ can be endowed with the norm $\|\nabla \cdot \|_{0,\mathcal{O}}$.
 For simplicity, we will omit the domain $\mathcal{O}$ from the notation of norms throughout the remainder of this paper. The scalar product and norm in $Q$ are the standard $L^2(\mathcal{O})$ inner product denoted by $(\cdot,\cdot)$ and its associated norm $\|\cdot\|_0$. Additionally, we define $\mathbb{V}_0:=\mathbb{H}_0^1(\mathcal{O})$, 
	\begin{align*}
		\mathbb{V}_{\mathrm{div}}&:=\left\{\boldsymbol{v}\in\mathbb{V}:\nabla\cdot\boldsymbol{v}=0\ \text{ in }\ \mathcal{O}\right\},\  	\mathbb{V}_{0,\mathrm{div}}:=\left\{\boldsymbol{v}\in\mathbb{V}_0:\nabla\cdot\boldsymbol{v}=0\ \text{ in }\ \mathcal{O}\right\}, \\
		\mathbb{H}_{\mathrm{div}}&:=\left\{\boldsymbol{v}\in\mathbb{L}^2(\mathcal{O}):\nabla\cdot\boldsymbol{v}=0\ \text{ in }\ \mathcal{O}\ \text{ and }\ v_{\boldsymbol{n}}\big|_{\Gamma}=0 \right\}, 
	\end{align*}
where $\bv$ is divergence-free in the sense of distributions, that is, $\int_{\mathcal{O}}\bv(\bx)\cdot\nabla \phi(\bx)d\bx=0$ for all $\phi\in C_c^{\infty}(\mathcal{O})$, and $v_{\boldsymbol{n}}\big|_{\Gamma},$ the restriction of $\bv\cdot\bn$ on $\Gamma$, is understood in the sense of trace (\cite[Theorem 1.2]{Te1}). We use the notation 
	\begin{align*}
		(p,q)&:=\int_{\mathcal{O}}p(\boldsymbol{x})q(\boldsymbol{x})d\boldsymbol{x}, \ \text{ for all }\ p,q\in L^2(\mathcal{O}),\\
		(\boldsymbol{v},\boldsymbol{w})&:=\int_{\mathcal{O}}\sum_{i=1}^dv_i(\boldsymbol{x})w_i(\boldsymbol{x})d\boldsymbol{x}, \ \text{ for all }\ \boldsymbol{v},\boldsymbol{w}\in\mathbb{H},\\
		(\nabla\boldsymbol{v},\nabla\boldsymbol{w})&:=\int_{\mathcal{O}}\sum_{i,j=1}^d\frac{\partial v_i}{\partial x_j}\frac{\partial w_i}{\partial x_j}d\boldsymbol{x}, \ \text{ for all }\ \boldsymbol{v},\boldsymbol{w}\in\mathbb{H}^1(\mathcal{O}).
	\end{align*}
	Moreover, we define 
	\begin{align*}
		a(\bu,\bv)&:=(\nabla\bu,\nabla\bv),\ d(\bv,p):=-(\nabla\cdot\bv,p),\ b(\bu,\bv,\bw):=((\bu\cdot\nabla)\bv,\bw),\nonumber\\
		c_r(\bu,\bv,\bw)&:=(|\bu|^{r-1}\bv,\bw),\ 	c_q(\bu,\bv,\bw):=(|\bu|^{q-1}\bv,\bw),
\end{align*}
for all $\bu,\bv,\bw\in\mathbb{V}\cap\mathbb{L}^{r+1}(\mathcal{O})$ and $p\in Q$. We denote 
\begin{align*}
\mbox{$c_r^0(\bu,\bv)=	c_r(\bu,\bu,\bv)=(|\bu|^{r-1}\bu,\bv)$\  and \ $c_q^0(\bu,\bv)=c_q(\bu,\bu,\bv)=(|\bu|^{q-1}\bu,\bv)$. }
\end{align*}
Since $\nabla\cdot\bu=0$ in $\mathcal{O}$ and $\bu\cdot\boldsymbol{n}=0$ on $\Gamma$, we immediately have 
\begin{align}\label{eqn-best}
	b(\bu,\bv,\bw)=-b(\bu,\bw,\bv)\ \text{ and }\ b(\bu,\bv,\bv)=0\ \text{ for all } \ \bu\in\mathbb{V}_{\mathrm{div}}, \ \bv,\bw\in\mathbb{V}.  
\end{align}
Using H\"older's, Gagliardo-Nirenberg's and Young's inequalities, we estimate $b(\bu,\bv,\bw)$ as 
\begin{align}
	|b(\bu,\bv,\bw)|&\leq\|\bu\|_{\mathbb{L}^4}\|\nabla\bv\|_0\|\bw\|_{\mathbb{L}^4}\leq C_g^2\|\bu\|_0^{1-\frac{d}{4}}\|\nabla\bu\|_0^{\frac{d}{4}}\|\nabla\bv\|_0\|\bw\|_0^{1-\frac{d}{4}}\|\nabla\bw\|_0^{\frac{d}{4}}\nonumber\\&\leq C_b\|\nabla\bu\|_0\|\nabla\bv\|_0\|\nabla\bw\|_0,
\end{align}
where $C_b = C_g^2 C_p^{\,2-\tfrac{d}{2}}$, where $C_g$ and $C_p$ denote the constants arising in the Gagliardo-Nirenberg and Poincaré inequalities, respectively. Therefore, we can define a bilinear operator $\mathcal{B}(\cdot,\cdot):\mathbb{V}\times \mathbb{V}\to \mathbb{V}^{\prime}$ such that $$\langle\mathcal{B}(\bu,\bv),\bw\rangle=b(\bu,\bv,\bw),$$ for all $\bu,\bv,\bw\in \mathbb{V}$. 

It can also be seen that 
\begin{align*}
	c_r^0(\bu,\bu)=\|\bu\|_{\mathbb{L}^{r+1}}^{r+1}\ \text{ for all }\ \bu\in \mathbb{L}^{r+1}(\mathcal{O}). 
\end{align*}
Moreover, we have 
\begin{align*}
	|c_r^0(\bu,\bv)|=|(|\bu|^{r-1}\bu,\bv)|\leq\|\bu\|_{\mathbb{L}^{r+1}}^{r}\|\bv\|_{\mathbb{L}^{r+1}},
\end{align*}
so that we can define an operator $\mathcal{C}_r(\cdot):\mathbb{L}^{r+1}(\mathcal{O})\to \mathbb{L}^{\frac{r+1}{r}}(\mathcal{O})$ such that 
\begin{align*}
	\mathcal{C}_r(\bu)=|\bu|^{r-1}\bu.
\end{align*}
From \cite[Section 2.4]{SGMTM}, we infer for all $\bu,\bv\in\L^{r+1}(\mathcal{O})$ and $r\geq 1$ that
\begin{align}\label{2.23}
	&\langle|\bu|^{r-1}\bu-|\bv|^{r-1}\bv,\bu-\bv\rangle\geq \frac{1}{2}\||\bu|^{\frac{r-1}{2}}(\bu-\bv)\|_{\H}^2+\frac{1}{2}\||\bv|^{\frac{r-1}{2}}(\bu-\bv)\|_{\H}^2\geq 0,
\end{align}
and 
\begin{align}\label{Eqn-mon-lip}
	&\langle\bu|\bu|^{r-1}-\bv|\bv|^{r-1},\bu-\bv\rangle
	\geq \frac{1}{2^{r-1}}\|\bu-\bv\|_{\L^{r+1}}^{r+1}.
\end{align}
For $\mathcal{C}_r(\bu)=|\bu|^{r-1}\bu$ such that $c_r^0(\bu,\bv)=\langle\mathcal{C}_r(\bu),\bv\rangle$, for every $\bv\in{\L}^{r+1}(\mathcal{O})$,  the Gateaux derivative is obtained as 
\begin{align}\label{Gaetu}
	\mathcal{C}_r'(\bu)\bv&=\left\{\begin{array}{cl}\bv,&\text{ for }r=1,\\ \left\{\begin{array}{cc}|\bu|^{r-1}\bv
	+(r-1)|\bu|^{r-3}(\bu\cdot\bv)\,\bu,&\text{ if }\bu\neq \boldsymbol{0},\\\boldsymbol{0},&\text{ if }\bu=\boldsymbol{0},\end{array}\right.&\text{ for } 1<r<3,\\ |\bu|^{r-1}\bv+(r-1)|\bu|^{r-3}(\bu\cdot\bv)\,\bu, &\text{ for }r\geq 3.\end{array}\right.
\end{align}
Similar properties hold true for the operator $c^0_q(\cdot,\cdot)$ also. 
	We now state the following result, which plays a key role in establishing the convergence of the nonlinear pumping and damping terms. 
\begin{lemma}[Brezis-Lions Lemma, {\cite[Lemma 1.3]{JLL}}]\label{Lem-Lions}
	Let $\mathcal{O} \subset \R^n$ be bounded open, and for $m\in\N$ and $1<q <\infty$, let $\varphi_m$ and $\varphi$ be functions in $\mathrm{L}^q(\mathcal{O})$.  Suppose that the sequence $\{\varphi_m\}_{m \in \N}$ is uniformly bounded in
	$L^q(\mathcal{O})$, that is, for some $C>0$ independent of $m,$
	\begin{align*}
		\|\varphi_m\|_{L^q(\mathcal{O})} \leq C,\ \ \mbox{for every}\ \ m\in\N \; \; \mbox{and} \ \ \varphi_m \to \varphi\;\; \mbox{a.e. in}\ \ \mathcal{O},\  \mbox{as}\ \ m \to \infty.
	\end{align*}
	Then, $\varphi_m \rightharpoonup  \varphi$ in $\mathrm{L}^q(\mathcal{O})$, as $m\to \infty$.
\end{lemma}

\subsection{The problem formulation}
For $\boldsymbol{\eta}\in\mathbb{L}^2(\Gamma_{\mathrm{S}}):=  L^2(\Gamma_{\mathrm{S}};\mathbb{R}^{d})$, we define 
\begin{align}\label{eqn-def-j}
	j(\boldsymbol{\eta})=\int_{\Gamma_{\mathrm{S}}}g|\boldsymbol{\eta}|\d\mathrm{S}.
\end{align}
It is clear that $j$ is a continuous functional defined on $\mathbb{L}^2(\Gamma_{\mathrm{S}})$.  Note that 
\begin{align}
\boldsymbol{\psi}\in	\partial 	j(\boldsymbol{\eta})\iff \begin{cases}
	\boldsymbol{\psi}=g\dfrac{\boldsymbol{\eta}}{|\boldsymbol{\eta}|},&\text{for  }\boldsymbol{\eta}\neq0,\\[6pt]
	|\boldsymbol{\psi}|\le g,&\text{for  }\boldsymbol{\eta}=0.
\end{cases}
\end{align}
The boundary condition in \eqref{eqn-boundary-1} is equivalent to requiring that the negative tangential stress belongs to the subdifferential of the friction functional at the tangential velocity, that is, $-\boldsymbol{\sigma}_{\boldsymbol{\tau}}\in\partial j(\bu_{\boldsymbol{\tau}})$. Moreover, it can be easily seen that 
\begin{align}
|	j(\boldsymbol{\eta_1})-j(\boldsymbol{\eta_2})|&=\int_{\Gamma_{\mathrm{S}}}g\left||\boldsymbol{\eta_1}|-|\boldsymbol{\eta_2}|\right|\d\mathrm{S}\leq \int_{\Gamma_{\mathrm{S}}}g\left[|\boldsymbol{\eta_1}-\boldsymbol{\eta_2}|\right]\d\mathrm{S} \nonumber\\&\leq \|g\|_{L^2(\Gamma_{\mathrm{S}})}\|\boldsymbol{\eta_1}-\boldsymbol{\eta_2}\|_{\mathbb{L}^2(\Gamma_{\mathrm{S}})}. 
\end{align}

Let $\lambda_0$ denote the smallest eigenvalue of the Steklov type eigenvalue problem
\begin{align} 
	\bu\in\mathbb{V}, \ \int_{\mathcal{O}}\nabla\bu:\nabla\bv \d\bx =\lambda\int_{\Gamma_{\mathrm{S}}}\bu_{\tau}\cdot\bv_{\tau} d S,\ \text{ for all }\ \bv\in\mathbb{V}. \end{align}
Since $|\Gamma_{\mathrm{D}}| > 0$,  Poincar\'e's  inequality holds on $\mathbb{V}$. Moreover, the trace operator $\bu \mapsto \bu_{\tau}|_{\Gamma_{\mathrm{S}}}$ is compact from $\mathbb{V}$ into $\mathbb{L}^2(\Gamma_{\mathrm{S}})$ (because the embedding $\mathbb{H}^1(\mathcal{O}) \hookrightarrow \mathbb{L}^2(\Gamma)$ is compact; see \cite[Theorem 2.21]{SMAOMS}). Hence, an application of the spectral theorem for compact self-adjoint operators yields a sequence of eigenvalues $\{\lambda_k\}_{k=0}^\infty$ such that
$
\lambda_k > 0, \ \lambda_k \to +\infty \ \text{as } k \to \infty.
$
In addition, the following trace inequality holds:
\begin{align}\label{eqn-trace} \|\bv_{\tau}\|_{\L^2(\Gamma_{\mathrm{S}})}\leq\lambda_0^{-1/2}\|\nabla\bv\|_{0}, \ \text{ for all }\ \bv\in\mathbb{V}. 
\end{align}



Let us now provide the mixed as well as variational formulations of our problem: 
	\begin{problem}[Mixed (saddle-point) formulation]\label{problem-1}
		Find $(\bu,p,\boldsymbol{\sigma}_{\boldsymbol{\tau}})\in\mathbb{V}\cap\mathbb{L}^{r+1}(\mathcal{O})\times Q\times\mathbb{L}^2(\Gamma_{\mathrm{S}})$ with $-\boldsymbol{\sigma}_{\boldsymbol{\tau}}\in\partial j(\bu_{\boldsymbol{\tau}})$ such that
		\begin{equation}\label{eqn-problem-1}
		\begin{aligned}
			\mu a(\bu,\bv)&+b(\bu,\bu,\bv)+\alpha(\bu,\bv)+\beta c_r^0(\bu,\bv)-\kappa c_q^0(\bu,\bv)+d(\bv,p)-(\boldsymbol{\sigma}_{\boldsymbol{\tau}},\boldsymbol{v}_{\boldsymbol{\tau}})_{\mathbb{L}^2(\Gamma_{\mathrm{S}})}\\&=\langle\f, \bv\rangle, \ \text{ for all }\ \bv\in\mathbb{V}\cap\mathbb{L}^{r+1}(\mathcal{O}),\\
			b(\bu,q)&=0, \ \text{ for all }\ q\in Q. 
		\end{aligned}
			\end{equation}
	\end{problem}
	
	\begin{problem}[Variational inequality with pressure]\label{problem-2}
		Find $(\bu,p)\in\mathbb{V}\cap\mathbb{L}^{r+1}(\mathcal{O})\times Q$ such that for all $\bv\in\mathbb{V}\cap\mathbb{L}^{r+1}(\mathcal{O})$ and all $q\in Q$,
		\begin{equation}\label{eqn-problem-2}
		\begin{aligned}
			\mu a(\bu,\bv-\bu)+b(\bu,\bu,\bv-\bu)+\alpha(\bu,\bv-\bu)&+\beta c_r^0(\bu,\bv-\bu)-\kappa c_q^0(\bu,\bv-\bu)\\+d(\bv-\bu,p)+j(\boldsymbol{v}_{\boldsymbol{\tau}})-j(\bu_{\boldsymbol{\tau}})&\geq \langle\f, \bv-\bu\rangle,\\
			d(\bu,q)&=0. 
		\end{aligned}
		\end{equation}
	\end{problem}
	
	Upon eliminating the pressure variable $p$, one arrives at the following reduced hemivariational inequality:
	
		\begin{problem}[Variational inequality without pressure]\label{problem-3}
		Find $\bu\in\mathbb{V}_{\mathrm{div}}\cap\mathbb{L}^{r+1}(\mathcal{O})$ such that 
		\begin{equation}\label{eqn-problem-3}
			\begin{aligned}
				\mu a(\bu,\bv-\bu)+b(\bu,\bu,\bv-\bu) & +\alpha(\bu,\bv-\bu)+\beta c_r^0(\bu,\bv-\bu)-\kappa c_q^0(\bu,\bv-\bu)\\+j(\boldsymbol{v}_{\boldsymbol{\tau}})-j(\bu_{\boldsymbol{\tau}})&\geq \langle\f, \bv-\bu\rangle, \ \mbox{ for all \ $\bv\in\mathbb{V}_{\mathrm{div}}\cap\mathbb{L}^{r+1}(\mathcal{O})$.}
			\end{aligned}
		\end{equation}
	\end{problem}
	\begin{lemma}[Equivalence between Problems \ref{problem-1} and \ref{problem-2}]
		Problems \ref{problem-1} and \ref{problem-2} are equivalent. Moreover, the unknowns $(y,p)$ in both formulations coincide.
	\end{lemma}
	\begin{proof}
		Let $(\bu,p,\boldsymbol{\sigma}_{\boldsymbol{\tau}})$ satisfy \eqref{eqn-problem-1}. Taking the difference of the relations for $\bv$ and $\bu$ gives 
			\begin{align}\label{eqn-est-1}
			\mu a(\bu,\bv-\bu)&+b(\bu,\bu,\bv-\bu)+\alpha(\bu,\bv-\bu)+\beta c_r^0(\bu,\bv-\bu)-\kappa c_q^0(\bu,\bv-\bu)\nonumber\\&+d(\bv-\bu,p)-(\boldsymbol{\sigma}_{\boldsymbol{\tau}},\boldsymbol{v}_{\boldsymbol{\tau}}-\bu_{\boldsymbol{\tau}})_{\mathbb{L}^2(\Gamma_{\mathrm{S}})}=\langle\f, \bv-\bu\rangle, \ \text{ for all }\ \bv\in\mathbb{V}\cap\mathbb{L}^{r+1}(\mathcal{O}). 
		\end{align}
		Since $-\boldsymbol{\sigma}_{\boldsymbol{\tau}}\in\partial j(\bu_{\boldsymbol{\tau}})$, we immediately have 
		\begin{align}
			(-\boldsymbol{\sigma}_{\boldsymbol{\tau}},\boldsymbol{v}_{\boldsymbol{\tau}}-\bu_{\boldsymbol{\tau}})_{\mathbb{L}^2(\Gamma_{\mathrm{S}})}\leq j(\bv_{\boldsymbol{\tau}})-j(\bu_{\boldsymbol{\tau}}). 
		\end{align}
		Inserting this into the identity \eqref{eqn-est-1} yields \eqref{eqn-problem-2}.

	Assume $(\bu,p)$ satisfies Problem \ref{problem-2}.
	Define the linear functional $\mathcal{L}$ on $\mathbb{V}\cap\mathbb{L}^{r+1}(\mathcal{O})$ by
	$$
	\mathcal{L}(\bv):=\mu a(\bu,\bv)+b(\bu,\bu,\bv)+\alpha(\bu,\bv)+\beta c_r^0(\bu,\bv)-\kappa c_q^0(\bu,\bv)+d(\bv,p)-\langle\f,\bv\rangle.
	$$
	By continuity of the forms and the data,  $\mathcal{L}$ is well-defined and continuous on $\mathbb{V}\cap\mathbb{L}^{r+1}(\mathcal{O})$.
	First, we show that $\mathcal{L}$ vanishes on the subspace of test functions with zero tangential trace. Let $\bw\in\mathbb{V}\cap\mathbb{L}^{r+1}(\mathcal{O})$ satisfy $\bw_{\boldsymbol\tau}=0$. Take $\bv=\bu\pm\bw$ in the variational inequality. Because $(\bu\pm\bw)_{\boldsymbol\tau}=\bu_{\boldsymbol\tau}$, the $j$-terms cancel and we obtain for both signs
	$$
	\mathcal{L}(\pm\bw)\ge 0,
	$$
	so that $\mathcal{L}(\bw)=0$. Thus $\mathcal{L}$ depends only on the tangential trace of its argument; that is, there exists a linear functional $\ell$ on the trace space $T:=\{\bv_{\boldsymbol\tau |_{\Gamma_S}}:\bv\in\mathbb V\cap\mathbb{L}^{r+1}(\mathcal{O})\}$ such that
$$
	\mathcal{L}(\bv)=\ell(\bv_{\boldsymbol\tau})\ \text{ for all }\ \bv\in\mathbb{V}\cap\mathbb{L}^{r+1}(\mathcal{O}).
	$$
	Next, rewrite the variational inequality in terms of $\ell$. For arbitrary $\bv\in\mathbb V\cap\mathbb{L}^{r+1}(\mathcal{O})$,
	$$
	\mathcal{L}(\bv-\bu)+j(\bv_{\boldsymbol\tau})-j(\bu_{\boldsymbol\tau})\geq 0
	\Longrightarrow
	j(\bv_{\boldsymbol\tau})-j(\bu_{\boldsymbol\tau})+\ell(\bv_{\boldsymbol\tau}-\bu_{\boldsymbol\tau})\ge 0.
	$$
	Therefore, $\bu_{\boldsymbol\tau}$ is a minimizer of the convex functional $\boldsymbol{\xi}\mapsto j(\boldsymbol{\xi})  +\ell (\boldsymbol{\xi})$ on the trace space $T$. By convex analysis, there exists a subgradient $\boldsymbol{\eta}\in T^*$ with
		$$
	-\boldsymbol{\eta}\in\partial j(\bu_{\boldsymbol\tau})\ \text{ and }\  \boldsymbol{\eta}=\ell\ \text{ in }T^*.
	$$
	Finally, identify $\boldsymbol{\eta}$ with a function $\boldsymbol\sigma_{\boldsymbol\tau}\in\mathbb L^2(\Gamma_{\mathrm S})$. 
	This identification follows from the fact that the trace space $T$ is continuously embedded in $\mathbb{L}^2(\Gamma_{\mathrm S})$. Under the assumptions of the problem, the dual space $T^*$ can therefore be represented by functions in $\mathbb{L}^2(\Gamma_{\mathrm S})$. In particular, if the trace operator maps $\mathbb V$ onto a Hilbert subspace of $\mathbb{L}^2(\Gamma_{\mathrm S})$, then by the Riesz representation theorem, there exists $\boldsymbol\sigma_{\boldsymbol\tau}\in\mathbb{L}^2(\Gamma_{\mathrm S})$ such that
	$$
	\boldsymbol{\eta}(\boldsymbol{\zeta})=(\boldsymbol\sigma_{\boldsymbol\tau},\boldsymbol{\zeta})_{\mathbb L^2(\Gamma_{\mathrm S})}
	\ \text{ for all } \ \boldsymbol{\zeta}\in T.
	$$
	With this identification, we have
		$$
	-\boldsymbol\sigma_{\boldsymbol\tau}\in\partial j(\bu_{\boldsymbol\tau})\ \text{ and }\
	\ell(\bv_{\boldsymbol\tau})=(\boldsymbol\sigma_{\boldsymbol\tau},\bv_{\boldsymbol\tau})_{\mathbb L^2(\Gamma_{\mathrm S})}\ \text{ for all }\ \bv\in\mathbb{V}.
	$$
		But $\ell(\bv_{\boldsymbol\tau})=\mathcal{L}(\bv)$ by definition, so for all $\bv\in\mathbb V\cap\mathbb{L}^{r+1}(\mathcal{O})$,
		$$
	\mu a(\bu,\bv)+b(\bu,\bu,\bv)+\alpha(\bu,\bv)+\beta c_r^0(\bu,\bv)-\kappa c_q^0(\bu,\bv)+d(\bv,p)-(\boldsymbol\sigma_{\boldsymbol\tau},\bv_{\boldsymbol\tau})_{\mathbb L^2(\Gamma_{\mathrm S})}
	=\langle\f,\bv\rangle.
	$$
		Together with $ -\boldsymbol\sigma_{\boldsymbol\tau}\in\partial j(\bu_{\boldsymbol\tau})$ and $d(\bu,q)=0$ this exactly yields Problem \ref{problem-1}.
		\end{proof}
		
			\section{Existence and uniqueness of weak solution}\label{sec3}\setcounter{equation}{0} 
		We now address the existence and uniqueness of solutions to Problem \ref{problem-3}. To this end, we introduce a regularized formulation of \eqref{eqn-problem-3}. Existence of a solution to the regularized problem is first established via a Galerkin approximation combined with compactness arguments. We then prove the existence and uniqueness of weak solutions to Problem \ref{problem-3} for $d\in\{2,3\}$ and $r\in[1,\infty)$. Finally, by exploiting the inf–sup condition (also known as the Ladyzhenskaya-Babuška-Brezzi (LBB) condition), we establish the existence and uniqueness of solutions to Problem \ref{problem-2} for $1 \leq r \leq \tfrac{2d}{d-2}$ when $d>2$, and for all $1 \leq r < \infty$ when $d=2$.

		\subsection{Regularized problem} \label{subsec:regprob}
		Since the functional $j(\boldsymbol{\eta})$ is not differentiable with respect to $\boldsymbol{\eta}$, it is natural to approximate it by a family of convex and differentiable functionals $j_{\varepsilon}(\boldsymbol{\eta})$, where $\varepsilon > 0$ is a small parameter. For this purpose, we employ the Yosida regularization method. The regularized functional is defined as
			$$
		j_{\varepsilon}(\boldsymbol{\eta})=\int_{\Gamma_{\mathrm{S}}} g\varrho_{\varepsilon}(\boldsymbol{\eta}) \d\mathrm{S},
		$$
		where
			$$
		\varrho_{\varepsilon}(\boldsymbol{\zeta}) = 
		\begin{cases}
			|\boldsymbol{\zeta}| - \tfrac{\varepsilon}{2}, & \text{if } |\boldsymbol{\zeta}| > \varepsilon, \\[6pt]
			\dfrac{|\boldsymbol{\zeta}|^2}{2\varepsilon}, & \text{if } |\boldsymbol{\zeta}| \leq \varepsilon,
		\end{cases}
		$$
			which corresponds to the Yosida regularization of $\varrho(\boldsymbol{\zeta}) = |\boldsymbol{\zeta}|$. The properties of $j_{\varepsilon}(\cdot)$ are established in \cite{NSa}.

		\begin{lemma}\label{lem-diff-Yosida}
			The Yosida regularization operator $j_{\varepsilon}$ is convex and Gâteaux differentiable in $\mathbb{L}^2(\Gamma_{\mathrm{S}})$. Moreover, for all $\boldsymbol{\eta}, \boldsymbol{\xi} \in \mathbb{L}^2(\Gamma_{\mathrm{S}})$, we have
		\begin{align}
			 \langle\nabla j_{\varepsilon}(\boldsymbol{\eta}),\boldsymbol{\xi}\rangle=\lim_{h\to 0}\frac{1}{h}\left[ j_{\varepsilon}(\boldsymbol{\eta}+h\boldsymbol{\xi})- j_{\varepsilon}(\boldsymbol{\eta})\right]=\int_{\Gamma_{\mathrm{S}}}g\rho_{\varepsilon}(\boldsymbol{\eta})\cdot\boldsymbol{\xi}\d\mathrm{S},  
			 \end{align} where 
			 \begin{align} \rho_{\varepsilon}(\boldsymbol{\zeta})=\partial \varrho_{\varepsilon}(\boldsymbol{\zeta})=\left\{\begin{array}{ll}\boldsymbol{\zeta}/|\boldsymbol{\zeta}|,&\ \text{ if } \ |\boldsymbol{\zeta}|>\varepsilon,\\ \boldsymbol{\zeta}/\varepsilon,&\ \text{ if }\ |\boldsymbol{\zeta}|\leq\varepsilon. \end{array}\right. \end{align}
			In addition, the following estimate holds:
		\begin{align}\label{eqn-j-est} | j_{\varepsilon}(\boldsymbol{\eta})- j(\boldsymbol{\eta})|\leq\frac{\varepsilon}{2}\|g\|_{L^1(\Gamma_{\mathrm{S}})},\ \text{ for all }\ \boldsymbol{\eta}\in \mathbb{L}^2(\Gamma_{\mathrm{S}}). \end{align}
				\end{lemma}
		Having introduced $j_{\varepsilon}(\boldsymbol{\eta})$, we now formulate the corresponding regularized problem associated with the variational inequality \eqref{eqn-problem-3}.

			\begin{problem}\label{problem-4}
				Find $\bu_{\varepsilon}\in\mathbb{V}_{\mathrm{div}}\cap\mathbb{L}^{r+1}(\mathcal{O})$ such that 
				\begin{equation}\label{eqn-problem-4}
					\begin{aligned}
						\mu a(\bu_{\varepsilon},\bv-\bu_{\varepsilon}) & +b(\bu_{\varepsilon},\bu_{\varepsilon},\bv-\bu_{\varepsilon})+\alpha(\bu_{\varepsilon},\bv-\bu_{\varepsilon})+\beta c_r^0(\bu_{\varepsilon},\bv-\bu_{\varepsilon})-\kappa c_q^0(\bu_{\varepsilon},\bv-\bu_{\varepsilon})\\ &+j_{\varepsilon}(\boldsymbol{v}_{\boldsymbol{\tau}})-j_{\varepsilon}(\bu_{\varepsilon\boldsymbol{\tau}})\geq \langle\f, \bv-\bu_{\varepsilon}\rangle, \ \mbox{ for all \ $\bv\in \mathbb{V}_{\mathrm{div}}\cap\mathbb{L}^{r+1}(\mathcal{O})$.}
					\end{aligned}
				\end{equation}
			\end{problem}
		Since $j_{\varepsilon}$ is differentiable, Lemma \ref{lem-diff-Yosida} together with standard results from convex analysis implies that the variational inequality problem \eqref{eqn-problem-4} is equivalent to the following variational problem:

				\begin{problem}\label{problem-5}
				Find $\bu_{\varepsilon}\in\mathbb{V}_{\mathrm{div}}\cap\mathbb{L}^{r+1}(\mathcal{O})$ such that 
				\begin{equation}\label{eqn-problem-5}
					\begin{aligned}
						\mu a(\bu_{\varepsilon},\bw)+b(\bu_{\varepsilon},\bu_{\varepsilon},\bw)+\alpha(\bu_{\varepsilon},\bw)&+\beta c_r^0(\bu_{\varepsilon},\bw)-\kappa c_q^0(\bu_{\varepsilon},\bw)\\+(g\rho_{\varepsilon}(\bu_{\varepsilon\boldsymbol{\tau}}),\bw_{\boldsymbol{\tau}})_{\mathbb{L}^2(\Gamma_{\mathrm{S}})}&=\langle\f, \bw\rangle, \ \mbox{ for all \ $\bw\in \mathbb{V}_{\mathrm{div}} \cap\mathbb{L}^{r+1}(\mathcal{O})$.}
					\end{aligned}
				\end{equation}
			\end{problem}
			Using Lemma \ref{lem-diff-Yosida} and the convexity of $j_{\varepsilon}$, we obtain
		\begin{align} \int_{\Gamma_{\mathrm{S}}}g\rho_{\varepsilon}(\bu_{\varepsilon\boldsymbol{\tau}})\cdot(\boldsymbol{v}_{\boldsymbol{\tau}}-\bu_{\varepsilon\boldsymbol{\tau}})\d\mathrm{S}\leq j_{\varepsilon}(\boldsymbol{v}_{\boldsymbol{\tau}})-j_{\varepsilon}(\bu_{\varepsilon\boldsymbol{\tau}}), \ \text{ for all }\ \bv\in\mathbb{V}_{\mathrm{div}}. 
		\end{align}
			It follows that any solution $\bu_{\varepsilon} \in \mathbb{V}_{\mathrm{div}}\cap\mathbb{L}^{r+1}(\mathcal{O})$ of \eqref{eqn-problem-5} is also a weak solution of the variational inequality \eqref{eqn-problem-4}.
			Conversely, by choosing test functions of the form $\bv = \bu_{\varepsilon} \pm \lambda \bw$ with $\bw \in \mathbb{V}_{\mathrm{div}}$ and arbitrary $\lambda > 0$ in \eqref{eqn-problem-4}, and then letting $\lambda \to 0$, we recover the variational problem \eqref{eqn-problem-5}.

\subsection{Existence of a weak solution to regularized problem} \label{subsec:existweak}
We begin by proving the existence of a weak solution to problem \eqref{eqn-problem-5} via a Galerkin approximation, and subsequently pass to the limit $\varepsilon\to 0$ to derive the existence of a weak solution to problem \eqref{eqn-problem-4}.
			
\begin{proposition}\label{prop-existence-ueps}	Given $\f\in\mathbb{V}^{\prime}$ and $g\in L^2(\Gamma_{\mathrm{S}})$,  there exists a weak solution $\bu_{\varepsilon}$ of the variational problem \eqref{eqn-problem-5} and the following estimate holds:
\begin{align}	&\frac{\mu}{2}\|\nabla\bu_{\varepsilon}\|_0^2+\alpha\|\bu_{\varepsilon}\|_0^2+\frac{\beta}{2}\|\bu_{\varepsilon}\|_{\mathbb{L}^{r+1}}^{r+1}\le \mathcal{K}, \label{eqn-ener-1} \\
&\mathcal{K}:=\left(\frac{r-q}{r+1}\right)\left(\frac{2(q+1)}{\beta(r+1)}\right)^{\frac{q+1}{r-q}}\kappa^{\frac{r+1}{r-q}}|\mathcal{O}|+\frac{1}{\mu\lambda_0}\|g\|_{L^2(\Gamma_{\mathrm{S}})}^2+\frac{1}{\mu}\|\f\|_{\mathbb{V}^{\prime}}^2.\label{calK}
\end{align}
\end{proposition}

\begin{proof}
We aim to prove the existence of solutions to the variational problem \eqref{eqn-problem-5} in the following steps:
\vskip 0.1cm
\noindent\textbf{Step 1:} \emph{Galerkin approximation.} 
Let $\{\boldsymbol{e}_k\}_{k=1}^{\infty}\subset \mathbb{H}^2(\mathcal{O})$  be an orthogonal basis in $\mathbb{V}_{\mathrm{div}}$ which is an orthonormal basis in $\mathbb{H}_{\mathrm{div}}$ \cite[Theorem 2.24]{JCBob3D}, \cite[Ch. 1, Sect. 2]{Te}. Since $\mathbb{H}^2(\mathcal{O})\hookrightarrow \mathbb{L}^p(\mathcal{O})$ for all $1\leq p<\infty$ and $1\leq d\leq 4$, we immediately have $\{\boldsymbol{e}_k\}_{k=1}^{\infty}\in \mathbb{L}^{r+1}(\mathcal{O})$ for all $1\leq r<\infty$. Fix an arbitrary $m\in\mathbb{N}$ and define the operator $\mathrm{P}_m$ from $\mathbb{V}_{\mathrm{div}}^{\prime}$ to $\mathrm{span}\{\boldsymbol{e}_1,\ldots,\boldsymbol{e}_m\}=:\mathbb{H}_{\mathrm{div},m} $ by 
\begin{align}\label{eqn-projection}
\mathrm{P}_m\bu^*:=\sum_{k=1}^{m}\langle\bu^*,\boldsymbol{e}_k\rangle\boldsymbol{e}_k,\quad \bu^*\in\mathbb{V}_{\mathrm{div}}^{\prime}.
\end{align}
We denote by $\mathrm{P}_m$ the restriction of the previously defined operator to $\mathbb{H}_{\mathrm{div}}$, without introducing a new notation. In particular, we have $\H\hookrightarrow\mathbb{V}_{\mathrm{div}}^{\prime}$, that is, every element $\bu\in\H_{\mathrm{div}}$ induces a functional $\bu^*\in\mathbb{V}_{\mathrm{div}}^{\prime}$ by 
$\langle\bu^*,\bv\rangle:=(\bu,\bv), \  \bv\in\mathbb{V}_{\mathrm{div}}.$
Thus the restriction of $\P_m$ to $\H_{\mathrm{div}}$ is given by 
\begin{align}
\P_m\bu:=\sum_{k=1}^{m}(\bu,\boldsymbol{e}_k)\boldsymbol{e}_k, \quad  \bu\in\H_{\mathrm{div}}.
\end{align}
Hence, in particular, $\P_m$ is the orthogonal projection from $\H_{\mathrm{div}}$ onto $\H_{\mathrm{div},m}$.  For each $m\in\N$, we search for approximate solutions of the form 
\begin{align}\label{eqn-finite-dim}
\bu_{\varepsilon}^m(x)=\sum_{k=1}^m \upalpha^m_k\boldsymbol{e}_k(x)=\sum_{k=1}^m (\bu_{\varepsilon}^m,\boldsymbol{e}_k)\boldsymbol{e}_k(x), \ \boldsymbol{e}_k\in\mathbb{H}_{\mathrm{div},m}, 
\end{align} 
where the coefficients $\upalpha^m_1,\ldots,\upalpha^m_m$ are solutions of the following $m$ equations:
\begin{equation}\label{eqn-Galerkin}
\begin{aligned}
\mu a(\bu_{\varepsilon}^m,\boldsymbol{e}_k)+b(\bu_{\varepsilon}^m,\bu_{\varepsilon}^m,\boldsymbol{e}_k)+\alpha(\bu_{\varepsilon}^m,\boldsymbol{e}_k)&+\beta c_r^0(\bu_{\varepsilon}^m,\boldsymbol{e}_k)-\kappa c_q^0(\bu_{\varepsilon}^m,\boldsymbol{e}_k)\\+(g\rho_{\varepsilon}(\bu_{\varepsilon\boldsymbol{\tau}}^m),\boldsymbol{e}_{k\boldsymbol{\tau}})_{\mathbb{L}^2(\Gamma_{\mathrm{S}})}&=\langle\f, \boldsymbol{e}_k\rangle, \quad  1\le k \le m.
\end{aligned}
\end{equation}
We define $\mathcal{B}^m=\mathrm{P}_m\mathcal{B}$, $\mathcal{C}^m_r=\mathrm{P}_m\mathcal{C}_r$, $\mathcal{C}^m_q=\mathrm{P}_m\mathcal{C}_q$, so that
$$\langle\mathcal{B}^m(\bu_{\varepsilon}^m),\bw\rangle=\langle\mathcal{B}(\bu_{\varepsilon}^m,\bu_{\varepsilon}^m),\bw\rangle=b(\bu_{\varepsilon}^m,\bu_{\varepsilon}^m,\bw)\ \text{ for all }\  \bw\in\H_{\mathrm{div},m}$$ and $$\langle\mathcal{C}^m_r(\bu_{\varepsilon}^m),\bw\rangle=\langle\mathcal{C}_r(\bu_{\varepsilon}^m),\bw\rangle=c_r^0(\bu_{\varepsilon}^m,\bw), \ \text{ for all }\ \bw\in\H_{\mathrm{div},m}.$$
\vskip 0.1cm
\noindent\textbf{Step 2:} \emph{Existence of solutions to finite dimensional problem.} 
Following the approach of \cite[Theorem 1.2, p. 110]{Te} and \cite[Theorem 5.1]{KKMTM-2025}, we establish existence results for the finite-dimensional system \eqref{eqn-Galerkin}. The equations \eqref{eqn-finite-dim}-\eqref{eqn-Galerkin} constitute a nonlinear system for the coefficients $\upalpha^m_1,\ldots,\upalpha^m_m$. By applying \cite[Lemma 1.4, Chapter 2]{Te}, we deduce the existence of a solution to this system.
We set $W=\mathrm{span}\{\boldsymbol{e}_1,\ldots,\boldsymbol{e}_m\}$, endowed with the scalar product
$(\cdot,\cdot)_W := (\nabla\cdot,\nabla\cdot),$
which is induced by $\mathbb{V}_{\mathrm{div}}$. The projection operator $\mathrm{P}=\mathrm{P}_m$ is then defined by
\begin{align}
(\mathrm{P}_m(\bu_{\varepsilon}),\bv_{\varepsilon})_W& = (\nabla\mathrm{P}_m(\bu_{\varepsilon}),\nabla\bv_{\varepsilon})
\nonumber\\ & := \mu (\nabla\bu_{\varepsilon},\nabla\bv_{\varepsilon})+ b(\bu_{\varepsilon},\bu_{\varepsilon},\bv_{\varepsilon}) + \alpha (\bu_{\varepsilon}, \bv_{\varepsilon})+\beta c_r^0(\bu_{\varepsilon},\bv_{\varepsilon})-\kappa c_q^0(\bu_{\varepsilon},\bv_{\varepsilon})\nonumber\\&\quad+(g\rho_{\varepsilon}(\bu_{\varepsilon\boldsymbol{\tau}}),\bv_{\varepsilon\boldsymbol{\tau}})-\left\langle\f,\bv_{\varepsilon}\right\rangle,
\end{align}
for all $\bu_{\varepsilon},\bv_{\varepsilon}\in W$. Since $\mathrm{P}_m$ acts on a finite-dimensional space, its continuity as a mapping from $W$ into itself follows immediately from the equivalence of norms in finite-dimensional spaces. Hence, to apply Lemma 1.4 of Chapter 2 in \cite{Te}, it suffices to verify that $$(\mathrm{P}_m(\bu_{\varepsilon}),\bu_{\varepsilon})_W>0, \ \text{ for } \ \|\bu_{\varepsilon}\|_W=k>0.$$ Here $\|\cdot\|_W$ denotes the norm on $W$. We now employ the Cauchy-Schwarz, H\"older's  and Young's inequalities along with the trace estimate \eqref{eqn-trace}, to derive 
\begin{align}
(\mathrm{P}_m(\bu_{\varepsilon}),\bu_{\varepsilon})_W&=\mu \|\nabla\bu_{\varepsilon}\|_{0}^2+ \alpha \|\bu_{\varepsilon}\|_{0}^2 +\beta\|\bu_{\varepsilon}\|_{\mathbb{L}^{r+1}}^{r+1}-\kappa\|\bu_{\varepsilon}\|_{\mathbb{L}^{q+1}}^{q+1}\nonumber\\&\quad+(g\rho_{\varepsilon}(\bu_{\varepsilon\boldsymbol{\tau}}),\bu_{\varepsilon\boldsymbol{\tau}})-\left\langle\f,\bu_{\varepsilon}\right\rangle
\nonumber\\ & \geq \mu \|\nabla\bu_{\varepsilon}\|_{0}^2+ \alpha \|\bu_{\varepsilon}\|_{0}^2 +\beta\|\bu_{\varepsilon}\|_{\mathbb{L}^{r+1}}^{r+1}-\kappa|\mathcal{O}|^{\frac{r-q}{r+1}}\|\bu_{\varepsilon}\|_{\mathbb{L}^{r+1}}^{q+1}\nonumber\\&\quad-\|g\|_{L^2(\Gamma_{\mathrm{S}})}\|\bu_{\varepsilon\boldsymbol{\tau}}\|_{\mathbb{L}^2(\Gamma_{\mathrm{S}})}-\|\f\|_{\mathbb{V}^{\prime}}\|\bu_{\varepsilon}\|_{\mathbb{V}}
\nonumber\\&\geq \frac{\mu}{2} \|\nabla\bu_{\varepsilon}\|_{0}^2 +\frac{\beta}{2}\|\bu_{\varepsilon}\|_{\mathbb{L}^{r+1}}^{r+1}-\left(\frac{r-q}{r+1}\right)\left(\frac{2(q+1)}{\beta(r+1)}\right)^{\frac{q+1}{r-q}}\kappa^{\frac{r+1}{r-q}}|\mathcal{O}|\nonumber\\&\quad-\frac{1}{\mu\lambda_0}\|g\|_{L^2(\Gamma_{\mathrm{S}})}^2-\frac{1}{\mu}\|\f\|_{\mathbb{V}^{\prime}}^2,
\end{align}
where we have used the fact that $|\rho_{\varepsilon}(\bu_{\varepsilon\boldsymbol{\tau}})|\leq 1$. 	It follows that $$(\mathrm{P}_m(\bu_{\varepsilon}),\bu_{\varepsilon})_W>0\ \text{ for }\ \|\bu_{\varepsilon}\|_{W}=k$$ and $k$ is sufficiently large, more precisely $$k>\sqrt{\frac{2}{\mu}}\left(\left(\frac{r-q}{r+1}\right)\left(\frac{2(q+1)}{\beta(r+1)}\right)^{\frac{q+1}{r-q}}\kappa^{\frac{r+1}{r-q}}|\mathcal{O}|+\frac{1}{\mu\lambda_0}\|g\|_{L^2(\Gamma_{\mathrm{S}})}^2+\frac{1}{\mu}\|\f\|_{\mathbb{V}^{\prime}}^2\right)^{1/2}.$$ Hence, the hypotheses of \cite[Lemma 1.4, Chapter 2]{Te} are satisfied and a solution $\bu_{\varepsilon}^m$  of \eqref{eqn-Galerkin} exists. 
								
\vskip 0.1cm
\noindent\textbf{Step 3:} \emph{Uniform energy estimates and convergences.} 
Multiplying \eqref{eqn-Galerkin} by $\upalpha_{k,\varepsilon}^m$ and then summing over $k=1,\ldots,m$, we find 
\begin{align}
&	\mu\|\nabla\bu_{\varepsilon}^m\|_0^2+\alpha\|\bu_{\varepsilon}^m\|_0^2+\beta\|\bu_{\varepsilon}^m\|_{\mathbb{L}^{r+1}}^{r+1}\nonumber\\&=\kappa\|\bu_{\varepsilon}^m\|_{\mathbb{L}^{q+1}}^{q+1}-(g\rho_{\varepsilon}^m(\bu_{\varepsilon\boldsymbol{\tau}}^m),\bu_{k\boldsymbol{\tau}}^m)_{\mathbb{L}^2(\Gamma_{\mathrm{S}})}+(\f,\bu_{\varepsilon}^m)
\nonumber\\&\leq\frac{\mu}{2} \|\nabla\bu^m_{\varepsilon}\|_{0}^2+\frac{\beta}{2}\|\bu_{\varepsilon}^m\|_{\mathbb{L}^{r+1}}^{r+1}+\left(\frac{r-q}{r+1}\right)\left(\frac{2(q+1)}{\beta(r+1)}\right)^{\frac{q+1}{r-q}}\kappa^{\frac{r+1}{r-q}}|\mathcal{O}|\nonumber\\&\quad+\frac{1}{\mu\lambda_0}\|g\|_{L^2(\Gamma_{\mathrm{S}})}^2+\frac{1}{\mu}\|\f\|_{\mathbb{V}^{\prime}}^2.
\end{align}
We therefore obtain the following uniform energy estimate:
\begin{equation}\label{eqn-ener-est}
\frac{\mu}{2}\|\nabla\bu_{\varepsilon}^m\|_0^2+\alpha\|\bu_{\varepsilon}^m\|_0^2+\frac{\beta}{2}\|\bu_{\varepsilon}^m\|_{\mathbb{L}^{r+1}}^{r+1}\le  \mathcal{K}.
\end{equation}
Thus, the sequence $\{\bu_{\varepsilon}^m\}_{m\in\mathbb{N}}$ is uniformly bounded in $\mathbb{V}_{\mathrm{div}}\cap\mathbb{L}^{r+1}(\mathcal{O})$. By the Banach-Alaoglu theorem, there exist a subsequence (still denoted by $\{\bu_{\varepsilon}^m\}_{m\in\mathbb{N}}$) and an element $\bu_{\varepsilon}\in\mathbb{V}_{\mathrm{div}}\cap \L^{r+1}(\mathcal{O})$  such that
\begin{equation}\label{eqn-conv-1}
\left\{
\begin{aligned}
\bu_{\varepsilon}^m &\xrightarrow{w} \bu_{\varepsilon}\ \text{ in }\ \mathbb{V}_{\mathrm{div}},\\
\bu_{\varepsilon}^m &\xrightarrow{w} \bu_{\varepsilon}\ \text{ in }\  \L^{r+1}(\mathcal{O}). 
\end{aligned}
\right. 
\end{equation}
Since the embedding $\mathbb{V}_{\mathrm{div}}\hookrightarrow \mathbb{H}_{\mathrm{div}}$ is compact, we also obtain the strong convergence
\begin{align}\label{eqn-conv-2}
\bu_{\varepsilon}^m \to \bu_{\varepsilon}\ \text{ in }\ \mathbb{H}_{\mathrm{div}}.
\end{align}
Furthermore, along a subsequence (again denoted by the same symbol), we have
\begin{align}\label{eqn-conv-3}
\bu_{\varepsilon}^m(\bx)\to \bu_{\varepsilon}(\bx)\ \text{ for a.e. }\ \bx\in\mathcal{O}.
\end{align}
Since for every $\bv\in\mathbb{H}^1(\mathcal{O})$, we have (\cite[Theorem]{SCBLR})
\begin{align*}
\|\bv\|_{\mathbb{L}^2(\Gamma)}\leq C\|\bv\|_{0}^{1/2}\|\bv\|_{\mathbb{H}^1}^{1/2},\ \text{ for all }\ \bv\in \mathbb{H}^1(\mathcal{O}). 
\end{align*}
Therefore, one can conclude that 
\begin{align}\label{eqn-conv-4}
\bu_{\varepsilon}^m \to \bu_{\varepsilon}\ \text{ in }\ \mathbb{L}^2(\Gamma). 
\end{align}
The above convergence also implies that 
\begin{align}\label{eqn-conv-5}
\bu_{\varepsilon}^m(\bx) \to \bu_{\varepsilon}(\bx)\ \text{ for a.e. }\ \bx\in\Gamma. 
\end{align}
		
\vskip 0.1cm
\noindent\textbf{Step 4:} \emph{Passage to limit.} 
Our next aim is to show that the limit $\bu_{\varepsilon}$ satisfies \eqref{eqn-problem-5}.  Using the convergences given in \eqref{eqn-conv-1}, we have as $m\to\infty,$ 
\begin{align}
|a(\bu_{\varepsilon}^m,\boldsymbol{e}_k)-a(\bu_{\varepsilon},\boldsymbol{e}_k)|&=|(\nabla(\bu_{\varepsilon}^m-\bu_{\varepsilon}),\nabla\boldsymbol{e}_k)|\to 0,\ \text{ as }\ m\to\infty.
\end{align}
The strong convergence in \eqref{eqn-conv-2} implies $(\bu_{\varepsilon}^m,\boldsymbol{e}_k)\to (\bu_{\varepsilon},\boldsymbol{e}_k)$, as $m\to\infty$. Using \eqref{eqn-best}, H\"older's  and Gagliardo-Nirenberg's inequalities, the bound \eqref{eqn-ener-est}, and the convergence \eqref{eqn-conv-2}, we find 
\begin{align}\label{eqn-b-conv}
&|b(\bu_{\varepsilon}^m,\bu_{\varepsilon}^m,\boldsymbol{e}_k)-b(\bu_{\varepsilon},\bu_{\varepsilon},\boldsymbol{e}_k)|
\nonumber\\&\leq|b(\bu_{\varepsilon}^m-\bu_{\varepsilon},\bu_{\varepsilon}^m,\boldsymbol{e}_k)|+|b(\bu_{\varepsilon},\bu_{\varepsilon}^m-\bu_{\varepsilon},\boldsymbol{e}_k)|\nonumber\\&= |b(\bu_{\varepsilon}^m-\bu_{\varepsilon},\boldsymbol{e}_k,\bu_{\varepsilon}^m)|+|b(\bu_{\varepsilon},\boldsymbol{e}_k,\bu_{\varepsilon}^m-\bu_{\varepsilon})|
\nonumber\\&\leq\left(\|\bu_{\varepsilon}^m\|_{\mathbb{L}^4}+\|\bu_{\varepsilon}\|_{\mathbb{L}^4}\right)\|\bu_{\varepsilon}^m-\bu_{\varepsilon}\|_{\mathbb{L}^4}\|\nabla\boldsymbol{e}_k\|_0\nonumber\\&\leq C_g^2C_p^{1-\frac{d}{4}}\left(\|\nabla\bu_{\varepsilon}^m\|_{0}+\|\nabla\bu_{\varepsilon}\|_{0}\right)^{1+\frac{d}{4}} \|\bu_{\varepsilon}^m-\bu_{\varepsilon}\|_{0}^{1-\frac{d}{4}} \|\nabla\boldsymbol{e}_k\|_0\nonumber\\&\to 0\ \text{ as }\ m\to\infty.
\end{align}
In view of \eqref{eqn-ener-est}, the sequence $\{\mathcal{C}_r(\bu_{\varepsilon}^m)\}_m$ remains bounded in $\L^{\frac{r+1}{r}}(\mathcal{O})$, since
$\|\mathcal{C}_r(\bu_{\varepsilon}^m)\|_{\L^{\frac{r+1}{r}}}
=\|\bu_{\varepsilon}^m\|_{\L^{r+1}}^r\leq C$; hence, the generalized Lebesgue dominated convergence theorem (Brezis--Lions lemma; see Lemma \ref{Lem-Lions}) gives
$$\mathcal{C}_r(\bu_{\varepsilon}^m) \rightharpoonup \mathcal{C}_r(\bu_{\varepsilon}) \ \text{ in } \ \L^{\frac{r+1}{r}}(\mathcal{O}).
$$
Therefore, we deduce 
\begin{align}\label{eqn-c-conv}
c_r^0(\bu_{\varepsilon}^m,\boldsymbol{e}_k)	=\langle \mathcal{C}_r(\bu_{\varepsilon}^m),\boldsymbol{e}_k\rangle\to 	\langle \mathcal{C}_r(\bu_{\varepsilon}),\boldsymbol{e}_k\rangle = c_r^0(\bu_{\varepsilon},\boldsymbol{e}_k)	\ \text{ as }\ m\to\infty. 
\end{align}
A similar argument shows that 
$$
c_q^0(\bu_{\varepsilon}^m,\boldsymbol{e}_k)	\to c_q^0(\bu_{\varepsilon},\boldsymbol{e}_k)	\ \text{ as }\ m\to\infty. 
$$
The convergence \eqref{eqn-conv-5} and Lebesgue's dominated convergence theorem imply
\begin{align}
(g\rho_{\varepsilon}(\bu_{\varepsilon\boldsymbol{\tau}}^m),\boldsymbol{e}_{k\boldsymbol{\tau}})_{\mathbb{L}^2(\Gamma_{\mathrm{S}})}\to (g\rho_{\varepsilon}(\bu_{\varepsilon\boldsymbol{\tau}}),\boldsymbol{e}_{k\boldsymbol{\tau}})_{\mathbb{L}^2(\Gamma_{\mathrm{S}})}\ \text{ as }\ m\to\infty. 
\end{align}
Using the above limits, 	passing limit as $m\to\infty$ in \eqref{eqn-Galerkin}, we arrive at 
\begin{equation*}
	\begin{aligned}
	\mu a(\bu_{\varepsilon},\boldsymbol{e}_k)+b(\bu_{\varepsilon},\bu_{\varepsilon},\boldsymbol{e}_k)+\alpha(\bu_{\varepsilon},\boldsymbol{e}_k)&+\beta c_r^0(\bu_{\varepsilon},\boldsymbol{e}_k)-\kappa c_q^0(\bu_{\varepsilon},\boldsymbol{e}_k)\\+(g\rho_{\varepsilon}(\bu_{\varepsilon\boldsymbol{\tau}}),\boldsymbol{e}_{k\boldsymbol{\tau}})_{\mathbb{L}^2(\Gamma_{\mathrm{S}})}&=\langle\f, \boldsymbol{e}_k\rangle,
	\end{aligned}
\end{equation*}
for all $k=1,2,\ldots.$ Therefore, for all $\bv\in \mathbb{V}_{\mathrm{div}}\cap\mathbb{L}^{r+1}(\mathcal{O})$, we have 
\begin{equation}\label{eqn-limit}
		\begin{aligned}
			\mu a(\bu_{\varepsilon},\bv)+b(\bu_{\varepsilon},\bu_{\varepsilon},\bv)+\alpha(\bu_{\varepsilon},\bv)&+\beta c^0_r(\bu_{\varepsilon},\bv)-\kappa c^0_q(\bu_{\varepsilon},\bv)\\+(g\rho_{\varepsilon}(\bu_{\varepsilon\boldsymbol{\tau}}),\bv_{\boldsymbol{\tau}})_{\mathbb{L}^2(\Gamma_{\mathrm{S}})}&=\langle\f, \bv\rangle,
		\end{aligned}
\end{equation}
which means that $\bu_{\varepsilon}$ is a weak solution of \eqref{eqn-problem-5}. The weakly lowersemicontinuity property of norms implies 
	\begin{align}
	&\frac{\mu}{2}\|\nabla\bu_{\varepsilon}\|_0^2+\alpha\|\bu_{\varepsilon}\|_0^2+\frac{\beta}{2}\|\bu_{\varepsilon}\|_{\mathbb{L}^{r+1}}^{r+1}\nonumber\\&\leq 	\frac{\mu}{2}\liminf_{m\to\infty}\|\nabla\bu_{\varepsilon}^m\|_0^2+\alpha\liminf_{m\to\infty}\|\bu_{\varepsilon}^m\|_0^2+\frac{\beta}{2}\liminf_{m\to\infty}\|\bu_{\varepsilon}^m\|_{\mathbb{L}^{r+1}}^{r+1}\nonumber\\&\leq\liminf_{m\to\infty}\left[\frac{\mu}{2}\|\nabla\bu_{\varepsilon}^m\|_0^2+\alpha\|\bu_{\varepsilon}^m\|_0^2+\frac{\beta}{2}\|\bu_{\varepsilon}^m\|_{\mathbb{L}^{r+1}}^{r+1}\right]
	\leq \mathcal{K},
	\end{align}
	which completes the proof. 
	\end{proof}
	
\subsection{Existence and uniqueness of solutions to Problem \ref{problem-3}} \label{subsec:wellp}
Since the regularized variational inequality problem \eqref{eqn-problem-4} is equivalent to problem \eqref{eqn-problem-5}, Proposition \ref{prop-existence-ueps} ensures that problem  \eqref{eqn-problem-4} possesses a solution  $\bu_{\varepsilon}\in\mathbb{V}_{\mathrm{div}}\cap \cap\mathbb{L}^{r+1}(\mathcal{O})$, which is uniformly bounded in $\mathbb{V}_{\mathrm{div}}\cap \cap\mathbb{L}^{r+1}(\mathcal{O})$. By passing to the limit as $\varepsilon\to 0$ in \eqref{eqn-problem-4}, we arrive at the following theorem establishing the existence of a solution $\bu$ to the variational inequality problem \eqref{eqn-problem-3}.
			
\begin{theorem}\label{thm-exists-unique}
	Given $\f\in\mathbb{V}^{\prime}$ and $g\in L^2(\Gamma_{\mathrm{S}})$,  there exists a weak solution $\bu$ of Problem \ref{problem-3} and the following estimate holds:
	\begin{equation}\label{eqn-ener-2}
	\frac{\mu}{2}\|\nabla\bu\|_0^2+\alpha\|\bu\|_0^2+\frac{\beta}{2}\|\bu\|_{\mathbb{L}^{r+1}}^{r+1}\le\mathcal{K},
	\end{equation}
	where $\mathcal{K}$ is as in \eqref{calK}.
	Moreover, for $r>3$ if 
	\begin{align}\label{eqn-uni-1}
	\alpha>  \eta_{1,r}+\eta_{2,r}+\eta_{3,r},
	\end{align}
	where 
	\begin{align*}
	\eta_{1,r}&=\left(\frac{1}{4\mu}\right)^{\frac{r-1}{r-3}}\left(\frac{r-3}{r-1}\right)\left(\frac{8}{\beta (r-1)}\right)^{\frac{2}{r-3}},\\	\eta_{2,r}&=\left(\frac{r-q}{r-1}\right)\left(\frac{2(q-1)}{\beta(r-1)}\right)^{\frac{q-1}{r-q}}\left( \kappa q2^{q-1}\right)^{\frac{r-1}{r-q}},\\
	\eta_{3,r}&=\left(\frac{r-q}{r-1}\right)\left(\frac{4(q-1)}{\beta(r-1)}\right)^{\frac{q-1}{r-q}}\left( \kappa q2^{q-1}\right)^{\frac{r-1}{r-q}}. 
	\end{align*}
	or 
	\begin{align}\label{eqn-uni-2}
	 \alpha> \eta_{2,r}+\eta_{3,r}+\frac{1}{4\mu}\ \text{ and }\ \beta\mu\geq 1,
	\end{align}
	then the uniqueness follows. If $r=3$ and $\beta\mu\geq 1$, then the uniqueness follows. If $r\in[1,3)$, then the uniqueness follows provided
	\begin{align}\label{eqn-uni-3}
	\alpha> \eta_4\frac{C_g^{\frac{8}{4-d}}2^{\frac{d}{4-d}}}{\alpha\mu^{\frac{d}{4-d}}}\mathcal{K}^{\frac{4}{4-d}} \ \text{ where }\ \eta_4=C_g^{\frac{8}{4-d}}\left(\frac{8}{4-d}\right)\left(\frac{4+d}{8\theta\mu}\right)^{\frac{4+d}{4-d}}.
	\end{align}
	\end{theorem}
	\begin{proof}
	We prove the existence and uniqueness of weak solution to problem \eqref{eqn-problem-3} in the following steps:
	\vskip 0.1cm
	\noindent\textbf{Step 1:} \emph{Existence.} 
	Let $\bu_{\varepsilon}\in\mathbb{V}_{\mathrm{div}}\cap\mathbb{L}^{r+1}(\mathcal{O})$ be a weak solution of the regularized problem \eqref{eqn-problem-4} which satisfies the estimate \eqref{eqn-ener-1}. Once again an application of the Banach-Alaoglu theorem yields the existence of a subsequence $\{\bu_{\varepsilon_k}\}_{k=1}^{\infty}$ such that the following convergences hold: 
	\begin{equation}\label{eqn-conv-6}
	\left\{
	\begin{aligned}
	\bu_{\varepsilon_k} &\xrightarrow{w} \bu\  \text{ in }\ \mathbb{V}_{\mathrm{div}},\
	\bu_{\varepsilon_k} \xrightarrow{w} \bu\ \text{ in }\  \L^{r+1}(\mathcal{O}),\\
	\bu_{\varepsilon_k}& \to \bu\  \text{ in }\ \mathbb{H}_{\mathrm{div}},\  	\bu_{\varepsilon_k} \to \bu \ \text{ in }\ \mathbb{L}^2(\Gamma), \ \text{ as }\ k\to\infty. 
	\end{aligned}
	\right.
	\end{equation}
	Moreover, along a further subsequence denoted by the same symbol such that 
	\begin{equation}\label{eqn-conv-7}
		\left\{
		\begin{aligned}
		\bu_{\varepsilon_k}(\bx)&\to \bu(\bx)\ \text{ for a.e. }\ \bx\in\mathcal{O},\\	\bu_{\varepsilon_k}(\bx) &\to \bu(\bx)\ \text{ for a.e. }\ \bx\in\Gamma, \  \text{ as }\  k\to\infty. 
		\end{aligned}
		\right.
	\end{equation}
Our next aim is to show that $\bu$ is a weak solution of the variational inequality problem \eqref{eqn-problem-3}. Using the weakly lower semicontinuous property of the norm, we have 
\begin{align}\label{eqn-a-lim}
\limsup_{k\to\infty}a(\bu_{\varepsilon_k},\bv-\bu_{\varepsilon_k})&=\limsup_{k\to\infty}(\nabla\bu_{\varepsilon_k},\nabla(\bv-\bu_{\varepsilon_k}))\nonumber\\&=(\nabla\bu,\nabla\bv)-\liminf_{k\to\infty}\|\nabla \bu_{\varepsilon_k}\|_0^2\nonumber\\&\leq (\nabla\bu,\nabla\bv)-\|\nabla\bu\|_0^2=a(\bu,\bv-\bu). 
\end{align}
Using the first weak convergence given in \eqref{eqn-conv-6}, we infer 
\begin{align}\label{eqn-f-lim}
\lim_{k\to\infty}\langle\f,\bv-\bu_{\varepsilon_k}\rangle= \langle\f,\bv-\bu\rangle. 
\end{align}
Moreover, using \eqref{eqn-j-est}, we have 
\begin{align}\label{eqn-j-lim}
	\lim_{k\to\infty}j_{\varepsilon_k}(\bv_{\boldsymbol{\tau}})=j(\bv_{\boldsymbol{\tau}}). 
\end{align}
Let us now consider the convergence of trilinear form. A calculation similar to \eqref{eqn-b-conv} yields 
\begin{align}\label{eqn-b-lim}
&|b(\bu_{\varepsilon_k},\bu_{\varepsilon_k},\bv-\bu_{\varepsilon_k})-b(\bu,\bu,\bv-\bu)|\nonumber\\&\leq |b(\bu_{\varepsilon_k},\bu_{\varepsilon_k},\bv)-b(\bu,\bu,\bv)|+|b(\bu_{\varepsilon_k},\bu_{\varepsilon_k},\bu_{\varepsilon_k})-b(\bu,\bu,\bu)|\nonumber\\&\leq |b(\bu_{\varepsilon_k}-\bu,\bv,\bu_{\varepsilon_k})|+|b(\bu,\bv,\bu_{\varepsilon_k}-\bu)|
\nonumber\\&\leq  C_g^2C_p^{1-\frac{d}{4}}\left(\|\nabla\bu_{\varepsilon_k}\|_{0}+\|\nabla\bu\|_{0}\right)^{1+\frac{d}{4}} \|\bu_{\varepsilon_k}-\bu\|_{0}^{1-\frac{d}{4}} \|\nabla\bv\|_0\nonumber\\&\to 0\ \text{ as }\ k\to\infty.
\end{align}
Next, we consider the convergence of the nonlinear damping and pumping terms.  Arguing similarly as in \eqref{eqn-c-conv}, we get 
\begin{align}\label{eqn-c-conv-1}
\lim_{k\to\infty}c_r^0(\bu_{\varepsilon_k},\bv)=c_r^0(\bu,\bv)\ \text{ and }\ 	\lim_{k\to\infty}c_q^0(\bu_{\varepsilon_k},\bv)=c_q^0(\bu,\bv),\ \text{ for all }\ \bv\in\L^{r+1}(\mathcal{O}). 
\end{align}
Using \eqref{eqn-c-conv-1} and the weakly lower semicontinuous property of the norm, we deduce 
\begin{align}
\limsup_{k\to\infty}c_r^0(\bu_{\varepsilon_k},\bv-\bu_{\varepsilon_k})&=\limsup_{k\to\infty}\left[c_r^0(\bu_{\varepsilon_k},\bv)-c_r^0(\bu_{\varepsilon_k},\bu_{\varepsilon_k})\right]\nonumber\\&=c_r^0(\bu,\bv)-\liminf_{k\to\infty}\|\bu_{\varepsilon_k}\|_{\L^{r+1}}^{r+1}\nonumber\\&\le c_r^0(\bu,\bv)-\|\bu\|_{\L^{r+1}}^{r+1}=c_r^0(\bu,\bv-\bu). 
\end{align}
Let us now consider 
\begin{align}
&	|c_q^0(\bu_{\varepsilon_k},\bv-\bu_{\varepsilon_k})-c_q^0(\bu,\bv-\bu)|\nonumber\\&\leq|c_q^0(\bu_{\varepsilon_k},\bv)-c^0(\bu,\bv)|+|c_q^0(\bu_{\varepsilon_k},\bu_{\varepsilon_k})-c_q^0(\bu,\bu)| \nonumber\\&\leq|c_q^0(\bu_{\varepsilon_k},\bv)-c_q^0(\bu,\bv)|+|\|\bu_{\varepsilon_k}\|_{\L^{q+1}}^{q+1}-\|\bu\|_{\L^{q+1}}^{q+1}|=:I_{1k}+I_{2k}. 
\end{align}
Note that $I_{1k}\to 0$ as $k\to\infty$ by using \eqref{eqn-c-conv-1}. We estimate $I_{2k}$
by applying H\"older's inequality together with suitable interpolation inequalities, and by exploiting the convergence property \eqref{eqn-conv-6}, as follows:
\begin{align*}
I_{2k}&=\int_{\mathcal{O}}\left[|\bu_{\varepsilon_k}(\bx)|^{q+1}-|\bu(\bx)|^{q+1}\right]\d\bx
\nonumber\\&=\int_{\mathcal{O}}\int_0^1\frac{\d}{\d\theta}|\theta\bu_{\varepsilon_k}(\bx)+(1-\theta)\bu(\bx)|^{q+1}\d\theta \d\bx\nonumber\\&= (q+1)\int_{\mathcal{O}}\int_0^1 |\theta\bu_{\varepsilon_k}(\bx)+(1-\theta)\bu(\bx)|^{q-1}\left(\theta\bu_{\varepsilon_k}(\bx)+(1-\theta)\bu(\bx)\right)\cdot\left(\bu_{\varepsilon_k}(\bx)-\bu(\bx)\right)\d\theta \d\bx\nonumber\\&\leq (q+1)\int_{\mathcal{O}}\left[|\bu_{\varepsilon_k}(\bx)|+|\bu(\bx)|\right]^{q}|\bu_{\varepsilon_k}(\bx)-\bu(\bx)|\d\bx\nonumber\\&\leq (q+1)\left(\|\bu_{\varepsilon_k}\|_{\L^{q+1}}+\|\bu\|_{\L^{q+1}}\right)^{q}\|\bu_{\varepsilon_k}-\bu\|_{\L^{q+1}}\nonumber\\&\leq (q+1)|\mathcal{O}|^{\frac{q(r-q)}{(r+1)(q+1)}}\left(\|\bu_{\varepsilon_k}\|_{\L^{r+1}}+\|\bu\|_{\L^{r+1}}\right)^q\left(\|\bu_{\varepsilon_k}\|_{\L^{r+1}}+\|\bu\|_{\L^{r+1}}\right)^{\frac{(r+1)(q-1)}{(r-1)(q+1)}}\|\bu_{\varepsilon_k}-\bu\|_{0}^{\frac{2(r-q)}{(r-1)(q+1)}}\nonumber\\&\to 0\ \text{ as }\ k\to\infty. 
\end{align*}
Therefore, we immediately have 
\begin{align}\label{eqn-c-lim}
\lim_{k\to\infty} c_q^0(\bu_{\varepsilon_k},\bv-\bu_{\varepsilon_k})=c_q^0(\bu,\bv-\bu). 
\end{align}
Using \eqref{eqn-j-est}, \eqref{eqn-def-j} and \eqref{eqn-conv-6}, we find 
\begin{align*}
|j_{\varepsilon_k}(\bu_{\varepsilon_k\boldsymbol{\tau}})-j(\bu_{\boldsymbol{\tau}})|&\leq|j_{\varepsilon_k}(\bu_{\varepsilon_k\boldsymbol{\tau}})-j(\bu_{\varepsilon_k\boldsymbol{\tau}})|+|j(\bu_{\varepsilon_k\boldsymbol{\tau}})-j(\bu_{\tau})|\nonumber\\&\leq \frac{\varepsilon_k}{2}\|g\|_{L^1(\Gamma_{\mathrm{S}})}+\int_{\Gamma_{\mathrm{S}}}g\left[|\bu_{\varepsilon_k\boldsymbol{\tau}}-\bu_{\tau}|\right]\d\mathrm{S}\nonumber\\&\leq \frac{\varepsilon_k}{2}\|g\|_{L^1(\Gamma_{\mathrm{S}})}+\|g\|_{L^2(\Gamma_{\mathrm{S}})}\|\bu_{\varepsilon_k\boldsymbol{\tau}}-\bu_{\tau}\|_{\mathbb{L}^2(\Gamma_{\mathrm{S}})}\nonumber\\&\to 0\ \text{ as }\ k\to\infty,
\end{align*}
so that 
\begin{align}\label{eqn-j-lim-1}
\lim_{k\to\infty}	j_{\varepsilon_k}(\bu_{\varepsilon_k\boldsymbol{\tau}})=j(\bu_{\boldsymbol{\tau}}). 
\end{align}
By substituting $\varepsilon_k$ in place of $\varepsilon$ in \eqref{eqn-problem-4} and taking $\limsup$ as $k \to +\infty$, we use the  limits \eqref{eqn-a-lim}-\eqref{eqn-j-lim-1} to conclude that the limit function $\bu$ satisfies the variational inequality problem, and therefore is a weak solution of \eqref{eqn-problem-3}. Moreover, the estimate \eqref{eqn-ener-2} is obtained directly from \eqref{eqn-ener-1} by using the weak lowersemicontinuity of norms.
			
\vskip 0.1cm
\noindent\textbf{Step 2:} \emph{Uniqueness.} 
Let us now prove the uniqueness of the weak solution. Let $\bu_1$ and $\bu_2$ be two weak solutions of \eqref{eqn-problem-3} with the forcing $\f$ and denote $\bw=\bu_1-\bu_2$. 
Taking $\bv=\bu_2\in\mathbb{V}_{\mathrm{div}}\cap\mathbb{L}^{r+1}(\mathcal{O})$ in the inequality \eqref{eqn-problem-3} for $\bu=\bu_1$, we find 
\begin{equation}\label{eqn-sol-1}
\begin{aligned}
\mu a(\bu_1,\bu_2-\bu_1)+b(\bu_1,\bu_1,\bu_2-\bu_1)+\alpha(\bu_1,\bu_2-\bu_1)&+\beta c_r^0(\bu_1,\bu_2-\bu_1)\\-\kappa c_q^0(\bu_1,\bu_2-\bu_1)+j(\bu_{2\boldsymbol{\tau}})-j(\bu_{1\boldsymbol{\tau}})&\geq \langle\f, \bu_2-\bu_1\rangle.
\end{aligned}
\end{equation}
Similarly,
\begin{equation}\label{eqn-sol-2}
\begin{aligned}
\mu a(\bu_2,\bu_1-\bu_2)+b(\bu_2,\bu_2,\bu_1-\bu_2)+\alpha(\bu_2,\bu_1-\bu_2)&+\beta c_r^0(\bu_2,\bu_1-\bu_2)\\-\kappa c_q^0(\bu_2,\bu_1-\bu_2)+j(\bu_{1\boldsymbol{\tau}})-j(\bu_{2\boldsymbol{\tau}})&\geq \langle\f, \bu_1-\bu_2\rangle.
\end{aligned}
\end{equation}
Adding  \eqref{eqn-sol-1} and \eqref{eqn-sol-1}, we deduce 
\begin{align}\label{eqn-sol-3}
	&	\mu a(\bw,\bw)+\alpha(\bw,\bw)+\beta\langle\mathcal{C}_r(\bu_1)-\mathcal{C}_r(\bu_2),\bu_1-\bu_2\rangle\nonumber\\&\leq-\left[\langle\mathcal{B}(\bu_1,\bu_1)-\mathcal{B}(\bu_2,\bu_2),\bw\rangle\right] +\kappa \langle\mathcal{C}_q(\bu_1)-\mathcal{C}_q(\bu_2),\bw\rangle\nonumber\\&=-\langle\mathcal{B}(\bw,\bu_2),\bw\rangle+\kappa \langle\mathcal{C}_q(\bu_1)-\mathcal{C}_q(\bu_2),\bw\rangle,
\end{align}
where we have used the fact that $\langle\B(\bu_1,\bw),\bw\rangle=0$. 	
					From \eqref{2.23}, we have  
					\begin{align}\label{eqn-ces}
						&\langle\mathcal{C}_r(\bu_1)-\mathcal{C}_r(\bu_2),\bu_1-\bu_2\rangle\geq \frac{1}{2}\||\bu_1|^{\frac{r-1}{2}}(\bu_1-\bu_2)\|_0^2+\frac{1}{2}\||\bu_2|^{\frac{r-1}{2}}(\bu_1-\bu_2)\|_0^2. 
					\end{align}
Now, we distinguish three cases depending on the value of $r$.\\
						\noindent	\textbf{Case 1:} $3<r<\infty$.
						Let us  estimate the term $|\langle\mathcal{B}(\bw,\bu_2),\bw\rangle|$ using    H\"older's and Young's inequalities as 
					\begin{align}\label{eqn-bdes}
							|\langle\mathcal{B}(\bw,\bu_2),\bw\rangle|&=|\langle\mathcal{B}(\bw,\bw),\bu_2\rangle|
					\leq\|\nabla\bw\|_0\|  |\bu_2|\bw \|_0\nonumber\\&\leq  \theta\mu\|\bw\|_{\V}^2+\frac{1}{4\theta\mu}\|  |\bu_2|\bw  \|_0^2,
					\end{align}
				for some $0<\theta\leq 1$. When $r>3$, an application of Hölder's and Young's inequalities gives the estimate of $\frac{1}{2\mu}\||\bu_2|\bw\|_0^2$  as (\cite[Theorem 2.5]{SGMTM}, \cite{MT2})
					\begin{align}\label{eqn-des}
						\frac{1}{4\theta\mu}	\||\bu_2|\bw\|_0^2&=	\frac{1}{4\theta\mu}\int_{\mathcal{O}}|\bu_2(\bx)|^2|\bw(\bx)|^2\d\bx\nonumber\\&=\frac{1}{4\theta\mu}\int_{\mathcal{O}}|\bu_2(\bx)|^2|\bw(\bx)|^{\frac{4}{r-1}}|\bw(\bx)|^{\frac{2(r-3)}{r-1}}\d\bx\nonumber\\&\leq\frac{1}{4\theta\mu}\left(\int_{\mathcal{O}}|\bu_2(\bx)|^{r-1}|\bw(\bx)|^2\d\bx\right)^{\frac{2}{r-1}}\left(\int_{\mathcal{O}}|\bw(\bx)|^2\d\bx\right)^{\frac{r-3}{r-1}}\nonumber\\&\leq\frac{\beta}{4}\int_{\mathcal{O}}|\bu_2(\bx)|^{r-1}|\bw(\bx)|^2\d\bx+\eta_{1,r}\int_{\mathcal{O}}|\bw(\bx)|^2\d\bx,
					\end{align}
					where 
					\begin{align}\label{eqn-rho-1}
						\eta_{\theta,r}=\left(\frac{1}{4\theta\mu}\right)^{\frac{r-1}{r-3}}\left(\frac{r-3}{r-1}\right)\left(\frac{8}{\beta (r-1)}\right)^{\frac{2}{r-3}}.
					\end{align}
					 Using \eqref{eqn-des} in \eqref{eqn-bdes}, we find 
					\begin{align}\label{eqn-bes}
					|\langle\mathcal{B}(\bw,\bu_2),\bw\rangle|&\leq\theta\mu\|\nabla\bw\|_{0}^2+\frac{\beta}{4}\||\bu_2|^{\frac{r-1}{2}}\bw\|_0^2+\eta_{1,r}\|\bw\|_0^2.
					\end{align}
					An application of  Taylor's formula (\cite[Theorem 7.9.1]{PGC}) yields 
					\begin{align}\label{eqn-est-c0-1}
							\kappa|\langle\mathcal{C}_q(\bu_1)-\mathcal{C}_q(\bu_2),\bw\rangle|&= \kappa\bigg|\bigg<\int_0^1\mathcal{C}_q^{\prime}(\theta\bu_1+(1-\theta)\bu_2)\d\theta\bw,\bw\bigg>\bigg|\nonumber\\&\leq \kappa q2^{q-1}\bigg<\int_0^1|\theta\bu_1+(1-\theta)\bu_2|^{q-1}\d\theta|\bw|,|\bw|\bigg>\nonumber\\&\leq \kappa q2^{q-1}\left<\left(|\bu_1|^{q-1}+|\bu_2|^{q-1}\right)|\bw|,|\bw|\right>\nonumber\\&=   \kappa q2^{q-1}\||\bu_1|^{\frac{q-1}{2}}\bw\|_0^2+ \kappa q2^{q-1}\||\bu_2|^{\frac{q-1}{2}}\bw\|_0^2. 
					\end{align}
					Using H\"older's inequality, we estimate $ \kappa q2^{q-1}\||\bu_1|^{\frac{q-1}{2}}\bw\|_0^2$ as 
					\begin{align}\label{eqn-est-c0-2}
						 \kappa q2^{q-1}\||\bu_1|^{\frac{q-1}{2}}\bw\|_0^2&= \kappa q2^{q-1}\int_{\mathcal{O}}|\bu_1(\bx)|^{q-1}|\bw(\bx)|^2\d\bx \nonumber\\&= \kappa q2^{q-1}\int_{\mathcal{O}}|\bu_1(\bx)|^{q-1}|\bw(\bx)|^{\frac{2(q-1)}{r-1}}|\bw(\bx)|^{\frac{2(r-q)}{r-1}}\d\bx \nonumber\\&\leq \kappa q2^{q-1}\bigg(\int_{\mathcal{O}}|\bu_1(\bx)|^{r-1}|\bw(\bx)|^2\d\bx\bigg)^{\frac{q-1}{r-1}}\bigg(\int_{\mathcal{O}}|\bw(\bx)|^2\d\bx\bigg)^{\frac{r-q}{r-1}}\nonumber\\&\leq\frac{\beta}{2}\int_{\mathcal{O}}|\bu_1(\bx)|^{r-1}|\bw(\bx)|^2\d\bx+\eta_{2,r}\int_{\mathcal{O}}|\bw(\bx)|^2\d\bx,
					\end{align}
					where 
					\begin{align}\label{eqn-rho-2}
						\eta_{2,r}=\left(\frac{r-q}{r-1}\right)\left(\frac{2(q-1)}{\beta(r-1)}\right)^{\frac{q-1}{r-q}}\left( \kappa q2^{q-1}\right)^{\frac{r-1}{r-q}}. 
					\end{align}
					A similar calculation yields 
					\begin{align}\label{eqn-est-c0-3}
						& \kappa q2^{q-1}\||\bu_2|^{\frac{q-1}{2}}\bw\|_0^2\leq\frac{\beta}{4}\int_{\mathcal{O}}|\bu_2(\bx)|^{r-1}|\bw(\bx)|^2\d\bx+\eta_{3,r}\int_{\mathcal{O}}|\bw(\bx)|^2\d\bx,
					\end{align}
					where 
					\begin{align}\label{eqn-rho-3}
						\eta_{3,r}=\left(\frac{r-q}{r-1}\right)\left(\frac{4(q-1)}{\beta(r-1)}\right)^{\frac{q-1}{r-q}}\left( \kappa q2^{q-1}\right)^{\frac{r-1}{r-q}}. 
					\end{align}
					Using \eqref{eqn-est-c0-2} and \eqref{eqn-est-c0-3} in \eqref{eqn-est-c0-1}, we deduce 
					\begin{align}\label{eqn-est-c0-4}
						\kappa\langle\mathcal{C}_q(\bu_1)-\mathcal{C}_q(\bu_2),\bw\rangle|&\leq\frac{\beta}{2}\||\bu_1|^{\frac{r-1}{2}}\bw\|_0^2+\frac{\beta}{4}\||\bu_2|^{\frac{r-1}{2}}\bw\|_0^2+(\eta_{2,r}+\eta_{3,r})\|\bw\|_0^2. 
					\end{align}
					Combining \eqref{eqn-sol-3}, \eqref{eqn-ces}, \eqref{eqn-bes} and \eqref{eqn-est-c0-4}, we finally get for $r>3,$
					\begin{align}\label{eqn-final-est}
					(1-\theta)\mu\|\nabla\bw\|_0^2+\left(\alpha-\eta_r\right)\|\bw\|_0^2\leq0, 
					\end{align}
					where $\eta_r=\eta_{\theta,r}+\eta_{2,r}+\eta_{3,r}$, and $\eta_{1,r}$, $\eta_{2,r}$ and $\eta_{3,r}$ are defined in \eqref{eqn-rho-1}, \eqref{eqn-rho-2} and \eqref{eqn-rho-3}, respectively.  Therefore, for $$(1-\theta)C_p\mu +\alpha>\eta_r,$$ where $C_p$ is the Poincar\'e constant, the uniqueness follows. By taking $\theta=1$, one obtains the condition \eqref{eqn-uni-1}.

					One can estimate $|\langle\B(\bw,\bw),\bu_2\rangle|$ in the following way also: 
					\begin{align}\label{2.26}
						&	|\langle\B(\bw,\bw),\bu_2\rangle|\nonumber\\&\leq\|\nabla\bw\|_{0}\||\bu_2|\bw\|_0
					\leq\theta\mu \|\nabla\bw\|_{0}^2+\frac{1}{4\theta\mu }\int_{\mathcal{O}}|\bu_2(\bx)|^2|\bw(\bx)|^2\d\bx\nonumber\\&= \theta\mu \|\nabla\bw\|_{0}^2+\frac{1}{4\theta\mu} \int_{\mathcal{O}}|\bw(\bx)|^2\left(|\bu_2(\bx)|^{r-1}+1\right)\frac{|\bu_2(\bx)|^2}{|\bu_2(\bx)|^{r-1}+1}\d\bx\nonumber\\&\leq \theta\mu \|\nabla\bw\|_{0}^2+\frac{1}{4\theta\mu}\int_{\mathcal{O}}|\bu_2(\bx)|^{r-1}|\bw(\bx)|^2\d\bx+\frac{1}{4\theta\mu}\int_{\mathcal{O}}|\bw(\bx)|^2\d\bx,
					\end{align}
					where 	we have used the fact that $\left\|\frac{|\bu_2|^2}{|\bu_2|^{r-1}+1}\right\|_{\L^{\infty}}<1$, for $r\geq 3$. Using \eqref{eqn-ces},  \eqref{eqn-est-c0-3} and \eqref{2.26} in \eqref{eqn-sol-3}, we deduce 
					\begin{align}
				(1-\theta)\mu\|\nabla\bw\|_{0}^2+\left(\alpha-\eta_{2,r}-\eta_{3,r}-\frac{1}{4\theta\mu}\right)\|\bw\|_{0}^2+\left(\frac{\beta}{4}-\frac{1}{4\theta\mu}\right)\||\bu_2|^{\frac{r-1}{2}}\bw\|_0^2\leq 0,
					\end{align}
					and the uniqueness follows, provided $$(1-\theta)C_p\mu + \alpha> \eta_{2,r}+\eta_{3,r}+\frac{1}{4\theta\mu}\ \text{ and }\ \beta\mu\geq 1.$$  By taking $\theta=1$, we arrive at the condition \eqref{eqn-uni-2}.
					
					\vskip 0.1cm
					\noindent	\textbf{Case 2:} \emph{$r=3$ with $\beta\mu\geq 1$.} 
					From \eqref{2.23}, we infer
					\begin{align}\label{231}
						\beta\langle\mathcal{C}_r(\bu_1)-\mathcal{C}_r(\bu_2),\bw\rangle\geq\frac{\beta}{2}\||\bu_1|\bw\|_0^2+\frac{\beta}{2}\||\bu_2|\bw\|_0^2. 
					\end{align}
					We estimate  $|\langle\mathcal{B}(\bw,\bu_2),\bw\rangle|$ using H\"older's and Young's inequalities as 
					\begin{align}\label{232}
					 |\langle\mathcal{B}(\bw,\bu_2),\bw\rangle|\leq\||\bu_2|\bw\|_0\|\nabla\bw\|_{0} \leq\frac{1}{\beta} \|\nabla\bw\|_{0}^2+\frac{\beta}{4}\||\bu_2|\bw\|_0^2.
					\end{align}
					Combining \eqref{eqn-sol-3}, \eqref{eqn-est-c0-4}, \eqref{231} and \eqref{232}, we deduce 
					\begin{align}\label{eqn-critical}
					\left(\mu-\frac{1}{\beta}\right)\|\nabla\bw\|_0^2+\alpha\|\bw\|_0^2\leq 0,
					\end{align}
				and the uniqueness follows,	provided $\beta\mu \geq 1$. 
					\vskip 0.1cm
					\noindent		\textbf{Case 3:} \emph{$r\in[1,3)$.} 
					For this case, one can estimate $|\langle\B(\bw,\bw),\bu_2\rangle|$ using  H\"older's, Gagliardo-Nirenberg's and Young's inequalities  as 
					\begin{align}\label{2.21}
						|\langle\B(\bw,\bw),\bu_2\rangle|&\leq \|\nabla\bw\|_0\|\bw\|_{\L^4}\|\bu_2\|_{\L^4}
						\nonumber\\&\leq C_g\|\nabla\bw\|_{0}\|\bw\|_0^{1-\frac{d}{4}}\|\nabla\bw\|_{0}^{\frac{d}{4}}\|\bu_2\|_{\L^4}
						\nonumber\\&\leq C_g\|\bw\|_0^{1-\frac{d}{4}}\|\nabla\bw\|_{0}^{1+\frac{d}{4}}\|\bu_2\|_{\L^4}\nonumber\\&\leq \theta\mu\|\nabla\bw\|_{0}^2+\eta_4\|\bu_2\|_{\L^4}^{\frac{8}{4-d}}\|\bw\|_0^2,
					\end{align}
					where 
					\begin{align}
						\eta_4=C_g^{\frac{8}{4-d}}\left(\frac{8}{4-d}\right)\left(\frac{4+d}{8\theta\mu}\right)^{\frac{4+d}{4-d}}.
					\end{align}
					Applying Gagliardo-Nirenberg's inequality, we estimate 
					\begin{align*}
						\|\bu_2\|_{\L^4}^{\frac{8}{4-d}}&\leq \left(C_g\|\bu_2\|_0^{1-\frac{d}{4}}\|\nabla\bu_2\|_{0}^{\frac{d}{4}}\right)^{\frac{8}{4-d}}=C_g^{\frac{8}{4-d}}\|\bu_2\|_0^2\|\nabla\bu_2\|_{0}^{\frac{2d}{4-d}}\nonumber\\&\leq C_g^{\frac{8}{4-d}}\frac{\mathcal{K}}{\alpha}\left(\frac{2\mathcal{K}}{\mu}\right)^{\frac{d}{4-d}}=\frac{C_g^{\frac{8}{4-d}}2^{\frac{d}{4-d}}}{\alpha\mu^{\frac{d}{4-d}}}\mathcal{K}^{\frac{4}{4-d}},
					\end{align*}
					where we use  \eqref{eqn-ener-2} and $\mathcal{K}$ is defined in \eqref{calK}. 
					Using \eqref{eqn-ces}, \eqref{eqn-est-c0-4} and \eqref{2.21} in \eqref{eqn-sol-3}, we deduce 
					\begin{align}
						(1-\theta)\mu\|\nabla\bw\|_0^2+\left(\alpha-\eta_4\frac{C_g^{\frac{8}{4-d}}2^{\frac{d}{4-d}}}{\alpha\mu^{\frac{d}{4-d}}}\mathcal{K}^{\frac{4}{4-d}}\right)\|\bw\|_0^2\leq 0,
					\end{align}
					and the uniqueness follows provided $$(1-\theta)C_p\mu +\alpha> \eta_4\frac{C_g^{\frac{8}{4-d}}2^{\frac{d}{4-d}}}{\alpha\mu^{\frac{d}{4-d}}}\mathcal{K}^{\frac{4}{4-d}}.$$  By taking $\theta=1$, we get the condition \eqref{eqn-uni-1}.
			\end{proof}
			
			\begin{remark}
				We note that when 
				$r\geq 3$, the parameter conditions are independent of the forcing term 
				$\f$ (see \eqref{eqn-uni-1} and \eqref{eqn-uni-2}). In contrast, for 
				$1\leq r<3$, uniqueness holds only under an additional smallness assumption on 
				$\f$ (see \eqref{eqn-uni-3}).
			\end{remark}
	
		\subsection{Existence and uniqueness of solutions to Problem \ref{problem-2}} \label{subsec:comWell}
	We next state the inf-sup condition (also referred to as the Ladyzhenskaya-Babuška-Brezzi (LBB) condition; see \cite[Theorem 2.2]{JSHNJW}), \cite[Lemma 2.2]{NSa}: there exists a constant $\vartheta>0$  such that 
		\begin{align}\label{eqn-inf-sup}
			\vartheta\|q\|_{0}\leq 	\sup_{\bv\in\mathbb{V}_0}\frac{d(\bv,q)}{\|\nabla\bv\|_{0}}, \ \text{ for all }\ q\in Q,
		\end{align}
	where $\mathbb{V}_0:=\mathbb{H}_0^1(\mathcal{O})$. 	We now turn to the analysis of Problem \ref{problem-2}.  From this point onward, we assume that $1 \leq r \leq \tfrac{2d}{d-2}$ for $d > 2$, and $1 \leq r < \infty$ for $d=2$. For brevity, we express this as $1 \leq r \leq \tfrac{2d}{(d-2)^+}$. With this assumption in place, the Sobolev embedding theorem yields the continuous embedding $\H^1(\mathcal{O}) \hookrightarrow \L^{r+1}(\mathcal{O})$. Consequently, Problem \ref{problem-2} reduces to finding $\bu \in \V$ and $p\in Q$.
		
	The preceding assumption will be instrumental in establishing the following result.
		\begin{theorem}\label{thm-equivalent}
		Suppose that the assumptions of Theorem \ref{thm-exists-unique} hold and $1\leq r\leq\frac{2d}{(d-2)^+}$. Then, Problem \ref{problem-2} has a unique solution.		
		\end{theorem}
		\begin{proof}
			Let $\bu\in\V_{\mathrm{div}}$ be  the unique solution to Problem 	\ref{problem-3}. Then   for all $ \bv\in \V_{0,\mathrm{div}}$, we infer from \eqref{eqn-problem-3} that 
			\begin{align*}
				\mu a(\bu,\bv)+b(\bu,\bu,\bv)+\alpha(\bu,\bv)+\beta c_r^0(\bu,\bv) -\kappa c_q^0(\bu,\bv)= \langle\f, \bv\rangle. 
			\end{align*}
			The inf-sup condition \eqref{eqn-inf-sup}, together with the standard functional-analytic result stated in \cite[Theorem 2.2]{JSHNJW}, guarantees the existence of $p\in Q$ for which
			\begin{align}\label{eqn-pressure-1}
				\mu a(\bu,\bv)+b(\bu,\bu,\bv)+\alpha(\bu,\bv)+\beta c_r^0(\bu,\bv)-\kappa c_q^0(\bu,\bv)+d(\bv,p)= \langle\f, \bv\rangle 
			\end{align}
			holds for every $\bv\in\V_0$.
			Now, let $\bv\in\V$ be arbitrary. By virtue of the inf-sup condition \eqref{eqn-inf-sup}, there exists $\bv_1\in\V_0$ satisfying
			\begin{align}\label{eqn-pressure-2}
				 d(\bv_1, q) = d(\bv, q)  \  \text{ for all } \ q \in Q.
			\end{align}
			Let $\bv_2 = \bv - \bv_1$. Then $\bv_2 \in \V$, and it satisfies
			$$
			 d(\bv_2, q) = 0 \ \text{ for all } \ q \in Q.
			$$
			It follows immediately from the preceding relation that $\bv_2\in\V_{\mathrm{div}}$. Taking $\bv_2$ as the test function in \eqref{eqn-problem-3}, we obtain
			\begin{align}
				\mu a(\bu,\bv-\bv_1-\bu)&+b(\bu,\bu,\bv-\bv_1-\bu)+\alpha(\bu,\bv-\bv_1-\bu)+\beta c_r^0(\bu,\bv-\bv_1-\bu)\nonumber\\&-\kappa c_q^0(\bu,\bv-\bv_1-\bu)+j(\boldsymbol{v}_{\boldsymbol{\tau}}-\boldsymbol{v}_{1\boldsymbol{\tau}})-j(\bu_{\boldsymbol{\tau}})\geq \langle\f, \bv-\bv_1-\bu\rangle,
			\end{align}
			It is important to observe that $\bv_1=0$ on $\Gamma_1$, and hence, $j(\bv_{\boldsymbol{\tau}}-\bv_{1,\boldsymbol{\tau}}) = j(\bv_{\boldsymbol{\tau}})$. Substitution of \eqref{eqn-pressure-1} and \eqref{eqn-pressure-2} into the preceding inequality yields
			\begin{align}
					&\mu a(\bu,\bv-\bu)+b(\bu,\bu,\bv-\bu)+\alpha(\bu,\bv-\bu)+\beta c_r^0(\bu,\bv-\bu)-\kappa c_q^0(\bu,\bv-\bu)\nonumber\\&\quad+j(\boldsymbol{v}_{\boldsymbol{\tau}})-j(\bu_{\boldsymbol{\tau}})\nonumber\\&\geq  \langle\f, \bv-\bv_1-\bu\rangle+ 	\mu a(\bu,\bv_1)+b(\bu,\bu,\bv_1)+\alpha(\bu,\bv_1)+\beta c_r^0(\bu,\bv_1)-\kappa c_q^0(\bu,\bv_1)
					\nonumber\\&=-d(\bv_1,p)+\langle\f,\bv-\bu\rangle
					\nonumber\\&=-d(\bv-\bu,p)+\langle\f,\bv-\bu\rangle, \ \text{ for all }\ \bv\in\V,
			\end{align}
		since $\bu\in\V_{\mathrm{div}}$. Thus, the first relation in \eqref{eqn-problem-2} is satisfied. The second relation follows immediately from $\bu\in\V_{\mathrm{div}}$.
		
		It remains to establish an estimate for the pressure variable $p$. To this end, we invoke the inf-sup condition \eqref{eqn-inf-sup} and the identity \eqref{eqn-pressure-1}, which together give
			\begin{align}
				\vartheta\|p\|_{Q}&\leq 	\sup_{\bv\in\V_0}\frac{1}{\|\bv\|_{\V}}[\langle\f,\bv\rangle -\mu a(\bu,\bv)-b(\bu,\bu,\bv)-\alpha (\bu,\bv)-\beta c_r^0(\bu,\bv)+\kappa c_q^0(\bu,\bv)]
		\nonumber\\&\leq\|\f\|_{\V^{\prime}}+\mu\|\nabla\bu\|_{0}+C_b\|\nabla\bu\|_{0}^2+\alpha\|\bu\|_{0}+C_s\beta \|\bu\|_{\L^{r+1}}^r+C_s\kappa|\mathcal{O}|^{\frac{r-q}{(q+1)(r+1)}}\|\bu\|_{\L^{r+1}}^q\nonumber\\&\leq  \|\f\|_{\V^{\prime}}+\sqrt{2\mu\mathcal{K}}+\frac{2C_b\mathcal{K}}{\mu}+\sqrt{\alpha\mathcal{K}}+C_s\beta\mathcal{K}^{\frac{r}{r+1}}+C_s\kappa|\mathcal{O}|^{\frac{r-q}{r+1}}\mathcal{K}^{\frac{q}{r+1}}.
			\end{align}
			In deriving the above estimate, we have used \eqref{eqn-ener-2} and the Sobolev embedding
			\begin{align}
				\|\bv\|_{\L^{r+1}}\leq C_s\|\bv\|_{\V}\ \text{ for all }\ \bv\in\V,
			\end{align}
			which is valid for $1\leq r\leq\frac{2d}{(d-2)^+}$.

			It remains to prove the uniqueness of the solution $(\bu,p)$ to Problem \ref{problem-2}. Since the velocity field $\bu$ satisfies \eqref{eqn-problem-3}, uniqueness of the solution to Problem \ref{problem-3} immediately determines $\bu$ uniquely. Let $(\bu,p_1)$ and $(\bu,p_2)$ be two solutions of Problem \ref{eqn-problem-2} associated with the same forcing term $\f$. From \eqref{eqn-pressure-1}, we obtain
			$$
		d(\bv, p_1 - p_2) = 0 \  \text{ for all }\  \bv \in \V_0. 
			$$
			Taking $q=p_1-p_2$ in the inf-sup condition \eqref{eqn-inf-sup}, we conclude that $\|p_1 - p_2\|_Q = 0$. Hence, $p_1 = p_2$ in $Q$.
		\end{proof}

\section{Discontinuous Galerkin Methods}\label{sec4}\setcounter{equation}{0} 
For clarity of presentation, we consider $\mathcal{O}$ to be an open, bounded polygonal domain. The domain is partitioned into a family of triangulations $\{\mathcal{T}_h\}_{h>0},$ and for each $h>0,$ $\mathcal{T}_h= \{K\}$, which are locally quasi-uniform, and shape-regular
 \cite[Section 3.2]{PGC1}. Here, $h$ denotes the mesh parameter. 
For two distinct neighboring elements $K_i, K_j \in \mathcal{T}_h$, denote by
$e = \partial K_i \cap \partial K_j$ $(i \neq j)$
their common interface. The diameters of an element $K$ and of an edge $e$ are written as $h_K$ and $h_e$, respectively. We define:
\begin{itemize}
	\item $\mathscr{E}_h$: the collection of all edges in the triangulation $\mathcal{T}_h$,
	\item  $\mathscr{E}_h^I$: the set of all interior edges,
	\item  $\mathscr{E}_h^{\mathrm{S}}$: the set of all edges lying on the boundary portion $\overline{\Gamma}_{\mathrm{S}}$, and
	\item $\mathscr{E}_h^\ast := \mathscr{E}_h^I \cup \mathscr{E}_h^{\mathrm{S}}$.
\end{itemize}
For vectors $\boldsymbol{v}$ and $\boldsymbol{n}$, let $\boldsymbol{v} \otimes \boldsymbol{n}$ denote the matrix whose $(i, j)$-th entry is $v_i n_j$. For two $d\times d$ matrix-valued variables $\boldsymbol{A}$ and $\boldsymbol{B}$, we define
$$
\boldsymbol{A} : \boldsymbol{B} = \sum_{i,j=1}^{d} A_{ij} B_{ij}.
$$
Let $e = \partial K_1 \cap \partial K_2$, and let $\boldsymbol{n}_1$ and $\boldsymbol{n}_2$ be the unit normal vectors on $e$ pointing outward from $K_1$ and $K_2$, respectively.
We define the average $\{\cdot\}$ and the jump $[\cdot]$ on $e$ for a scalar $q$, a vector $\boldsymbol{v}$, and a matrix $\boldsymbol{A}$ as follows:
\begin{align*}
	\{q\}&=\frac{1}{2}\left(q|_{\partial K_1}+q|_{\partial K_2}\right), \ [q]=q|_{\partial K_1}\boldsymbol{n}_1+q|_{\partial K_2}\boldsymbol{n}_2,\\
	\{\boldsymbol{v}\}&=\frac{1}{2}\left(\boldsymbol{v}|_{\partial K_1}+\boldsymbol{v}|_{\partial K_2}\right), \ [\boldsymbol{v}]=\boldsymbol{v}|_{\partial K_1}\cdot\boldsymbol{n}_1+\boldsymbol{v}|_{\partial K_2}\cdot\boldsymbol{n}_2,\\
	\{\boldsymbol{A}\}&=\frac{1}{2}\left(\boldsymbol{A}|_{\partial K_1}+\boldsymbol{A}|_{\partial K_2}\right), \ [\boldsymbol{A}]=\boldsymbol{A}|_{\partial K_1}\boldsymbol{n}_1+\boldsymbol{A}|_{\partial K_2}\boldsymbol{n}_2.
\end{align*}
We also define a matrix-valued jump $\llbracket\cdot\rrbracket$ for a vector $\boldsymbol{v}$ by
$$
\llbracket\boldsymbol{v}\rrbracket = \boldsymbol{v}|_{\partial K_1} \otimes \boldsymbol{n}_1 + \boldsymbol{v}|_{\partial K_2} \otimes \boldsymbol{n}_2
$$
on $e$.
If $e$ lies on the boundary $\partial \mathcal{O}$, the above definitions are modified accordingly:
\begin{align*}
	\{q\}=q, \ 	\{\boldsymbol{v}\}=	\boldsymbol{v},\ \{\boldsymbol{A}\}=\boldsymbol{A}\ \text{ and } \ [q]=q\boldsymbol{n},\ [\boldsymbol{v}]=	\boldsymbol{v}\cdot \boldsymbol{n},\ [\boldsymbol{A}]=\boldsymbol{A}\boldsymbol{n},\ \llbracket\boldsymbol{v}\rrbracket=\boldsymbol{v}\otimes\boldsymbol{n}. 
\end{align*}
It follows immediately from a straightforward calculation that
\begin{align*}
	&\sum_{K\in\mathcal{T}_h}\int_{\partial K}q\boldsymbol{v}\cdot\boldsymbol{n}\d\mathrm{S}=\int_{\mathscr{E}_h^{I}}[q]\cdot\{\boldsymbol{v}\}\d\mathrm{S}+\int_{\mathscr{E}_h}\{q\}[\boldsymbol{v}]\d\mathrm{S},\\
	&\sum_{K\in\mathcal{T}_h}\int_{\partial K}\boldsymbol{v}\cdot\boldsymbol{A}\boldsymbol{n}\d\mathrm{S}=\int_{\mathscr{E}_h^{I}}[\boldsymbol{A}]\cdot\{\boldsymbol{v}\}\d\mathrm{S}+\int_{\mathscr{E}_h}\{\boldsymbol{A}\}:\llbracket\boldsymbol{v}\rrbracket\d\mathrm{S}. 
\end{align*}
Let $k \geq 1$ be an integer. We define the finite element space for the velocity as
\begin{equation}\label{eq:Vh}
\mathbb{V}_h = \{\boldsymbol{v} \in \mathbb{L}^2(\mathcal{O}) : \boldsymbol{v}|_K \in [P_k(K)]^d \  \text{ for all } \ K \in \mathcal{T}_h \},
\end{equation}
and the finite element space for the pressure as
\begin{equation}\label{eq:Qh}
Q_h = \{q \in L^2_0(\mathcal{O}) : q|_K \in P_{k-1}(K) \ \text{ for all }\ K \in \mathcal{T}_h \},
\end{equation}
where $P_k(K)$ denotes the set of all polynomials of degree at most $k$ on the element $K$. We assume that these spaces satisfy the standard inverse inequalities (see \cite[Section 3.2]{PGC1}).

\subsection{DG methods for the variational inequality}
In this subsection, we assume that $\bu \in \V\cap \H^2(\mathcal{O})$  and $p\in Q\cap H^1(\mathcal{O})$.
We define the standard interior penalty forms for the CBFeD system. The forms $a_h^*(\cdot, \cdot)$, $d_h(\cdot, \cdot)$, $j(\cdot)$, and $\mathcal{F}_h(\cdot)$ correspond, respectively, to the DG discretization of the viscous term, the divergence term, the friction term, and the right-hand side of the system \eqref{eqn-CBF}:
\begin{align*}
	a_h^*(\bv,\bw)&=\sum_{K\in\mathcal{T}_h}\int_K\nabla\bv:\nabla\bw\d\bx - \int_{\mathscr{E}_h^*}\{\nabla\bv\}:\llbracket\bw\rrbracket\d\S,\\
	\widehat{\alpha}_h(\bv,\bw)&=\sum_{e\in\mathscr{E}_h^*}\gamma h_e^{-1}\int_e\llbracket\bv\rrbracket:\llbracket\bw\rrbracket\d\S,\\
	d_h(\bv,q)&=-\sum_{K\in\mathcal{T}_h}\int_Kq\nabla\cdot\bv\d\bx+\int_{\mathscr{E}_h^*}\{q\}[\bv]\d\S,\\
	j(\bv_{\boldsymbol{\tau}})&=\int_{\mathscr{E}_h^{\S}}g|\bv_{\boldsymbol{\tau}}|\d\S,\ \mathcal{F}_h(\bv)=\sum_{K\in\mathcal{T}_h}\int_{K}\f\cdot\bv\d\bx,
\end{align*}
where  $\gamma>0$ is a parameter that will be specified later. 

As discussed in \cite[Section 3.2]{FJWHWY}, we now briefly outline the derivation of the consistency term in $a_h(\cdot, \cdot)$ and the origin of the inequality. For $\bu \in \H^2(\mathcal{O}),$ $\nabla\bu$ is in $H^1(\mathcal{O}; \mathbb{R}^{d\times d})$, and we have $[\nabla \bu] = 0$ on the interior edges \cite{PErnDG}. For an arbitrary $\bv \in \V_h$, multiplying \eqref{eqn-CBF}$_1$
by $\bv-\bu$ and then integrating over $\mathcal{O}$,  we obtain
\begin{align}\label{eqn-dg-form}
-&\mu	\int_{\mathcal{O}}\Delta\bu\cdot(\bv-\bu)\d\bx +
\int_{\mathcal{O}}	(\bu \cdot \nabla)\bu \cdot(\bv-\bu) \d\bx +  \alpha \int_{\mathcal{O}}\bu \cdot(\bv-\bu) \d\bx + \beta|\bu|^{r-1}\bu \cdot(\bv-\bu) \d\bx \nonumber\\&\quad- \kappa \int_{\mathcal{O}}|\bu|^{q-1}\bu  \cdot(\bv-\bu)\d\bx +\int_{\mathcal{O}} \nabla p \cdot(\bv-\bu)\d\bx  =  \int_{\mathcal{O}}\f \cdot(\bv-\bu)\d\bx. 
	\end{align}
For the first term in \eqref{eqn-dg-form}, integrating over an element $K$,  applying integration by parts, and then summing over all elements, we obtain 
\begin{align}\label{eqn-dg-form-1}
	-\mu	\int_{\mathcal{O}}\Delta\bu\cdot(\bv-\bu)\d\bx&=-\mu\sum_{K\in\mathcal{T}_h}\int_K \Delta\bu\cdot(\bv-\bu)\d\bx\nonumber\\&=\mu\sum_{K\in\mathcal{T}_h}\int_K\nabla\bu:\nabla(\bv-\bu)\d\bx-\mu\sum_{K\in\mathcal{T}_h}\int_{\partial K}\frac{\partial\bu}{\partial\bn}\cdot(\bv-\bu)\d\S. 
\end{align}
Let us rewrite the edge integral term in \eqref{eqn-dg-form-1} as 
\begin{align}\label{eqn-dg-form-2}
	-\mu\sum_{K\in\mathcal{T}_h}\int_{\partial K}\frac{\partial\bu}{\partial\bn}\cdot(\bv-\bu)\d\S&=-\mu\int_{\mathscr{E}_h}\{\nabla\bu\}:\llbracket\bv-\bu\rrbracket\d\S-\mu\int_{\mathscr{E}_h^I}[\nabla\bu]\cdot\{\bv-\bu\}\d\S\nonumber\\&=-\mu\int_{\mathscr{E}_h^*}\{\nabla\bu\}:\llbracket\bv-\bu\rrbracket\d\S-\mu\int_{\mathscr{E}_h^{\S}}\frac{\partial\bu}{\partial \bn}\cdot(\bv-\bu)\d\S. 
\end{align}
For the second term in \eqref{eqn-dg-form-2}, using the definition of 
$\boldsymbol{\sigma}_{\boldsymbol{\tau}}$ and the boundary conditions \eqref{eqn-boundary}, we obtain
\begin{align}\label{eqn-dg-form-3}
	-\mu\int_{\mathscr{E}_h^{\S}}\frac{\partial\bu}{\partial \bn}\cdot(\bv-\bu)\d\S&= -\mu\int_{\mathscr{E}_h^{\S}}\left(\frac{\partial u_{\boldsymbol{n}}}{\partial \bn}\bn+\frac{\partial\bu_{\boldsymbol{\tau}}}{\partial \bn}\right)\cdot\left[(v_{\bn}-u_{\bn})\bn+(\bv_{\boldsymbol{\tau}}-\bu_{\boldsymbol{\tau}})\right]\d\S\nonumber\\&=-\mu\int_{\mathscr{E}_h^{\S}}\frac{\partial u_{\boldsymbol{n}}}{\partial \bn}(v_{\bn}-u_{\bn})\d\S-\mu\int_{\mathscr{E}_h^{\S}}\frac{\partial\bu_{\boldsymbol{\tau}}}{\partial \bn}\cdot(\bv_{\boldsymbol{\tau}}-\bu_{\boldsymbol{\tau}})\d\S\nonumber\\&= -\int_{\mathscr{E}_h^{\S}}\mu\frac{\partial\bu_{\boldsymbol{\tau}}}{\partial \bn}\cdot(\bv_{\boldsymbol{\tau}}-\bu_{\boldsymbol{\tau}})\d\S= \int_{\mathscr{E}_h^{\S}}\boldsymbol{\sigma}_{\boldsymbol{\tau}}(\bu)\cdot(\bu_{\boldsymbol{\tau}}-\bv_{\boldsymbol{\tau}})\d\S\nonumber\\&\leq \int_{\mathscr{E}_h^{\S}}g(|\bv_{\boldsymbol{\tau}}|-|\bu_{\boldsymbol{\tau}}|)\d\S. 
\end{align}
Using \eqref{eqn-dg-form-2} and \eqref{eqn-dg-form-3} in \eqref{eqn-dg-form-1}, we get
\begin{align}\label{eqn-dg-form-4}
		-\mu	\int_{\mathcal{O}}\Delta\bu\cdot(\bv-\bu)\d\bx&\leq \mu\sum_{K\in\mathcal{T}_h}\int_K\nabla\bu:\nabla(\bv-\bu)\d\bx -\mu\int_{\mathscr{E}_h^*}\{\nabla\bu\}:\llbracket\bv-\bu\rrbracket\d\S\nonumber\\&\quad+\int_{\mathscr{E}_h^{\S}}g(|\bv_{\boldsymbol{\tau}}|-|\bu_{\boldsymbol{\tau}}|)\d\S.
\end{align}

We employ the following variant of the Lesaint-Raviart upwinding scheme \cite{PLPAR}, introduced in \cite{VGBMRMFW}, to discretize the nonlinear convection term in \eqref{eqn-dg-form}. The superscripts int and ext denote, respectively, the traces of a function on a side of  $K$  taken from the interior and exterior of the edge  $e$. When the side of  $K$  lies on the boundary  $\partial\mathcal{O}$, the exterior trace is set to zero. Let us define, for all $\bz_h,\bu_h,\bv_h,\bw_h\in\V_h$,
\begin{align*}
	&b_h(\bz_h;\bu_h;\bu_h,\bw_h)\nonumber\\&:=\sum_{K\in\mathcal{T}_h}\int_K\left((\bu_h\cdot\nabla)\bv_h\cdot\bw_h+\frac{1}{2}(\operatorname{div}\bu_h)\bv_h\cdot\bw_h\right)\d\bx-\frac{1}{2}\int_{\mathscr{E}_h^*}[\bu_h]\{\bv_h\cdot\bw_h\}\d\S\nonumber\\&\quad+\sum_{K\in\mathcal{T}_h}\int_{\partial K_{-}^{\bz_h}\backslash\mathscr{E}_h^{\S}}|\{\bu_h\}\cdot\bn|\left(\bv_h^{\mathrm{int}}-\bv_h^{\mathrm{ext}}\right)\cdot\bw_h^{\mathrm{int}}\d\S,
\end{align*}
and 
\begin{align*}
	b_h^{NL}(\bz_h;\bu_h;\bv_h,\bw_h)=\sum_{K\in\mathcal{T}_h}\int_{\partial K_{-}^{\bz_h}\backslash\mathscr{E}_h^{\S}}|\{\bu_h\cdot\bn\}|\left(\bv_h^{\mathrm{int}}-\bv_h^{\mathrm{ext}}\right)\cdot\bw_h^{\mathrm{int}}\d\S,
\end{align*}
where 
\begin{align*}
	\partial K_{-}^{\bz_h}=\left\{\bx\in\partial K:\bz_h(\bx)\cdot\bn<0\right\},
\end{align*}
and the superscript $\bz_h$ indicates the dependence of $\partial K_{-}^{\bz_h}$ on $\bz_h$.

The nonlinear reaction terms are purely local (no inter-element coupling), and so their DG counterparts are simply sums of element integrals:
\begin{align*}
	c_{rh}(\bu_h;\bv_h)=\sum_{K\in\mathcal{T}_h}\int_K|\bu_h|^{r-1}\bu_h\cdot\bv_h\d\bx,\\
	c_{qh}(\bu_h;\bv_h)=\sum_{K\in\mathcal{T}_h}\int_K|\bu_h|^{q-1}\bu_h\cdot\bv_h\d\bx.
\end{align*}

Let $a_h(\cdot, \cdot)$ denote one of the DG forms corresponding to interior penalty scheme (IPG method including incomplete (IIPG) and symmetric (SIPG)/non-symmetric (NIPG) interior penalty), LDG method etc.\ of the diffusion term. A detailed formulation of the DG method can be found in \cite[pp.\ 894--895]{FJWHWY} and the references therein. For brevity, we do not repeat it here and instead state the following DG formulation:

\begin{problem} \label{eq:prb-discform}
Find $(\bu_h, p_h)\in \mathbb{V}_h \times Q_h $ such that 
\begin{equation*}
	\begin{aligned}
	& \mu a_h(\bu_h  ,\bv_h-\bu_h)  + b_h(\bu_h; \bu_h;\bu_h, \bv_h -\bu_h) +\alpha (\bu_h, \bv_h-\bu_h)+\beta c_{rh}(\bu_h; \bv_h-\bu_h) \\ & \quad \quad  - \kappa c_{qh}(\bu_h; \bv_h-\bu_h) + d_h(\bv_h -\bu_h, p_h)+ j(\bv_{h\tau}) - j(\bu_{h\tau}) \ge \mathcal{F}_h(\bv_h-\bu_h) \quad \text{for all } \bv_h \in \mathbb{V}_h, \\
	& d_h(\bu_h, q_h)=0 \quad \text{ for all } q_h \in Q_h.
	\end{aligned}
\end{equation*}
\end{problem}
In the next section, we discuss the existence and uniqueness of solution of the above problem, and in Section~\ref{sec6} an error estimate of the solutions of this problem and Problem~\ref{problem-2} is established.

\section{Stability, Existence and Uniqueness}\label{sec5}\setcounter{equation}{0}
In this section, we discuss the stability and solvability results for the DG formulation of the CBFeD model. In contrast to the results in \cite[Section 4]{FJWHWY}, the main new difficulty is the presence of the nonlinear damping and drag terms
\begin{align*}
	\alpha \bu +\beta |\bu|^{r-1}\bu \ \text{ and } \ -\kappa |\bu|^{q-1}\bu,
\end{align*}
where $\alpha>0, \beta>0$ and $\kappa>0.$ To proceed further, define the broken Sobolev norms on $\mathbb{V}(h)=\mathbb{V}_h+\mathbb{H}^2(\mathcal{O})\cap \mathbb{V}$ and $Q(h)=Q_h+Q\cap H^1(\mathcal{O})$ as
\begin{align*}
	& |\bv|_{*, e}^2:= \sum_{e\in \mathscr{E}_h^\ast} \frac{1}{h_e}\| \llbracket \bv \rrbracket\|_{0, e}^2, \quad \tnorm{\bv}_1^2:= \sum_{K\in \mathcal{T}_h}|\bv|_{1,K}^2 + |\bv|_{*, e}^2,\\
	& \tnorm{\bv}^2:= \tnorm{\bv}_1^2+\sum_{K\in \mathcal{T}_h}h_K^2 |\bv|_{2,K}^2, \text{ and }\|q\|_{\ast}^2:=\|q\|_0^2 + \sum_{K\in \mathcal{T}_h}h_K^2 |q|_{1,K}^2,
\end{align*}
 where
 \begin{align*}
 	|\cdot|_{m,K}:=\left(\sum_{|k|= m}\int_{K}|\partial^{\boldsymbol{k}}\cdot(\boldsymbol{x})|^2d\boldsymbol{x}\right)^{1/2}. 
 \end{align*}
 
 One can easily check the following consistency condition satisfied by the solution $\bu$ of \eqref{eqn-problem-2}:
 \begin{equation}\label{eqn-problem-2-constnc}
 	\begin{aligned}
 		& \mu a_h(\bu,\bv_h-\bu)+b_h(\bu;\bu; \bu,\bv_h-\bu)+\alpha(\bu,\bv_h-\bu)+\beta c_{rh}(\bu;\bv_h-\bu)\\ 
 		& \qquad \quad -\kappa c_{qh}(\bu;\bv_h-\bu)+d_h(\bv_h-\bu,p)+j(\boldsymbol{v}_{h \boldsymbol{\tau}})-j(\bu_{\boldsymbol{\tau}})\geq \mathcal{F}_h(\bv_h-\bu),\\
 		& d_h(\bu,q_h)=0, 
 	\end{aligned}
 \end{equation}
 for all $ \bv_h \in \V(h), \ q_h \in Q_h.$
 
We next state the discrete version of inf-sup condition (also referred to as the Ladyzhenskaya-Babuška-Brezzi (LBB) condition) see \cite{FJWHWY, BMR06, BMR08}: there exists a constant $\vartheta>0$  such that 
 \begin{align}\label{eqn-dis-inf-sup}
 	\vartheta \|q_h\|_0\leq 	\sup_{\bv_h\in\widehat{\V}_h}\frac{d_h(\bv_h,q_h)}{\tnorm{\bv_h}_1 }, \ \text{ for all }\ q_h\in Q_h.
 \end{align}

 Next, we list properties satisfied by the discrete operators $a_h(\cdot, \cdot)$, $b_h(\cdot;\cdot;\cdot,\cdot)$ and $d_h(\cdot, \cdot)$:
 
\begin{property}[\cite{VGBMRMFW_MathComp, FJWHWY}] \label{prop:A}
	There exist positive constants $\gamma_0,  C, C_0,  C_1, $ and $C_2$ independent of $h$ such that 
  \begin{itemize}
 	\item[(a)] $a_h(\bv_h , \bv_h) \ge \gamma_0 \tnorm{\bv_h}_1^2 $  for all $ \bv_h \in \mathbb{V}_h$;
 	\item[(b)] $ |a_h(\bu_h, \bv_h)| \le C \tnorm{\bu_h} \tnorm{\bv_h}, \text{ and } |d_h(\bv_h, q_h)| \le C \tnorm{\bv_h}\|q_h\|_\ast$
 	for all $\bu_h, \ \bv_h \in \mathbb{V}_h,$ and for all $q_h\in Q(h)$;
 	\item[(c)] for all $\bu_h, \bv_h, \bw_h, \bz_h \in \mathbb{V}_h,$
 	\begin{align*}
 		&	b_h(\bu_h;\bu_h; \bv_h, \bv_h) \ge 0, \quad  | b_h(\bz_h; \bu_h; \bv_h, \bw_h)| \le C_0 \tnorm{\bu_h}_1 \tnorm{\bv_h}_1 \tnorm{\bw_h}_1,\\
 		& | b_h^{NL}(\bz_h; \bw_h; \bz_h, \bv_h) - b_h^{NL}(\bw_h; \bw_h; \bz_h, \bv_h) | \le C_1 \|\bz_h - \bw_h\|_1 \|\bz_h\|_1 \|\bv_h\|_1;	
 	\end{align*}
 	\item[(d)] and for all $\bv_h \in \V_h,$  
 	\begin{align*}
 		\left\vert \mathcal{F}(\bv_h) + j (\bv_{h\tau}) \right\vert \le C_2 \left( \|\f\|_0+ \|g\|_{L^2(\Gamma_{\mathrm S})}  \right).
 	\end{align*}
 \end{itemize}
\end{property}
For an explicit expression of $\gamma_0,$ we refer \cite[Lemma 1, Section 4]{FJWHWY}.

\begin{lemma}\label{lem:discoerc}
There exists a constant $C_3=C_3(\kappa,\beta, |\mathcal{O}|)>0$ independent of $h,$ such that 
\begin{align*}
	\mu a_h(\bv_h, \bv_h)+\alpha (\bv_h, \bv_h) + \beta c_{rh}(\bv_h;\bv_h) -\kappa c_{qh} (\bv_h; \bv_h) \ge \mu\gamma_0 \tnorm{\bv_h}_1^2 + \frac{\beta}{2}\|\bv_h\|_{\mathbb{L}^{r+1}}^{r+1}-C_3,
\end{align*}
for all $\bv_h \in \mathbb{V}_h.$
\end{lemma}
 
 \begin{proof}
 Note that
 \begin{align*}
 	\kappa c_{qh}(\bv_h; \bv_h)= \kappa \sum_{K\in\mathcal{T}_h}\int_K |\bv_h|^{q-1}\bv_h\cdot\bv_h\d\bx & = \kappa \sum_{K\in\mathcal{T}_h}\int_K |\bv_h|^{q+1}\d\bx, 
 \end{align*}
 and an application of H\"{o}lder's inequality leads to
 \begin{align*}
 \kappa	\sum_{K\in\mathcal{T}_h}\int_K |\bv_h|^{q+1} d\boldsymbol{x} \le \kappa \|\bv_h\|_{\mathbb{L}^{r+1}}^{q+1} |\mathcal{O}|^{\frac{r-q}{r+1}} \le \frac{\beta}{2} \|\bv_h\|_{\mathbb{L}^{r+1}}^{r+1} + C_3,
 \end{align*}
 for some $C_3=C_3(\kappa,\beta, |\mathcal{O}|)>0.$
 Therefore, utilizing this along with Property~\ref{prop:A}(a), we obtain 
 \begin{align*}
 	\mu a_h(\bv_h, \bv_h) & +\alpha (\bv_h, \bv_h) + \beta c_{rh}(\bv_h;\bv_h) -\kappa c_{qh} (\bv_h; \bv_h) \\
 	& = \mu a_h(\bv_h, \bv_h)+\alpha \|\bv_h\|^2 +\beta  \|\bv_h\|_{\mathbb{L}^{r+1}}^{r+1} - \kappa \|\bv_h\|_{\mathbb{L}^{q+1}}^{q+1} \\
 	& \ge \mu\gamma_0 \tnorm{\bv_h}_1^2 +\frac{\beta}{2}\|\bv_h\|_{\mathbb{L}^{r+1}}^{r+1} - C_3.
 \end{align*}
 This completes the proof.
 \end{proof}

Before proceeding further, let us define a subspace $\widehat{\V}_h$ of $\V_h$ as
\begin{align*}
	\widehat{\V}_h := \left\lbrace \bv_h \in \V_h : [\bv_h]_e=0 \ \text{ for all } \ e\in \mathscr{E}_h^{\ast}, \ \bv_h\mid_e=\boldsymbol{0}\  \text{ for all } \ e \in \mathscr{E}_h^{\S}  \right\rbrace.
\end{align*}
This subspace is useful to derive error estimates using Raviart-Thomas interpolation operator. Note that for any $\bv_h \in \widehat{\V}_h,$ the integral term on slip boundary is not included in this subspace. The polar set of $\widehat{\V}_h$ is defined as
\begin{align*}
	{\widehat{\V}_h}^o:=\left\lbrace \Phi \in \V_h'\ : \ \langle\Phi,\bv_h\rangle=0 \text{ for all } \bv_h \in \widehat{\V}_h  \right\rbrace.
\end{align*}

 \begin{theorem}\label{thm-dgstab}
 	Let $\f\in\L^2(\mathcal{O})$ and $g\in L^2(\Gamma_{\mathrm{S}})$ be such that 
 	\begin{align} \label{eq:thmassump1}
 	\mu\gamma_0-C_5>0\  \text{ and } \ 	\frac{\mathcal{R}_0}{\mu\gamma_0-C_5} (2C_0+C_1)<1,
 	\end{align}
 where $\gamma_0$ is the constant introduced in Property \ref{prop:A},
 \begin{align}\label{eqn-r0}
 	\mathcal{R}_0=\min \left\lbrace \frac{2C_2}{\mu\gamma_0} \left( \|\f\|_0 + \|g\|_{L^2(\Gamma_{\mathrm S})}\right) +2\sqrt{\frac{C_3}{\mu\gamma_0}},  \frac{\mu\gamma_0}{2C_0} \right\rbrace,
 \end{align}
 and 
 \begin{align*}
 	C_5= \frac{r-q}{r-1}\left(  \frac{2(q-1)}{\beta(r-1)} \right)^{\frac{q-1}{r-q}} \left( \kappa q2^{q-1} \right)^{\frac{r-1}{r-q}}.
 \end{align*}
 Then Problem~\ref{eq:prb-discform} admits a unique solution $(\bu_h, p_h) \in \mathcal{D}_h \times Q_h,$ with $\|p_h\|_0  \le C,$ for some $C=C (\|\f\|_0, \mathcal{R}_0),$ where
 \begin{align*}
 \mathcal{D}_h = \left\lbrace \bv_h \in \V_h\ : \ \tnorm{\bv_h}_1 \le \mathcal{R}_0 \right\rbrace.
 \end{align*}
 \end{theorem}

 \begin{proof}
 Define the space $\V_h^\sigma = \left\lbrace \bv_h \in \V_h  \, :\, d_h(\bv_h, q_h)=0 \text{ for all } q_h \in Q_h \right\rbrace.$ Then, we can find an equivalent problem corresponding to Problem \ref{eq:prb-discform} as: 
 find $\bu_h\in \V_h^\sigma $ such that 
 \begin{equation} \label{eq:prb-discform-div0}
 	\begin{aligned}
 		& \mu a_h(\bu_h  ,\bv_h-\bu_h)  + b_h(\bu_h; \bu_h;\bu_h, \bv_h -\bu_h) +\alpha (\bu_h, \bv_h-\bu_h)+\beta c_{rh}(\bu_h; \bv_h-\bu_h) \\ & \qquad \quad  -\kappa c_{qh}(\bu_h; \bv_h-\bu_h) + j(\bv_{h\tau}) - j(\bu_{h\tau}) \ge \mathcal{F}_h(\bv_h-\bu_h) \  \text{ for all }\  \bv_h \in \mathbb{V}_h^\sigma. 
 	\end{aligned}
 \end{equation}
We discuss the existence and uniqueness of the above problem in the following way. We first consider an auxiliary problem: for a given $\bu_h \in \V_h^\sigma,$ find $\bw_h \in \V_h^\sigma $ such that 
 \begin{equation} \label{eq:prb-discform-div0-wh}
	\begin{aligned}
		& \mu a_h(\bw_h  ,\bv_h-\bw_h)  +\alpha (\bw_h, \bv_h-\bw_h)+\beta c_{rh}(\bw_h; \bv_h-\bw_h)-\kappa c_{qh}(\bw_h; \bv_h-\bw_h) \\ & \qquad \quad   + j(\bv_{h\tau}) - j(\bw_{h\tau}) \ge \mathcal{F}_h(\bv_h-\bw_h) - b_h(\bu_h; \bu_h;\bu_h, \bv_h -\bw_h) \ \text{ for all } \ \bv_h \in \mathbb{V}_h^\sigma. 
	\end{aligned}
\end{equation}
Since $\mu a_h(\cdot, \cdot) + \alpha (\cdot, \cdot) +\beta c_{rh} (\cdot, \cdot) -\kappa c_{qh} (\cdot, \cdot)$ is coercive due to Lemma \ref{lem:discoerc}, and $j(\cdot),$ $\mathcal{F}_h(\cdot)$ and $b_h (\cdot; \cdot; \cdot, \cdot)$ are continuous, the above auxiliary problem admits a unique solution $\bw_h \in \V_h^\sigma$ for a given $\bu_h \in \V_h^\sigma.$ 
 
 Thus, we can define a map $\mathcal{S}_h: \V_h^\sigma \rightarrow \V_h^\sigma$ by $\mathcal{S}_h(\bu_h) = \bw_h.$ Now, to prove the existence of a unique solution $\bu_h$ of \eqref{eq:prb-discform-div0}, it is enough to show that $\mathcal{S}_h$ has a unique fixed point $\bu_h \in \V_h^\sigma$: $\mathcal{S}_h(\bu_h)=\bu_h$. To establish that, we will prove that map $\mathcal{S}_h$ is a self map and a contraction. To do so, let us define the ball $\mathcal{D}_h^\sigma \subset \V_h^\sigma$ as 
 \begin{align*}
 \mathcal{D}_h^\sigma  = \left\lbrace \bv_h \in \V_h^\sigma\ : \ \tnorm{\bv_h}_1 \le \mathcal{R}_0 \right\rbrace, 
 \end{align*}
where $R_0$ is defined in \eqref{eqn-r0}. 
Let $\bu_h \in \mathcal{D}_h^\sigma,$ that is, $\bu_h\in \V_h^\sigma$ with $\tnorm{\bu_h}_1\le \mathcal{R}_0$ be given. Then there exists a solution $\bw_h \in \V_h^\sigma$ of \eqref{eq:prb-discform-div0-wh} such that $\bw_h = \mathcal{S}_h(\bu_h).$ 
 Choosing $\bv_h=\boldsymbol{0}$ in \eqref{eq:prb-discform-div0-wh}, we have 
 \begin{equation*}
	\begin{aligned}
		 \mu a_h(\bw_h  ,\bw_h)  +\alpha (\bw_h, \bw_h) & +  \beta c_{rh}(\bw_h; \bw_h)  -\kappa c_{qh}(\bw_h; \bw_h)  + j(\bw_{h\tau}) \\
		 &\le \mathcal{F}_h(\bw_h) - b_h(\bu_h; \bu_h;\bu_h, \bw_h).
	\end{aligned}
\end{equation*}
Substituting $\bv_h=2\bw_h$ again in \eqref{eq:prb-discform-div0-wh}, we obtain 
 \begin{equation*}
	\begin{aligned}
		 \mu a_h(\bw_h  ,\bw_h)  +\alpha (\bw_h, \bw_h) & +\beta c_{rh}(\bw_h; \bw_h)  -\kappa c_{qh}(\bw_h; \bw_h) + j(\bw_{h\tau}) \\
		 & \ge \mathcal{F}_h(\bw_h) - b_h(\bu_h; \bu_h;\bu_h, \bw_h),
	\end{aligned}
\end{equation*}
which in combination with the above expression leads to 
 \begin{equation*}
	\begin{aligned}
		 \mu a_h(\bw_h  ,\bw_h)  +\alpha (\bw_h, \bw_h)&+\beta c_{rh}(\bw_h; \bw_h)  -\kappa c_{qh}(\bw_h; \bw_h) \\
		&=   \mathcal{F}_h(\bw_h)- j(\bw_{h\tau}) - b_h(\bu_h; \bu_h;\bu_h, \bw_h),
	\end{aligned}
\end{equation*}
for all $ \bv_h \in \mathbb{V}_h^\sigma.$ Utilizing Lemma \ref{lem:discoerc}, and continuity of $j(\cdot)$ and $\mathcal{F}(\cdot),$ one can easily obtain 
\begin{align*}
	\mu\gamma_0 \tnorm{\bw_h}_1^2 \le C_2 \left( \|\f\|_0 + \|g\|_{L^2(\Gamma_{\mathrm S})}\right) \tnorm{\bw_h}_1 + C_0 \tnorm{\bu_h}_1^2 \tnorm{\bw_h}_1+C_3,
\end{align*}
which further implies 
\begin{align*}
	\mu\gamma_0 \tnorm{\bw_h}_1^2 - \left( C_2 \left( \|\f\|_0 + \|g\|_{L^2(\Gamma_{\mathrm S})}\right)  + C_0 \tnorm{\bu_h}_1^2 \right) \tnorm{\bw_h}_1 - C_3 \le 0.
\end{align*}
Since $\mu, \gamma_0>0,$ we have 
\begin{align*}
	\tnorm{\bw_h}_1 & \le  \frac{  C_2 \left( \|\f\|_0 + \|g\|_{L^2(\Gamma_{\mathrm S})}\right)  + C_0 \tnorm{\bu_h}_1^2 + \sqrt{  \left( C_2 \left( \|\f\|_0 + \|g\|_{L^2(\Gamma_{\mathrm S})}\right)  + C_0 \tnorm{\bu_h}_1^2\right)^2 + 4 \mu\gamma_0 C_3  } }{2\mu\gamma_0}  \\
	& \le \frac{  C_2 \left( \|\f\|_0 + \|g\|_{L^2(\Gamma_{\mathrm S})}\right)  + C_0 \tnorm{\bu_h}_1^2 +   \left( C_2 \left( \|\f\|_0 + \|g\|_{L^2(\Gamma_{\mathrm S})}\right)  + C_0 \tnorm{\bu_h}_1^2\right) + 2\sqrt{\mu \gamma_0 C_3  } }{2\mu\gamma_0} \\
	& \le \frac{C_2 \left( \|\f\|_0 + \|g\|_{L^2(\Gamma_{\mathrm S})}\right)  + C_0 \tnorm{\bu_h}_1^2}{\mu\gamma_0} + \sqrt{\frac{C_3}{\mu\gamma_0}} \\
	& \le \mathcal{R}_0,
\end{align*}
due to the choice of $\mathcal{R}_0$ and as $\bu_h \in \mathcal{D}_h^\sigma.$ Thus, the map $\mathcal{S}_h: \mathcal{D}_h^\sigma \rightarrow \mathcal{D}_h^\sigma$ is a self map. 

We now establish the contraction property. To this end, for given
$\bu_{1h},\bu_{2h}\in\mathcal{D}_h^\sigma$, let
$\bw_{1h}=\mathcal{S}_h(\bu_{1h})$ and $\bw_{2h}=\mathcal{S}_h(\bu_{2h})$
denote the corresponding solutions of
\eqref{eq:prb-discform-div0-wh}, respectively. Thus, we have 
\begin{equation*} 
	\begin{aligned}
		& \mu a_h(\bw_{1h}  ,\bv_h-\bw_{1h})  +\alpha (\bw_{1h}, \bv_h-\bw_{1h})+ \beta c_{rh}(\bw_{1h}; \bv_h-\bw_{1h})- \kappa c_{qh}(\bw_{1h}; \bv_h-\bw_{1h}) \\ & \qquad \quad   + j(\bv_{h\tau}) - j(\bw_{1h\tau}) \ge \mathcal{F}_h(\bv_h-\bw_{1h}) - b_h(\bu_{1h}; \bu_{1h};\bu_{1h}, \bv_h -\bw_{1h}) \quad \text{ for all } \bv_h \in \mathbb{V}_h^\sigma, 
	\end{aligned}
\end{equation*}
and
\begin{equation*} 
	\begin{aligned}
		& \mu a_h(\bw_{2h}  ,\bv_h-\bw_{2h})  +\alpha (\bw_{2h}, \bv_h-\bw_{2h})+\beta c_{rh}(\bw_{2h}; \bv_h-\bw_{2h})- \kappa c_{qh}(\bw_{2h}; \bv_h-\bw_{2h}) \\ & \qquad \quad   + j(\bv_{h\tau}) - j(\bw_{2h\tau}) \ge \mathcal{F}_h(\bv_h-\bw_{2h}) - b_h(\bu_{2h}; \bu_{2h};\bu_{2h}, \bv_h -\bw_{2h}) \quad \text{ for all } \bv_h \in \mathbb{V}_h^\sigma. 
	\end{aligned}
\end{equation*}
Choose $\bv_h=\bw_{2h}$ in the first and $\bv_h=\bw_{1h}$ in second equation above, and then add them to obtain 
\begin{align*}
&	\mu a_h   (\bw_{1h}- \bw_{2h}, \bw_{1h}- \bw_{2h}) +\alpha (\bw_{1h}- \bw_{2h}, \bw_{1h}- \bw_{2h}) \\
	& \quad+ \beta \left( c_{rh} (\bw_{1h}, \bw_{1h}- \bw_{2h}) -  c_{rh}(\bw_{2h}, \bw_{1h}- \bw_{2h})\right) \nonumber\\&\quad - \kappa  \left( c_{qh} (\bw_{1h}, \bw_{1h}- \bw_{2h}) -  c_{qh}(\bw_{2h}, \bw_{1h}- \bw_{2h})\right) \\
	&    \le b_h(\bu_{1h}; \bu_{1h};\bu_{1h}, \bw_{2h} -\bw_{1h}) - b_h(\bu_{2h}; \bu_{2h};\bu_{2h}, \bw_{2h} -\bw_{1h}).
\end{align*}
Note that
\begin{align*}
\beta \left( c_{rh} (\bw_{1h}, \bw_{1h}- \bw_{2h}) - c_{rh}(\bw_{2h}, \bw_{1h}- \bw_{2h})\right) \ge \frac{\beta}{2} \int_\mathcal{O} \left( |\bw_{1h}|^{r-1} + |\bw_{2h}|^{r-1} \right) | \bw_{1h}- \bw_{2h}|^2 d\boldsymbol{x},
\end{align*}
and 
\begin{align*}
\kappa  & \left( |\bw_{1h}|^{q-1}- |\bw_{2h}|^{q-1}, \bw_{1h}- \bw_{2h}\right)  \\
& \le \kappa q 2^{q-1} \int_\mathcal{O} \left( |\bw_{1h}|^{q-1}  + |\bw_{2h}|^{q-1} \right) | \bw_{1h}- \bw_{2h}|^2 d\boldsymbol{x} \\
& \le \frac{\beta}{2} \int_\mathcal{O} \left( |\bw_{1h}|^{r-1} + |\bw_{2h}|^{r-1} \right) | \bw_{1h}- \bw_{2h}|^2 d\boldsymbol{x} + C_5 \tnorm{ \bw_{1h}- \bw_{2h}}_1^2,
\end{align*}
obtained by applying Young's inequality on the terms $\kappa q2^{q-1}|\bw_{1h}|^{q-1}$ and $\kappa q2^{q-1}|\bw_{2h}|^{q-1}$, where 
\begin{align*}
	C_5=\frac{r-q}{r-1}\left(  \frac{2(q-1)}{\beta(r-1)} \right)^{\frac{q-1}{r-q}} \left( \kappa q2^{q-1} \right)^{\frac{r-1}{r-q}}.
\end{align*}
Therefore, we have
\begin{align*}
-\kappa & \big( c_{qh} (\bw_{1h},  \bw_{1h}- \bw_{2h}) - c_{qh}(\bw_{2h}, \bw_{1h}- \bw_{2h})\big) \\
&  \ge - \frac{\beta}{2} \int_\mathcal{O} \left( |\bw_{1h}|^{r-1} + |\bw_{2h}|^{r-1} \right) | \bw_{1h}- \bw_{2h}|^2 d\boldsymbol{x} - C_5 \tnorm{ \bw_{1h}- \bw_{2h}}_1^2.
\end{align*}
Utilizing this, along with the coercivity of $a_h(\cdot, \cdot),$ one can easily deduce 
\begin{align*}
&	(\mu\gamma_0 - C_5) \tnorm{ \bw_{1h}- \bw_{2h}}_1^2 \nonumber\\& \le \left\vert b_h(\bu_{1h}; \bu_{1h};\bu_{1h}, \bw_{2h} -\bw_{1h}) - b_h(\bu_{2h}; \bu_{2h};\bu_{2h}, \bw_{2h} -\bw_{1h}) \right\vert\\
	& \le \left\vert b_h(\bu_{1h}; \bu_{2h};\bu_{1h}-\bu_{2h}, \bw_{2h} -\bw_{1h}) - b_h(\bu_{1h}; \bu_{1h}- \bu_{2h};\bu_{1h}, \bw_{2h} -\bw_{1h}) \right\vert \\
	& \quad + b_h^{\mathrm{NL}}(\bu_{1h}; \bu_{2h};\bu_{1h}, \bw_{2h} -\bw_{1h}) - b_h^{\mathrm{NL}}(\bu_{2h}; \bu_{2h};\bu_{1h}, \bw_{2h} -\bw_{1h}) \\
	& \le C_0  \tnorm{\bu_{1h}-\bu_{2h}}_1  \tnorm{\bw_{2h} -\bw_{1h}}_1 (\tnorm{\bu_{1h}}_1 + \tnorm{\bu_{2h}}_1) \\
	& \quad + C_1 \tnorm{\bu_{1h}-\bu_{2h}}_1 \tnorm{\bu_{1h}}_1 \tnorm{\bw_{2h} -\bw_{1h}}_1.
\end{align*}
Under assumption \eqref{eq:thmassump1}, this further implies
\begin{align*}
\tnorm{ \bw_{1h}- \bw_{2h}}_1 \le \frac{\mathcal{R}_0}{\mu\gamma_0-C_5} (2C_0+C_1) \tnorm{\bu_{1h}-\bu_{2h}}_1. 
\end{align*}
By the assumption \eqref{eq:thmassump1}, the map $\mathcal{S}_h: \V_h^\sigma \rightarrow \V_h^\sigma$ by $\mathcal{S}_h(\bu_h) = \bw_h$ on $\mathcal{D}_h^\sigma$ is a contraction. Thus, it has a unique fixed point due to the Banach fixed point theorem, that is, the discrete problem admits a unique solution.

Next, we aim to show the existence of a $p_h,$ and to do that we first define a map $\mathcal{G}_h: Q_h \rightarrow {\widehat{\V}_h}^o$ by 
\begin{align*}
	\langle \mathcal{G}_hq_h,\bv_h\rangle_{\widehat{\V}_h^{\prime},\widehat{\V}_h }=d_h(\bv_h, q_h) \text{ for all } \bv_h \in \widehat{\V}_h, \ q_h \in Q_h.
\end{align*} 
Let $\bu_h$ be solution as obtained in previous step, and set a map 
\begin{align*}
\Phi(\bv_h)=\mu a_h(\bu_h, \bv_h) +\alpha (\bu_h, \bv_h)+b_h(\bu_h;\bu_h;\bu_h, \bv_h) +\beta c_{rh}(\bu_h; \bv_h)-\kappa c_{qh}(\bu_h; \bv_h) - \mathcal{F}_h(\bv_h).
\end{align*}
Note that for any $\bv_h\in\widehat{\V}_h$, we have $\Phi(\bv_h)=0$. Hence, $\Phi\in{\widehat{\V}_h}^{o}$. Furthermore, since $\mathcal{G}_h:Q_h\to{\widehat{\V}_h}^{o}$ is an isomorphism and the discrete inf--sup condition \eqref{eqn-dis-inf-sup} holds, it follows from \cite[Theorem 2.2]{JSHNJW} that there exists
$p_h\in Q_h$ such that $\mathcal{G}_h p_h=\Phi,$ that is,
\begin{align*}
	d_h(\bv_h, p_h)= \mu a_h(\bu_h, \bv_h) +\alpha (\bu_h, \bv_h)+b_h(\bu_h;\bu_h;\bu_h, \bv_h) +\beta c_{rh}(\bu_h; \bv_h)-\kappa c_{qh}(\bu_h; \bv_h) - \mathcal{F}_h(\bv_h),
\end{align*}
for all $\bv_h \in \widehat{\V}_h.$ Also, the discrete inf-sup condition \eqref{eqn-dis-inf-sup} implies
\begin{align*}
	\vartheta \|p_h\|_0 & \le \sup_{\bv_h \in \widehat{\V}_h} \frac{d_h(\bv_h, p_h)}{\tnorm{\bv_h}_1} \\
	& =  \sup_{\bv_h \in \widehat{\V}_h} \frac{\mu a_h(\bu_h, \bv_h) +\alpha (\bu_h, \bv_h)+b_h(\bu_h;\bu_h;\bu_h, \bv_h) +\beta c_{rh}(\bu_h; \bv_h)-\kappa c_{qh}(\bu_h; \bv_h) - \mathcal{F}_h(\bv_h)}{\tnorm{\bv_h}_1} \\
	& \le \mu \tnorm{\bu_h}_1+ \alpha\|\bu_h\|_0 + C_0 \tnorm{\bu_h}_1^2 + \beta C_s \tnorm{\bu_h}_1^r +\kappa C_s |\mathcal{O}|^{\frac{r-q}{(q+1)(r+1)}}\tnorm{\bu_h}_1^q+\|\f\|_{0}\\
	& \le \|\f\|_{0}+C(\mathcal{R}_0).
\end{align*}
The proof of the uniqueness of $p_h$ follows from a similar argument as used in the proof of Theorem~\ref{thm-equivalent}. 
\end{proof}

\section{Error Estimates}\label{sec6}\setcounter{equation}{0} 
This section is devoted to establish error estimates between the pairs $(\bu,p)$ and $(\bu_h, p_h),$ solutions of Problem~\ref{problem-2} and Problem~\ref{eq:prb-discform}, respectively. To do that, we first estimate the terms $\tnorm{\bu_h - \mathbf{R}_h\bu}_1$ and $\|p_h-q_h\|_0,$ where $\mathbf{R}_h\bu$ is Raviart-Thomas type interpolation of $\bu$ and $q_h$ is classical interpolation of the pressure $p$ satisfying
\begin{align} \label{eq:interpol CR}
	d_h(\mathbf{R}_h\bu- \bu, \zeta_h)=0 \ \text{ for all } \ \bu \in \V(h),\  \zeta_h \in Q_h.
\end{align}

\begin{lemma}[Raviart-Thomas interpolation operator \cite{RavTho75, BMR08, BMR06}] \label{lem:Rav-Th}
 The Raviart-Thomas interpolation operator $\mathbf{R}_h$ belongs to $\mathcal{L}(\H^1(\mathcal{O}); \widehat{\V}_h)$ and for all $\bv \in \H^1(\mathcal{O}),$ it satisfies 
\begin{align*}
	& (i) \ \int_K q_h \nabla \cdot (\mathbf{R}_h\bv - \bv) d\boldsymbol{x}=0 \text{ for all } K \in \mathcal{T}_h, \ q_h \in P_{k-1}(K), \\
	& (ii) \ \int_e q_h(\mathbf{R}_h\bv - \bv)\cdot \boldsymbol{n}  d\mathrm{S}=0 \text{ for all } e \in \mathscr{E}_h, \ q_h\in P_{k-1}(e), \\
	& (iii) \ \| \mathbf{R}_h \bv - \bv\|_{0,K} + h_K \|\nabla (\mathbf{R}_h\bv -\bv)\|_{0,K} \le C h_K \|\bv\|_{\H^1}, \\
	& (iv) \ \mathbf{R}_h\bv|_e\cdot \boldsymbol{n}\in P_{k-1}(e) \text{ for all } e \in \mathscr{E}_h, \text{ and } \tnorm{\mathbf{R}_h\bv}_1\le C \|\nabla \bv\|_0,
\end{align*} 
for some positive constant $C$ independent of $h_K.$
\end{lemma}

\begin{proposition} \label{pps-interim est}
Let $(\bu,p)\in \V\times Q$ and $(\bu_h, p_h)\in \V_h\times Q_h$ be the solutions of Problems~\ref{problem-2} and  \ref{eq:prb-discform}, as discussed in Theorems~\ref{thm-equivalent} and \ref{thm-dgstab}, respectively. 
Assume $\bu\in\H^2(\mathcal{O})$, $\bu_{\boldsymbol{\tau}}\in\H^2(\Gamma_S)$, and $p\in H^1(\mathcal{O})$. Further, assume \eqref{eq:thmassump1} along with 
\begin{align} \label{eq:thmassump2}
\beta rC_p \mathcal{K}^{r-1}+\kappa q C_p \mathcal{K}^{q-1} \le \frac{\mu\gamma_0}{4},
\end{align}
where $C_p$ is the Poincar\'{e} constant and $\mathcal{K}$ is as in \eqref{calK}. 
Then, for any $(\bv_h, q_h) \in \mathcal{D}_h \times Q_h,$ there exists a constant $C>0,$ independent of $h,$  such that 
\begin{align*}
	\tnorm{\bu_h-\bv_h}_1 + \|p_h-q_h\|_0 \le C \left( \tnorm{\bu -\bv_h}_1 + \|q_h-p\|_0 +  \|\bu_{\boldsymbol{\tau}}-\bv_{h\boldsymbol{\tau}}\|_{0,\Gamma_S}^{1/2}\right).
\end{align*}
\end{proposition}

\begin{proof}
We start with the observation that
\begin{equation} \label{eq:split of ah}
\begin{aligned}
	\mu a_h(\bu_h-\bv_h, \bu_h-\bv_h) & = \mu a_h(\bu_h-\bu, \bu_h-\bv_h) + \mu a_h(\bu-\bv_h, \bu_h-\bv_h) \\
	& = \mu\Big( a_h(\bu_h, \bu_h-\bv_h) - a_h(\bu, \bu_h-\bv_h)  \Big) + \mu  a_h(\bu-\bv_h, \bu_h-\bv_h) .
\end{aligned} 
\end{equation}
Note that Problem~\ref{eq:prb-discform} implies
	\begin{equation} \label{eq:d-u-v-l}
		\begin{aligned}
			 \mu a_h(\bu_h  ,\bu_h-\bv_h) & \le   b_h(\bu_h; \bu_h;\bu_h, \bv_h -\bu_h) +\alpha (\bu_h, \bv_h-\bu_h)+\beta c_{rh}(\bu_h; \bv_h-\bu_h) \\ & \quad  - \kappa c_{qh}(\bu_h; \bv_h-\bu_h) + d_h(\bv_h -\bu_h, p_h)+ j(\bv_{h\tau}) - j(\bu_{h\tau})+ \mathcal{F}_h(\bu_h-\bv_h) .
		\end{aligned}
	\end{equation}
Choose $\bv_h=\bu_h$ in \eqref{eqn-problem-2-constnc} to obtain
\begin{equation*}
	\begin{aligned}
		\mu a_h(\bu,\bu_h-\bu)+b_h(\bu,\bu,\bu_h-\bu)+\alpha(\bu,\bu_h-\bu)&+\beta c_{rh}(\bu,\bu_h-\bu)-\kappa c_{qh}(\bu,\bu_h-\bu)\\+d_h(\bu_h-\bu,p)+j(\bu_{h\boldsymbol{\tau}})-j(\bu_{\boldsymbol{\tau}})&\geq \langle\f, \bu_h-\bu\rangle,
	\end{aligned}
\end{equation*}
and then take $\bv_h=2\bu - \bv_h$ in \eqref{eqn-problem-2-constnc} to get
\begin{equation*}
	\begin{aligned}
		\mu a_h(\bu,\bu-\bv_h)+b_h(\bu,\bu,\bu-\bv_h)+\alpha(\bu,\bu-\bv_h)&+\beta c_{rh}(\bu,\bu-\bv_h)-\kappa c_{qh}(\bu,\bu-\bv_h)\\+d_h(\bu-\bv_h,p)+j(2\bu_{\boldsymbol{\tau}}-\boldsymbol{v}_{h\boldsymbol{\tau}})-j(\bu_{\boldsymbol{\tau}})&\geq \langle\f, \bu-\bv_h\rangle. 
	\end{aligned}
\end{equation*}
Thus, a combination of the above two inequalities lead to
\begin{equation} \label{eq:d-u-v-r}
\begin{aligned}
\mu a_h(\bu,\bu_h-\bv_h)+b_h(\bu,\bu,\bu_h-\bv_h)+\alpha(\bu,\bu_h-\bv_h)&+\beta c_{rh}(\bu,\bu_h-\bv_h)-\kappa c_{qh}(\bu,\bu_h-\bv_h)\\+d_h(\bu_h-\bv_h,p)+j(2\bu_{\boldsymbol{\tau}}-\boldsymbol{v}_{h\boldsymbol{\tau}})+j(\bu_{h\boldsymbol{\tau}})-2j(\bu_{\boldsymbol{\tau}})&\geq \mathcal{F}_h(\bu_h-\bv_h).
\end{aligned}
\end{equation}
Utilizing \eqref{eq:d-u-v-l} and \eqref{eq:d-u-v-r}, we infer 
\begin{equation} \label{eq:est term1}
\begin{aligned}
 \mu\Big( & a_h(\bu_h,  \bu_h-\bv_h) - a_h(\bu, \bu_h-\bv_h)  \Big)  \\
 & \le \Big( b_h(\bu_h; \bu_h;\bu_h, \bv_h -\bu_h) +\alpha (\bu_h, \bv_h-\bu_h)+\beta c_{rh}(\bu_h; \bv_h-\bu_h) - \kappa c_{qh}(\bu_h; \bv_h-\bu_h)\\ 
 & \quad   + d_h(\bv_h -\bu_h, p_h)+ j(\bv_{h\tau}) - j(\bu_{h\tau})\Big) + \Big( b_h(\bu;\bu;\bu_h-\bv_h)+\alpha(\bu,\bu_h-\bv_h) \\
 & \quad  +\beta c_{rh}(\bu,\bu_h-\bv_h)-\kappa c_{qh}(\bu,\bu_h-\bv_h)+d_h(\bu_h-\bv_h,p)+j(2\bu_{\boldsymbol{\tau}}-\boldsymbol{v}_{h\boldsymbol{\tau}})-j(\bu_{\boldsymbol{\tau}}) \Big)  \\
 & = \alpha (\bu_h-\bu, \bv_h-\bu_h) + \Big( b_h(\bu_h; \bu_h;\bu_h, \bv_h -\bu_h) - b_h(\bu;\bu;\bu,\bv_h -\bu_h) \Big) + d_h (\bv_h-\bu_h, p_h-p) \\
 & \quad + \beta \Big( c_{rh}(\bu_h; \bv_h-\bu_h) - c_h(\bu; \bv_h-\bu_h)\Big) -\kappa \Big( c_{qh}(\bu_h; \bv_h-\bu_h) - c_{qh}(\bu; \bv_h-\bu_h)\Big) \\
 & \quad  +j(2\bu_{\boldsymbol{\tau}}-\boldsymbol{v}_{h\boldsymbol{\tau}})+j(\boldsymbol{v}_{h\boldsymbol{\tau}})-2j(\bu_{\boldsymbol{\tau}}).
\end{aligned}
\end{equation}
We first observe that
\begin{align} \label{eq:err alpha est}
| \alpha (\bu_h-\bu, \bv_h-\bu_h)| \le \varepsilon \tnorm{\bu_h - \bv_h}_1^2 + C_\varepsilon \|\bu-\bv_h\|_0^2,
\end{align}
for any $\varepsilon>0$ and for some $C=C(\varepsilon)>0.$ Now, addition and subtraction implies
\begin{align*}
	b_h(\bu_h; &  \bu_h;\bu_h, \bv_h -\bu_h) - b_h(\bu;\bu;\bu,\bv_h -\bu_h) \\
	& = b_h(\bu_h;  \bu_h; \bu_h-\bv_h, \bv_h -\bu_h) + b_h(\bu_h;  \bu;\bv_h-\bu, \bv_h -\bu_h) \\
	& \quad + b_h(\bu_h;  \bu_h-\bv_h;\bv_h, \bv_h -\bu_h)+ b_h(\bu_h; \bv_h-\bu;\bv_h, \bv_h -\bu_h).
\end{align*}
Therefore, applying Property~\ref{prop:A} and \eqref{eq:thmassump1}, we infer
\begin{equation} \label{eq:err conv est}
\begin{aligned}
\Big\vert b_h(\bu_h; & \bu_h;\bu_h, \bv_h -\bu_h) - b_h(\bu;\bu;\bu,\bv_h -\bu_h) \Big\vert  \\
& \le C_0 (\tnorm{\bu_h}_1+\tnorm{\bv_h}_1) \tnorm{\bu_h-\bv_h}_1^2 + C_0 (\tnorm{\bu}_1+\tnorm{\bv_h}_1) \tnorm{\bu-\bv_h}_1\tnorm{\bv_h-\bu_h}_1 \\
& \le \frac{\gamma_0}{2} \tnorm{\bu_h-\bv_h}_1^2 +\varepsilon \tnorm{\bu_h-\bv_h}_1^2 + C_{\varepsilon, \mathcal{R}_0} \tnorm{\bu - \bv_h}_1^2.
\end{aligned}
\end{equation}
To estimate the term involving $d_h(\cdot, \cdot),$ we observe 
\begin{align*}
d_h & (\bv_h-\bu_h, p_h-p) \\
&= d_h (\bv_h-\bu_h, p_h-q_h) + d_h (\bv_h-\bu_h, q_h-p) \\
&= d_h (\bv_h-\bu, p_h-q_h)+d_h (\bu, p_h-q_h)-d_h (\bu_h, p_h-q_h)+d_h (\bv_h-\bu_h, q_h-p) \\
&  = d_h (\bv_h-\bu, p_h-q_h)+ d_h (\bv_h-\bu_h, q_h-p),
\end{align*}
as 
$d_h (\bu, p_h-q_h)=0$, and $d_h (\bu_h, p_h-q_h)=0$ due to \eqref{eq:interpol CR}, \eqref{eqn-problem-2-constnc}, and Problem~\ref{eq:prb-discform}, respectively. Therefore, we estimate the terms as 
\begin{align} \label{eq:err d_h est}
\Big\vert d_h (\bv_h-\bu_h, q_h-p) \Big\vert \le \tnorm{\bv_h-\bu_h}_1 \|q_h-p\|_0 \le \varepsilon \tnorm{\bv_h-\bu_h}_1^2 + C_\varepsilon \|q_h-p\|_0^2,
\end{align}
and
\begin{align*}
\left\vert d_h (\bv_h-\bu, p_h-q_h)\right\vert \le  \tnorm{\bv_h-\bu}_1\|p_h-q_h\|_0\le \frac{\vartheta}{2} \|p_h-q_h\|_0^2 +C_\vartheta \tnorm{\bv_h-\bu}_1^2.
\end{align*}
The terms involving damping can be estimated as follows: let $F(\boldsymbol{\xi}):=|\boldsymbol{\xi}|^{r-1}\boldsymbol{\xi},$ and then apply mean value theorem to obtain
\begin{equation} \label{eq:err beta est}
\begin{aligned}
\beta \big\vert  & c_{rh}(\bu_h; \bv_h-\bu_h) - c_{rh}(\bu; \bv_h-\bu_h)\big\vert \\
& = \beta \left \vert \left\langle\int_0^1 DF \left(\theta \bu_h +(1-\theta )\bu \right) (\bu_h -\bu)d\theta, \bv_h-\bu_h  \right\rangle \right \vert \\
& \le \beta r \left\langle \int_0^1 \left\vert \theta \bu_h + (1-\theta) \bu \right\vert^{r-1} |\bu_h-\bu|, |\bv_h-\bu_h| \right\rangle \\
& \le \beta r \left( \|\bu_h\|_{\L^{r+1}} + \|\bu\|_{\L^{r+1}} \right)^{r-1} \|\bu_h-\bu\|_{\L^{r+1}}\|\bv_h-\bu_h\|_{\L^{r+1}} \\
& \le \beta r C_p \mathcal{K}^{r-1} \tnorm{\bu_h-\bu}_1\tnorm{\bv_h-\bu_h}_1 \\
& \le \beta r C_p \mathcal{K}^{r-1}\tnorm{\bv_h-\bu_h}_1^2 + \varepsilon \tnorm{\bv_h-\bu_h}_1^2 + C_\varepsilon \tnorm{\bu-\bv_h}_1^2.
\end{aligned}
\end{equation}
Similarly,
\begin{align} \label{eq:err kappa est}
	\kappa \big\vert  & c_{qh}(\bu_h; \bv_h-\bu_h) - c_{qh}(\bu; \bv_h-\bu_h)\big\vert \le \kappa q C_p \mathcal{K}^{q-1}\tnorm{\bv_h-\bu_h}_1^2 + \varepsilon \tnorm{\bv_h-\bu_h}_1^2 + C_\varepsilon \tnorm{\bu-\bv_h}_1^2.
\end{align}
Now, we estimate
\begin{equation} \label{eq:est j term}
\begin{aligned}
	 j(2\bu_{\boldsymbol{\tau}}-\boldsymbol{v}_{h\boldsymbol{\tau}})+j(\boldsymbol{v}_{h\boldsymbol{\tau}})-2j(\bu_{\boldsymbol{\tau}})  & =  \int_{\mathscr{E}_h^{\S}}g|2\bu_{\boldsymbol{\tau}}-\bv_{h\boldsymbol{\tau}}|\d\S + \int_{\mathscr{E}_h^{\S}}g|\bv_{h\boldsymbol{\tau}}|\d\S-2 \int_{\mathscr{E}_h^{\S}}g|\bu_{\boldsymbol{\tau}}|\d\S  \\
	& \le 2 \int_{\mathscr{E}_h^{\S}}g|\bu_{\boldsymbol{\tau}}-\bv_{h\boldsymbol{\tau}}|\d\S  \le 2 \|g\|_{L^2(\Gamma_S)}\|\bu_{\boldsymbol{\tau}}-\bv_{h\boldsymbol{\tau}}\|_{0, \Gamma_S}.
\end{aligned}
\end{equation}
Also, it is easy to see that
\begin{align} \label{eq:est term2}
\vert	\mu  a_h(\bu-\bv_h, \bu_h-\bv_h)\vert \le \varepsilon \tnorm{\bu_h-\bv_h}_1^2 +C_\varepsilon\tnorm{\bu-\bv_h}_1^2. 
\end{align}
Combining the estimates \eqref{eq:err alpha est}--\eqref{eq:est j term} in \eqref{eq:est term1}, and then using it along with \eqref{eq:est term2} in \eqref{eq:split of ah}, we obtain 
\begin{align*}
&	\mu\gamma_0 \tnorm{\bu_h-\bv_h}_1^2 \\& \le \mu a_h (\bu_h-\bv_h, \bu_h-\bv_h) \\
	 & \le \left( 6\varepsilon + \frac{\gamma_0}{4} + \beta rC_p \mathcal{K}^{r-1}+\kappa q C_p \mathcal{K}^{q-1} \right)\tnorm{\bu_h -\bv_h}_1^2 + C_\varepsilon \|\bu -\bv_h\|_0^2 \\
	 & \quad + C_{\varepsilon, \mathcal{R}_0, \vartheta}  \tnorm{\bu -\bv_h}_1^2 +C_\varepsilon \|q_h - p\|_0^2+ 2 \|g\|_{L^2(\Gamma_S)} \|\bu_{\boldsymbol{\tau}}-\bv_{h\boldsymbol{\tau}}\|_{0, \Gamma_S}+\frac{\vartheta}{2} \|p_h-q_h\|_0^2.
\end{align*}
Choosing $\varepsilon= \frac{\mu\gamma_0}{24},$ and since $\beta rC_p \mathcal{K}^{r-1}+\kappa q C_p \mathcal{K}^{q-1} \le \frac{\mu\gamma_0}{4}$ due to \eqref{eq:thmassump2}, we have 
\begin{align*}
\tnorm{\bu_h-\bv_h}_1^2 \le C \left( \tnorm{\bu -\bv_h}_1^2 + \|q_h-p\|_0^2 +  \|\bu_{\boldsymbol{\tau}}-\bv_{h\boldsymbol{\tau}}\|_{0, \Gamma_S}\right)+\frac{\vartheta}{2} \|p_h-q_h\|_0^2.
\end{align*}
\noindent \textbf{Step 2.} In this step, we estimate $\|p_h-q_h\|_0.$ 
Let $\bw_h \in \widehat{\V}_h$ be arbitrary. Choosing $\bv_h = \bu+\bw_h$ in \eqref{eqn-problem-2-constnc},
\begin{equation*}
	\begin{aligned}
		\mu a_h(\bu,\bw_h)+b_h(\bu;\bu;\bu,\bw_h)+\alpha(\bu,\bw_h)+\beta c_{rh}(\bu,\bw_h)-\kappa c_{qh}(\bu,\bw_h)+d_h(\bw_h,p)\geq \langle\f, \bw_h\rangle,
	\end{aligned}
\end{equation*}
and then choosing $\bv = \bu-\bw_h$ in \eqref{eqn-problem-2-constnc}
\begin{equation*}
	\begin{aligned}
		\mu a_h(\bu,\bw_h)+b_h(\bu;\bu;\bu,\bw_h)+\alpha(\bu,\bw_h)+\beta c_{rh}(\bu,\bw_h)-\kappa c_{qh}(\bu,\bw_h)+d_h(\bw_h,p)\leq \langle\f, \bw_h\rangle.
	\end{aligned}
\end{equation*}
Therefore, we find 
\begin{equation*}
	\begin{aligned}
		\mu a_h(\bu,\bw_h)+b_h(\bu;\bu;\bu,\bw_h)+\alpha(\bu,\bw_h)+\beta c_h(\bu,\bw_h)-\kappa c_h^0(\bu,\bw_h)+d_h(\bw_h,p)= \langle\f, \bw_h\rangle,
	\end{aligned}
\end{equation*}
Similarly, from Problem~\ref{eq:prb-discform}, we have
\begin{equation*}
	\begin{aligned}
		& \mu a_h(\bu_h  ,\bw_h)  + b_h(\bu_h; \bu_h;\bu_h, \bw_h) +\alpha (\bu_h, \bw_h)+\beta c_{rh}(\bu_h; \bw_h) \\ & \qquad \quad  - \kappa c_{qh}(\bu_h; \bw_h) + d_h(\bw_h, p_h)= \mathcal{F}_h(\bw_h).
	\end{aligned}
\end{equation*}
A combination of last two equalities leads to 
\begin{align*}
&\mu a_h(\bu  -\bu_h  ,\bw_h)  + \left( b_h(\bu;\bu;\bu,\bw_h) - b_h(\bu_h; \bu_h;\bu_h, \bw_h)\right) + \alpha (\bu-\bu_h, \bw_h) \\
& \quad +\beta \left( c_{rh}(\bu,\bw_h) - c_{rh}(\bu_h,\bw_h) \right) - \kappa\left( c_{qh}(\bu,\bw_h) - c_{qh}(\bu_h,\bw_h) \right) + d_h (\bw_h, p-p_h)=0,
\end{align*}
for all $\bw_h \in \widehat{\V}_h.$ Utilizing the discrete inf-sup condition \eqref{eqn-dis-inf-sup}, we infer
\begin{align*}
	\vartheta & \|p_h-q_h\|_0 \le \sup_{\bw_h \in \widehat{\V}_h} \frac{d_h(\bw_h, p_h-q_h)}{\tnorm{\bw_h}_1} \\
	&\le \sup_{\bw_h \in \widehat{\V}_h} \frac{d_h(\bw_h, p_h-p) + d_h(\bw_h, p-q_h)}{\tnorm{\bw_h}_1} \\
	& =  \sup_{\bw_h \in \widehat{\V}_h} \frac{  \begin{aligned} &   \mu a_h(\bu  -\bu_h  ,\bw_h)  + \left( b_h(\bu;\bu;\bu,\bw_h) - b_h(\bu_h; \bu_h;\bu_h, \bw_h)\right) + \alpha (\bu-\bu_h, \bw_h) \\ &  +\beta \left( c_{rh}(\bu,\bw_h) - c_{rh}(\bu_h,\bw_h) \right) - \kappa\left( c_{qh}(\bu,\bw_h) - c_{qh}(\bu_h,\bw_h) \right)       \end{aligned}  }{\tnorm{\bw_h}_1}.
\end{align*}
Now, estimating the terms as performed in the previous step, one can easily deduce that 
\begin{align*}
\vartheta \|p_h-q_h\|_0 \le C \left( \mu +\beta r \mathcal{K}^{r-1} +\kappa q \mathcal{K}^{q-1} \right) \tnorm{\bu - \bu_h}_1 +\frac{\vartheta}{2} \|p_h-q_h\|_0+ C \|p-q_h\|_0.
\end{align*}
This estimate, together with Step 1, completes the proof.
\end{proof}

Finally, we are ready to state and prove the main result on error estimate. 

\begin{theorem} \label{thm-dgerr}
Let the assumptions \eqref{eq:thmassump1} and \eqref{eq:thmassump2} hold. Consider the lowest-order case $k=1$ for the finite element spaces $\V_h$ and $Q_h$ defined in \eqref{eq:Vh}--\eqref{eq:Qh}. Let $(\bu,p)\in \V\times Q$ and $(\bu_h, p_h)\in \V_h\times Q_h$ be the solutions of Problems~\ref{problem-2} and  \ref{eq:prb-discform}, respectively. 
Assume $\bu\in\H^2(\mathcal{O})$, $\bu_{\boldsymbol{\tau}}\in\H^2(\Gamma_S)$, and $p\in H^1(\mathcal{O})$. Then there exists a positive constant independent of $h$ such that
	\begin{align*}
		\tnorm{\bu-\bu_h}_1 + \|p-p_h\|_0 \le C h\left( \|\bu\|_{\H^2} +\|p\|_{H^1} + \|\bu_{\boldsymbol{\tau}}\|_{\H^2(\Gamma_S)}^{1/2} \right).
	\end{align*}
\end{theorem}

\begin{proof}
	The result follows from Proposition \ref{pps-interim est} using $(\mathbb{R}_h\bu, q_h)$ as Crouziex-Raviart type interpolation of $\bu$ and $q_h$ is classical interpolation of the pressure $p,$ and triangle inequality
	\begin{align*}
		\tnorm{\bu-\bu_h}_1 + \|p-p_h\|_0 & \le \tnorm{\bu-\mathbb{R}_h\bu}_1+\tnorm{\bu_h-\mathbb{R}_h\bu}_1 + \|p-q_h\|_0+\|q_h-p_h\|_0\\
		& \le C \left( \tnorm{\bu - \mathbb{R}_h\bu}_1 + \|q_h-p\|_0 +  \|\bu_{\boldsymbol{\tau}}-\mathbb{R}_h\bu_{\boldsymbol{\tau}}\|_{0, \Gamma_S}^{1/2}\right) \\
		& \le C h (\|\bu\|_{\H^2} + \|p\|_{H^1} + \|\bu_{\boldsymbol{\tau}}\|_{\H^2(\Gamma_S)}^{1/2}),
	\end{align*}
	where the last inequality is obtained using the properties of Raviart-Thomas interpolation operator as stated in Lemma \ref{lem:Rav-Th}. 
\end{proof}

\section{Numerical Examples}\label{sec7}\setcounter{equation}{0} 
In this section, we present the results of several numerical experiments that illustrate the theoretical findings established in the preceding sections. The computations are performed using an Uzawa-type algorithm, which is given in Algorithm \ref{alg:uzawa_picard_cbfed}.

\begin{algorithm}
	\caption{DG--Uzawa--Picard algorithm for the stationary CBFeD problem}
	\label{alg:uzawa_picard_cbfed}
	\begin{algorithmic}[1]
		\STATE Choose tolerances $\varepsilon_{\mathrm{in}}>0$ and $\varepsilon_{\mathrm{out}}>0$, Uzawa parameter $\rho>0$, and maximum inner and outer iteration numbers $\ell_{\max}$ and $n_{\max}$.
		
		\STATE Select an initial velocity guess $
		\bu_h^{0}\in \mathbb{V}_h,$
		and an initial multiplier $
		\lambda_h^{0}\in \Lambda,$ where
		\begin{align*}
			\Lambda=\{\mu\in L^2(\Gamma_S): |\mu(\boldsymbol{x})|\le 1 \ \text{a.e. on } \Gamma_S\}.
		\end{align*}
		
		\STATE Set $n=0$.
		
		\WHILE{$n\le n_{\max}$}
		
		\STATE Set $\bu_h^{n,0}=\bu_h^n,$ $\lambda_h^{n,0}=\lambda_h^n,$ and $\ell=0.$
		
		\WHILE{$\ell<\ell_{\max}$}
		
		\STATE Given $(\bu_h^{n, \ell},\lambda_h^{n, \ell})$, solve for
		\begin{align*}
			(\bu_h^{n, \ell+1},p_h^{n, \ell+1})\in \mathbb{V}_h\times Q_h
		\end{align*}
		such that
		\begin{equation*}
			\begin{aligned}
				& \mu a_h(\bu_h^{n,\ell+1}  ,\bv_h)  + b_h(\bu_h^{n,\ell}; \bu_h^{n,\ell+1};\bu_h^n, \bv_h) +\alpha (\bu_h^{n,\ell+1}, \bv_h)+\beta (|\bu_h^{n,\ell}|^{r-1}\bu_h^{n,\ell+1}; \bv_h) \\ & \quad \quad  - \kappa (|\bu_h^{n,\ell}|^{q-1}\bu_h^{n,\ell+1}; \bv_h) + d_h(\bv_h, p_h^{\ell+1})+  \int_{\Gamma_S}
				\lambda_h^{n,\ell} g \bv_{h\tau} \, dS  = \mathcal{F}_h(\bv_h) \quad \text{for all } \bv_h \in \mathbb{V}_h, \\
				& d_h(\bu_h^{n,\ell+1}, q_h)=0 \quad \text{ for all } q_h \in Q_h.
			\end{aligned}
		\end{equation*}
		
		\STATE Update the Lagrange multiplier:
		\begin{align*}
			\lambda_h^{n, \ell+1}
			=
			P_{\Lambda}
			\left(
			\lambda_h^{n,\ell}
			+
			\rho g (\bu_h^{n,\ell+1})_\tau
			\right),
		\end{align*}
		where $	P_{\Lambda}(\mu)=		\max\{-1,\min\{1,\mu\}\}.$
		
		\STATE Compute the relative error
		\begin{align*}
			\mathrm{Err}^{(n,\ell)}_{\mathrm{in}}
			=\tnorm{\bu_h^{n,\ell+1}-\bu_h^{n,\ell}}_{1,h},
		\end{align*}
		where $\varepsilon_0>0$ is a small number preventing division by zero.
		
		\IF{$\mathrm{Err}^{(n,\ell)}_{\mathrm{in}}<\varepsilon_{\mathrm{in}}$}
		
		\STATE Set
		\begin{align*}
			\bu_h^{n+1}=\bu_h^{n, \ell+1}, \quad p_h^{n+1}=p_h^{n,\ell+1}, \quad \lambda_h^{n+1}=\lambda_h^{n,\ell+1},
		\end{align*}
		and exit the inner loop.
		
		\ENDIF
		
		\STATE Set
		\begin{align*}
			\bu_h^{n,\ell}\leftarrow \bu_h^{n,\ell+1},
			\quad
			\lambda_h^{n,\ell}\leftarrow \lambda_h^{n,\ell+1}, \text{ and } \ell\leftarrow \ell+1.
		\end{align*}
		\ENDWHILE
		
		\STATE Compute the outer error
		\begin{align*}
			\mathrm{Err}^n_{\mathrm{out}}= \tnorm{\bu_h^{n+1} - \bu_h^n}_{1,h}.
		\end{align*}
		
		\IF{$\mathrm{Err}^n_{\mathrm{out}}<\varepsilon_{\mathrm{out}}$}
		
		\STATE Set computed solution $\bu_h=\bu_h^{n+1},$ $p_h=p_h^{n+1},$ $\lambda_h=\lambda_h^{n+1},$ and stop.
		
		\ENDIF
		
		\STATE $n\leftarrow n+1$.
		\ENDWHILE
	\end{algorithmic}
\end{algorithm}

The following data are used in the numerical computations:

\begin{table}[H]
	\centering
	\setlength{\heavyrulewidth}{0.12em}
	\begin{tabular*}{\textwidth}{@{\extracolsep{\fill}}ccccccccc@{}}
		\toprule
		Example & $\mu$ & $\alpha$ & $\beta$ & $r$ & $\kappa$ & $q$ & $\rho$ & $\gamma$ \\
		\midrule
		1 & 0.5 & 1.0 & 2.0 & 3.0 & 1.0 &  2.0 & 1.0 & 3.0 \\
		2 & 1.5 & 0.8 & 1.2 & 4.0 & 0.2 & 3.0 & 2.0 & 4.0 \\
		\bottomrule
	\end{tabular*}
	\caption{Setup for numerical experiments.}
	\label{tab:setup}
\end{table}

\begin{example}
	In this example, we consider the unit square domain $\mathcal{O}=(0,1) \times (0,1)$ with slip boundary as $\Gamma_{\mathrm S}=\{(x_1,x_2): 0<x_1<1, \ x_2=1\}$ and rest of the part boundary is considered as homogeneous Dirichlet boundary. The parameters are chosen as given in Table~\ref{tab:setup}. Here, we choose the given force function to be computed using
	\begin{align*}
	\f=	-\mu\Delta\bu^0 + (\bu^0 \cdot \nabla)\bu^0 + \alpha\bu^0 + \beta|\bu^0|^{r-1}\bu^0 - \kappa|\bu^0|^{q-1}\bu^0 + \nabla p^0,
	\end{align*}
	where 
	\begin{align*}
		& \bu^0(x_1,x_2)=(8x_1^2(x_1-1)^2x_2(x_2-1)(2x_2-1), -8x_1(x_1-1)(2x_1-1)x_2^2(x_2-1)^2), \\
		& p^0(x_1,x_2)=8(2x_1-1)(2x_2-1).
	\end{align*}
	For the case of homogeneous Dirichlet boundary problem, note that $(\bu^0, p^0)$ satisfies the equations \eqref{eqn-CBF} with the chosen $\f$ above, however it need not be an exact solution when we consider a slip boundary. On $\Gamma_{\mathrm{S}}$ unit normal vector is $\boldsymbol{n}=(0,1),$ and using equation \eqref{eqn-boundary-0} for tangential stress $\boldsymbol{\sigma}_{\boldsymbol{\tau}}$,
a direct computation shows that 
	\begin{align*}
		\max_{\overline{\Gamma_{\mathrm S}}} |\boldsymbol{\sigma}_{\boldsymbol{\tau}}| = \max_{0 \le x_1 \le 1} |8\mu x_1^2 (x_1-1)^2 | = 0.5\mu=0.25.
	\end{align*}
	Here, $\boldsymbol{\tau}$ represents the unit tangential vector obtained by rotating unit normal vector $\boldsymbol{n}$ counterclockwise for an angle of $\frac{\pi}{2}$ radians. Thus from \eqref{eqn-boundary}--\eqref{eqn-boundary-1}, it can be observed that, if $g(\boldsymbol{x})>0.25$ for all $\boldsymbol{x}\in \Gamma_{\mathrm S},$ then no-slip occurs and hence in this case $(\bu^0, p^0)$ is the solution of \eqref{eqn-CBF} with homogeneous Dirichlet boundary conditions, that is, with $\bu^0=\boldsymbol{0}$ on $\Gamma$. Otherwise, slip occurs. To perform computation, we choose different $g$ and the corresponding outcomes are provided here. First, we take $g=1.0>0.25,$ and as discussed, for this chosen $g,$ no-slip occurs, and hence $(\bu^0, p^0)$ is the exact solution of \eqref{eqn-CBF} with $\bu^0=\boldsymbol{0}$ on $\Gamma$. The computed velocity $\bu_h$, its magnitude, the tangential velocity along the slip boundary $x_2=1$, and the computed pressure are shown in the top row of Figure~\ref{fig:ex-1_slip_results}. To compute the errors, we consider the computed solution at mesh size $128 \times 128$ as reference solution $(\bu_h^\ast, p_h^\ast).$ The computed errors and orders of convergence with $P_1$ polynomials are tabulated in Table~\ref{tab:ex-1_dg_combined}. The figures and the table demonstrate the expected outcomes of the numerical experiments, thereby supporting the theoretical results. Furthermore, we consider $g=0.1$, $g=0.01$, and 	$g=0.1(\sin(x_1)+1)$, with slip occurring on $x_2=1$. 	The corresponding results, presented in Figure~\ref{fig:ex-1_slip_results} 	and Table~\ref{tab:ex-1_dg_combined}, support the theoretical results.

\begin{figure}
	\centering
	\includegraphics[width=\textwidth]{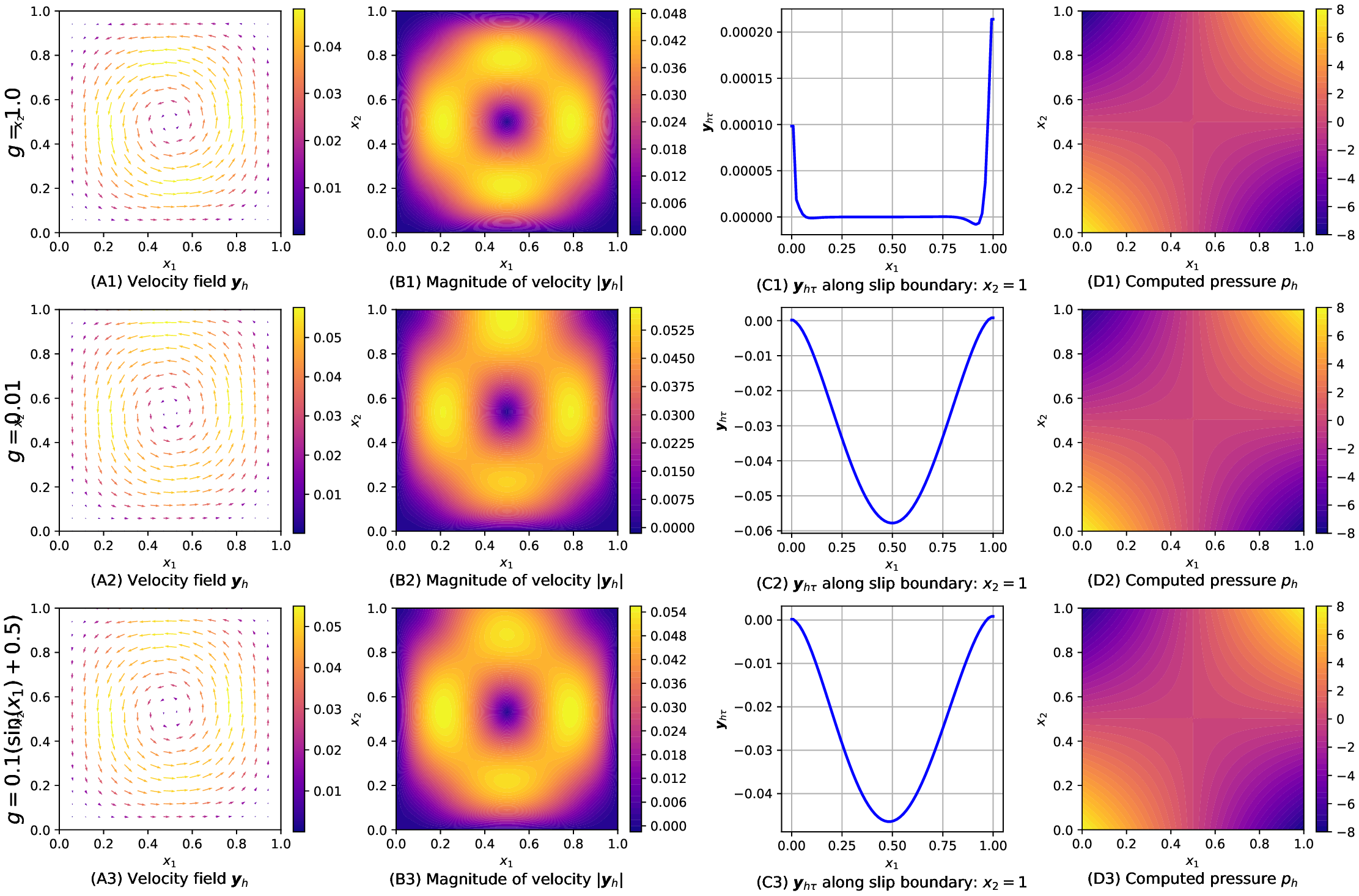}
	\caption{(Example 1) Velocity field, magnitude of velocity, tangential velocity along the slip boundary, and pressure for different slip boundary conditions.}
	\label{fig:ex-1_slip_results}
\end{figure}

\end{example}

\begin{table}[htbp]
	\centering
	\small
	\setlength{\tabcolsep}{5pt}
	\renewcommand{\arraystretch}{1.15}
	\setlength{\heavyrulewidth}{0.12em}
	\begin{tabular*}{\textwidth}{@{\extracolsep{\fill}}c cc cc cc cc@{}}
		\toprule
		& \multicolumn{2}{c}{$g=1$}
		& \multicolumn{2}{c}{$g=0.1$}
		& \multicolumn{2}{c}{$g=0.01$}
		& \multicolumn{2}{c}{$g=0.1(1+\sin x_1)$} \\
		\cmidrule(lr){2-3}
		\cmidrule(lr){4-5}
		\cmidrule(lr){6-7}
		\cmidrule(lr){8-9}
		
		Mesh
		& Error & OC
		& Error & OC
		& Error & OC
		& Error & OC \\
		\midrule
		
		\multicolumn{9}{c}{\textit{Velocity error: $\tnorm{\bu_h^\ast-\bu_h}_1$}} \\
		\midrule
		
		$4\times4$
		& $8.44380\mathrm{e}{-01}$ & --
		& $8.30640\mathrm{e}{-01}$ & --
		& $8.42220\mathrm{e}{-01}$ & --
		& $8.30100\mathrm{e}{-01}$ & -- \\
		
		$8\times8$
		& $4.00360\mathrm{e}{-01}$ & 1.08
		& $3.87010\mathrm{e}{-01}$ & 1.10
		& $3.86360\mathrm{e}{-01}$ & 1.12
		& $3.86500\mathrm{e}{-01}$ & 1.10 \\
		
		$16\times16$
		& $1.87890\mathrm{e}{-01}$ & 1.09
		& $1.81790\mathrm{e}{-01}$ & 1.09
		& $1.80710\mathrm{e}{-01}$ & 1.10
		& $1.81750\mathrm{e}{-01}$ & 1.09 \\
		
		$32\times32$
		& $9.03280\mathrm{e}{-02}$ & 1.06
		& $8.77250\mathrm{e}{-02}$ & 1.05
		& $8.71340\mathrm{e}{-02}$ & 1.05
		& $8.77300\mathrm{e}{-02}$ & 1.05 \\
		
		\midrule
		\multicolumn{9}{c}{\textit{Pressure error: $\|p_h^\ast-p_h\|_0$}} \\
		\midrule
		
		$4\times4$
		& $3.02090\mathrm{e}{-01}$ & --
		& $2.19740\mathrm{e}{-01}$ & --
		& $2.22500\mathrm{e}{-01}$ & --
		& $2.20070\mathrm{e}{-01}$ & -- \\
		
		$8\times8$
		& $1.18060\mathrm{e}{-01}$ & 1.36
		& $9.60180\mathrm{e}{-02}$ & 1.19
		& $9.74220\mathrm{e}{-02}$ & 1.19
		& $9.63620\mathrm{e}{-02}$ & 1.19 \\
		
		$16\times16$
		& $4.53450\mathrm{e}{-02}$ & 1.38
		& $4.22790\mathrm{e}{-02}$ & 1.18
		& $4.26270\mathrm{e}{-02}$ & 1.19
		& $4.23160\mathrm{e}{-02}$ & 1.19 \\
		
		$32\times32$
		& $1.80770\mathrm{e}{-02}$ & 1.33
		& $1.52380\mathrm{e}{-02}$ & 1.47
		& $1.53280\mathrm{e}{-02}$ & 1.48
		& $1.52410\mathrm{e}{-02}$ & 1.47 \\
		\bottomrule
	\end{tabular*}
	\caption{(Example 1) Velocity and pressure errors and observed orders of
		convergence for different slip boundary conditions.}
	\label{tab:ex-1_dg_combined}
\end{table}

\begin{figure}
	\centering
	\includegraphics[width=\textwidth]{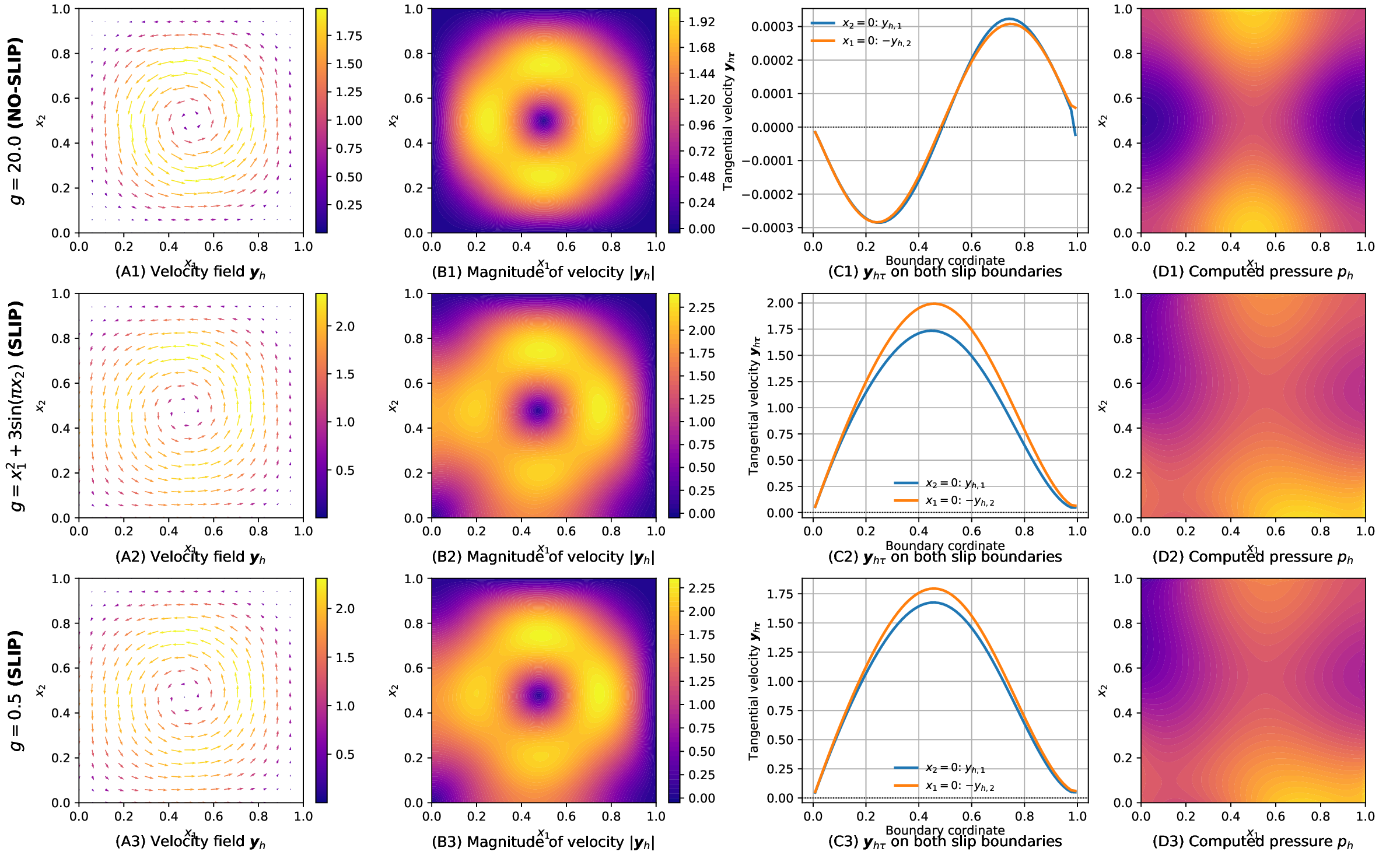}
	\caption{(Example 2) Velocity field, magnitude of velocity, tangential velocity along the slip boundary, and pressure for different slip boundary conditions.}
	\label{fig:ex-2_slip_results}
\end{figure}

\begin{example}
In this example, we consider the unit square domain $\mathcal{O}=(0,1)\times(0,1)$, where the left and bottom boundaries are treated as slip boundaries, that is,  the slip boundary is $\Gamma_{\mathrm S}=\Gamma_{\mathrm S_1}\cup \Gamma_{\mathrm S_2},$ where
\begin{align*}
\Gamma_{\mathrm S_1}=\{(0,x_2):0<x_2<1\}, \text{ and } \Gamma_{\mathrm S_2}=\{(x_1,0):0<x_1<1\},
\end{align*}
while homogeneous Dirichlet boundary conditions are imposed on the remaining part (the top and right boundaries) of the boundary. The model parameters are chosen as listed in Table~\ref{tab:setup}. To assess the accuracy and convergence of the proposed finite element method, we employ the method of manufactured solutions. The body force $\f$ is determined from the prescribed functions $(\bu^0,p^0)$ through
\begin{align*}
	\f
	&= -\mu\Delta\bu^0
	+(\bu^0\cdot\nabla)\bu^0
	+\alpha\bu^0
	+\beta|\bu^0|^{r-1}\bu^0
	-\kappa|\bu^0|^{q-1}\bu^0
	+\nabla p^0,
\end{align*}
where
\begin{align*}
	\bu^0(x_1,x_2)
	&=\bigl(-\cos(2\pi x_1)\sin(2\pi x_2)+\sin(2\pi x_2),\sin(2\pi x_1)\cos(2\pi x_2)-\sin(2\pi x_1)\bigr),\\
	p^0(x_1,x_2)
	&=2\pi\bigl(\cos(2\pi x_2)-\sin(2\pi x_1)\bigr).
\end{align*}
Note that on $\Gamma_{\mathrm S_1},$ the unit outward normal vector is $(-1,0)$ while the unit normal vector on boundary $\Gamma_{\mathrm S_2}$ is $(0,-1).$ Now, again an easy computation leads to 
\begin{align*}
	\max_{\overline{\Gamma_{\mathrm S_1}}} |\boldsymbol{\sigma}_{\boldsymbol{\tau}}| = \max_{0 \le x_2 \le 1} |2\pi\mu (1-\cos(2\pi x_2)) | = 4\pi\mu=6\pi,
\end{align*}
and 
\begin{align*}
	\max_{\overline{\Gamma_{\mathrm S_2}}} |\boldsymbol{\sigma}_{\boldsymbol{\tau}}| = \max_{0 \le x_1 \le 1} |2\pi\mu (\cos(2\pi x_1)-1) | = 4\pi\mu=6\pi.
\end{align*}
To validate the proposed DG scheme, we perform computations for different choices of the slip bound $g$. We first consider the constant slip bound $g=20.0>6\pi$ and hence we expect a no-slip boundary. In particular, $(\bu^0,p^0)$ is exact solution of \eqref{eqn-CBF} with homogeneous Dirichlet boundary condition, that is, with $\bu^0=\boldsymbol{0}$ on $\Gamma.$ The computed velocity field $\bu_h$, its magnitude, the tangential velocity along both $\Gamma_{\mathrm S_1}$ and $\Gamma_{\mathrm S_2}$, and computed pressure $p_h$ are displayed in the first row of Figure~\ref{fig:ex-2_slip_results}. To compute the errors, we consider the computed solution at mesh sizes $128 \times 128$ as reference solution $(\bu_h^\ast, p_h^\ast).$ The corresponding errors and experimental orders of convergence are reported in the first part of Table~\ref{tab:ex-2_dg_combined}. It can be observed that the convergence orders that are in good agreement with the theoretical predictions.

Next, we consider the slip bounds $g=0.5<6\pi$ and $g=x_1^2+3\sin(\pi x_2)<6\pi$, for which slip occurs along the designated slip boundaries. The corresponding numerical solutions are presented in second and third rows of Figure~\ref{fig:ex-2_slip_results}, while the computed errors and convergence orders are summarized in the remaining part of Table~\ref{tab:ex-2_dg_combined}. These numerical results further confirm the theoretical findings and demonstrate the effectiveness and robustness of the proposed DG method under both no-slip and slip regimes.
\end{example}
	
\begin{table}[htbp]
	\centering
	\small
	\setlength{\tabcolsep}{5pt}
	\renewcommand{\arraystretch}{1.15}
	\setlength{\heavyrulewidth}{0.12em}
	\begin{tabular*}{\textwidth}{@{\extracolsep{\fill}}c cc cc cc@{}}
		\toprule
		& \multicolumn{2}{c}{No-slip: $g=20.0$}
		& \multicolumn{2}{c}{Slip: $g=x_1^2+3\sin(\pi x_2)$}
		& \multicolumn{2}{c}{Slip: $g=0.5$} \\
		\cmidrule(lr){2-3}
		\cmidrule(lr){4-5}
		\cmidrule(lr){6-7}
		
		Mesh
		& Error & OC
		& Error & OC
		& Error & OC \\
		\midrule
		
		\multicolumn{7}{c}{\textit{Velocity error: $\tnorm{\bu_h^\ast-\bu_h}_1$}} \\
		\midrule
		
		$4\times4$
		& $3.98446\mathrm{e}{+00}$ & --
		& $3.12613\mathrm{e}{+00}$ & --
		& $3.08205\mathrm{e}{+00}$ & -- \\
		
		$8\times8$
		& $2.10314\mathrm{e}{+00}$ & 0.92
		& $1.80632\mathrm{e}{+00}$ & 0.79
		& $1.77809\mathrm{e}{+00}$ & 0.79 \\
		
		$16\times16$
		& $1.08827\mathrm{e}{+00}$ & 0.95
		& $1.00011\mathrm{e}{+00}$ & 0.85
		& $9.83154\mathrm{e}{-01}$ & 0.85 \\
		
		$32\times32$
		& $5.89729\mathrm{e}{-01}$ & 0.88
		& $5.99056\mathrm{e}{-01}$ & 0.74
		& $5.87088\mathrm{e}{-01}$ & 0.74 \\
		
		\midrule
		\multicolumn{7}{c}{\textit{Pressure error: $\|p_h^\ast-p_h\|_0$}} \\
		\midrule
		
		$4\times4$
		& $2.68853\mathrm{e}{+00}$ & --
		& $2.02757\mathrm{e}{+00}$ & --
		& $1.94454\mathrm{e}{+00}$ & -- \\
		
		$8\times8$
		& $1.50688\mathrm{e}{+00}$ & 0.84
		& $1.09849\mathrm{e}{+00}$ & 0.88
		& $1.05199\mathrm{e}{+00}$ & 0.89 \\
		
		$16\times16$
		& $9.39074\mathrm{e}{-01}$ & 0.68
		& $7.82768\mathrm{e}{-01}$ & 0.49
		& $7.58856\mathrm{e}{-01}$ & 0.47 \\
		
		$32\times32$
		& $6.37836\mathrm{e}{-01}$ & 0.56
		& $4.15390\mathrm{e}{-01}$ & 0.91
		& $3.99284\mathrm{e}{-01}$ & 0.93 \\
		\bottomrule
	\end{tabular*}
	\caption{(Example 2) Velocity and pressure errors and observed orders of
		convergence for different boundary conditions.}
	\label{tab:ex-2_dg_combined}
\end{table}

%
%
%
%
%
%
%

\medskip
\noindent
\textbf{Acknowledgments:}  Support for M. T. Mohan's research received from the National Board of Higher Mathematics (NBHM), Department of Atomic Energy, Government of India (Project No. 02011/13/2025/NBHM(R.P)/R\&D II/1137). Wasim Akram is supported by NBHM (National Board of Higher Mathematics, Department of Atomic Energy) postdoctoral fellowship, No. 0204/16(1)(2)/2024/R\&D-II/10823.
	
	\medskip\noindent	{\bf  Declarations:} 

\noindent 	{\bf  Ethical Approval:}   Not applicable 

\noindent  {\bf   Competing interests: } The author declare no competing interests.

\noindent 	{\bf   Funding: } NBHM, India: No. 02011/13/2025/NBHM(R.P)/R\&D II/1137 (M. T. Mohan) and No. 0204/16(1)(2)/2024/R\&D-II/10823 (W. Akram).

\noindent 	{\bf   Availability of data and materials: } Not applicable.


\begin{thebibliography}	{99}	
		
		
		
		
	
	\bibitem{WAMTM1} W. Akram and M.~T. Mohan, \emph{Mixed finite element method for a hemivariational inequality of stationary convective Brinkman--Forchheimer extended Darcy equations,}  J. Sci. Comput. {\bf 108} (2026), no.~1, Paper No. 10, 50 pp.
	
	\bibitem{WAMTM2} W. Akram and M.~T. Mohan, \emph{Optimal control of a hemivariational inequality of stationary convective Brinkman--Forchheimer extended Darcy equations with numerical approximation,} Comput. Methods Appl. Mech. Engrg. {\bf 452} (2026), Part B, Paper No. 118755, 31 pp.
		
		
		
		
		
		
		
		
		
	
	
	
	
		
		
		
		
		
		
		
		
		
		
		
		
		
		
		
		
		
		
		
		
		\bibitem{SCBLR}	S. C. Brenner and L. R. Scott, \emph{The mathematical theory of finite element methods}, Springer-Verlag, New York, 2002.
		
		
		
		
		
		
		
		
		
		
		
		
		
		
		
		
		\bibitem{PGC} P. G. Ciarlet, \emph{Linear and nonlinear functional analysis with applications}, SIAM Philadelphia, Philadelphia (2013).	
		
		\bibitem{PGC1}  P. G. Ciarlet,  \emph{The finite element method for elliptic problems}, North-Holland, Amsterdam, 1978. 
		
		
		
		
		
		
		
		
		
	\bibitem{JKDJK} 	J. K. Djoko and  J. Koko, \emph{Numerical methods for the Stokes and Navier–Stokes equations driven by threshold slip boundary conditions}, Comput. Methods Appl.Mech. Engrg., {\bf  305} (2016), 936--958.
	
		
		
		
		
		
		
		
		
		
		
		
		
		
		
		
		
		
		
		
	\bibitem{HFu} 	H. Fujita, \emph{Flow problems with unilateral boundary conditions}, College de France, Lecons, 1993.
	
		\bibitem{HFu1} H. Fujita, \emph{A mathematical analysis of motions of viscous incompressible fluid under leak or slip boundary conditions}, RIMS K\^oky\^uroku,  888 (1994), 199--216. 
		
		\bibitem{HFuHK}  H. Fujita, H. Kawarada, \emph{Variational inequalities for the Stokes equation with boundary conditions of friction type}, GAKUTO Int. Ser. Math. Sci. Appl., {\bf 11} (1998), 15--33.
		
		\bibitem{HFu2} H. Fujita, \emph{Non-stationary Stokes flows under leak boundary conditions of friction type}, J. Comput. Math. \textbf{19}, 1–-8 (2001).
		
		\bibitem{HFu3} H. Fujita, \emph{A coherent analysis of Stokes flows under boundary conditions of friction type}, J. Comput. Appl. Math. \textbf{149}, 57–69 (2002) 
		
		
		
		
		
		
		
		
		\bibitem{SGMTM} S. Gautam and M. T. Mohan, \emph{On the convective Brinkman--Forchheimer equations}, Dyn. Partial Differ. Equ., 
		{\bf 22}(3) (2025),  191--233. 
		
		 \bibitem{SGKKMTM}  S. Gautam, K. Kinra and M. T. Mohan, \emph{Feedback stabilization of convective Brinkman--Forchheimer extended Darcy equations}, Appl. Math. Optim., {\bf  91} (2025), no. 1, Paper No. 25, 75 pp.
		
		\bibitem{VGPR}  V. Girault and P.-A. Raviart,  \emph{Finite element methods for Navier--Stokes equations, theory and algorithms}, Springer-Verlag, Berlin, 1986. 
		
		\bibitem{VGBR} V.  Girault, V., B. Riviere, \emph{DG approxiamtion of coupled Navier–Stokes and Darcy equations by Beaver– Joseph–Saffman interface condition}, SIAM J. Numer. Anal. \textbf{47}, 2052–2089 (2009)
		
		\bibitem{VGBMRMFW_MathComp} V. Girault, B.~M. Rivi\`ere and M.~F. Wheeler, \emph{A discontinuous Galerkin method with nonoverlapping domain decomposition for the Stokes and Navier--Stokes problems}, Math. Comp. {\bf 74} (2005), no.~249, 53--84; MR2085402
		
		\bibitem{VGBMRMFW} V. Girault, B.~M. Rivi\`ere and M.~F. Wheeler, \emph{A splitting method using discontinuous Galerkin for the transient incompressible Navier--Stokes equations}, M2AN Math. Model. Numer. Anal. {\bf 39}(6) (2005),  1115--1147. 
		
		
		
		
		
		
		
		
		\bibitem{KWH} K. W. Hajduk and J. C. Robinson, \emph{Energy equality for the 3D critical convective Brinkman--Forchheimer equations}, J. Differential Equations, \textbf{263}(11) (2017), 7141--7161.
		
		
		
		
		
		
	\bibitem{WHKCFJ} W. Han, K. Czuprynski and F. Jing, \emph{Mixed finite element method for a hemivariational inequality of stationary Navier--Stokes equations},	J. Sci. Comput., {\bf 89}(1) (2021), Paper No. 8, 22 pp.
	
	
	\bibitem{WHFJYY}  	 W. Han, F. Jing and Y.  Yao, \emph{Stabilized low-order mixed finite element methods for a Navier--Stokes hemivariational inequality}, BIT, {\bf 63}(4) (2023),  Paper No. 46, 31 pp.
		
	\bibitem{WHSM1} W. Han, S. Mig\'orski, and M. Sofonea, \emph{A class of variational-hemivariational inequalities with applicationsto frictional contact problems,} SIAM J. Math. Anal., {\bf 46} (2014), 3891--3912. 
		
		
	\bibitem{WHSM} W. Han, S. Mig\'orski and M. Sofonea, \emph{Advances in variational and hemivariational inequalities: theory, numerical analysis, and applications}, Springer, Cham, 2015. 
		
		
		
	\bibitem{WHHQLM}  W. Han, H. Qiu and L. Mei, 
	\emph{On a Stokes hemivariational inequality for incompressible fluid flows with damping},
	Nonlinear Anal. Real World Appl., {\bf 79} (2024), Paper No. 104131, 16 pp.
		
		
		
	\bibitem{WHMS} W. Han and M. Sofonea, \emph{Numerical analysis of hemivariational inequalities in contact mechanics}, Acta Numer., {\bf 28}(2019), 175--286. 
		
	\bibitem{WHMS1}  W. Han, M. Sofonea and M. Barboteu,  \emph{Numerical analysis of elliptic hemivariational inequalities}, SIAMJ. Numer. Anal., {\bf 55}(2017), 640--663 (2017)
		
	\bibitem{WHMS2}  W. Han,  M. Sofonea and D. Danan, \emph{Numerical analysis of stationary variational-hemivariational inequalities,} Numer. Math., {\bf 139}  (2018), 563--592. 
		
	\bibitem{WHCW} W. Han,  and C. Wang, \emph{Numerical analysis of a parabolic hemivariational inequality for semipermeablemedia,} J. Comput. Appl. Math., {\bf 389} (2021), 113326.
		
	
		
		
	\bibitem{JSHNJW} 	J. S. Howell,  and  N. J. Walkington, \emph{Inf-sup conditions for twofold saddle point problems,}	Numer. Math., {\bf  118}(4) (2011),  663--693.
		
		
		
		
	
	\bibitem{JJSGMTM1} J. Jindal, S. Gautam, and  M. T. Mohan, \emph{A domain hemivariational inequality for 2D and 3D convective Brinkman--Forchheimer extended Darcy equations,} submitted, (2026), \href{https://doi.org/10.48550/arXiv.2603.28051}{https://doi.org/10.48550/arXiv.2603.28051}.
	
	\bibitem{JJSGMTM2} J. Jindal, S. Gautam, and  M. T. Mohan, \emph{Well-posedness of a boundary hemivariational inequality for stationary and non-stationary 2D and 3D convective Brinkman--Forchheimer equations,} submitted, (2026), \href{https://doi.org/10.48550/arXiv.2508.17093}{https://doi.org/10.48550/arXiv.2508.17093}.
	
	\bibitem{JJMTM} J. Jindal, and  M. T. Mohan, \emph{A quasi-variational--hemivariational inequality for the Convective Brinkman--Forchheimer extended Darcy equations with Bingham fluids,} submitted, (2026), \href{https://doi.org/10.48550/arXiv.2606.27283}{https://doi.org/10.48550/arXiv.2606.27283}.
	
		
	\bibitem{FJWHWY}  F. Jing, W. Han, W. Yan and F. Wang, \emph{Discontinuous Galerkin methods for a stationary Navier--Stokes problem with a nonlinear slip boundary condition of friction type}, \emph{J. Sci. Comput.}, {\bf  76} (2018), 888--912.
			
				
	\bibitem{FJWHGZ} F. Jing,  W. Han and G. Zhou, \emph{Numerical analysis of a class of hemivariational inequalities governed by fluid-fluid coupled flow}, Nonlinear Anal. Real World Appl., {\bf 85} (2025), Paper No. 104366, 20 pp.
				
			
	
	

		
		
		
		
		
		
	\bibitem{TKa}	T. Kashiwabara, \emph{On a finite element approximation of the Stokes equations under a slip boundary condition of the friction type}, Jpn. J. Indust. Appl.Math., {\bf  30} (2013), 227--261.
		
	\bibitem{TKa1}	T. Kashiwabara, \emph{Finite element method for Stokes equations under leak boundary condition of friction type}, SIAM J. Numer. Anal., {\bf 51} (2013), 2445--2469.
		
	
	\bibitem{KKMTM-RMP} K. Kinra and M.~T. Mohan, \emph{Asymptotic analysis of the 2D and 3D convective Brinkman--Forchheimer equations in unbounded domains}, Rev. Math. Phys. {\bf 37} (2025), no.~6, Paper No. 2450044, 52 pp.; MR4929736
	
	\bibitem{KKMTM-2025} K. Kinra and M.T. Mohan, \emph{Stochastic convective Brinkman--Forchheimer equations on general unbounded domains}, Submitted, 2025, \url{https://arxiv.org/pdf/2007.09376}. 
		
		
		 
		 
		  
		
		
	\bibitem{PLPAR}	P. Lesaint and P.-A. Raviart, \emph{On a finite element method for solving the neutron transport equation, in Mathematical aspects of finite elements in partial differential equations} (Proc. Sympos., Math. Res. Center, Univ. Wisconsin, Madison, Wis., 1974), pp. 89--123, Academic Press, New York-London.
		
		
		
		
		
		
		
		
		
		
		
		
		
		
	\bibitem{YLKL}  Y. Li and  K. Li, \emph{Uzawa iteration method for Stokes type variational inequality of the second kind}, Acta Math. App. Sin., {\bf 17} (2011) 303--316.
		
		
	\bibitem{YLRA} Y. Li and R. An, \emph{Semi-discrete stabilized finite element methods for Navier--Stokes equations with nonlinear slip boundary conditions based on regularization procedure}, Numer. Math., {\bf 117} (2011), 1--36.
		
	\bibitem{YLRA1}  Y. Li and  R. An, \emph{Penalty finite element method for Navier--Stokes equations with nonlinear slip boundary conditions}, Internat. J. Numer. Methods Fluids, {\bf  69}(2012), 550--566.
		
	\bibitem{YLRA2} Y. Li and  R. An, \emph{Two-level variational multiscale finite element methods for Navier--Stokes type variational inequality problem}, J. Comput. Appl. Math., {\bf 290} (2015), 656--669. 
		
		
	\bibitem{YLKTL} Y. Li and K. T. Li, \emph{Existence of the solution to stationary Navier–Stokes equations with nonlinear slip boundary conditions}, J. Math. Anal. Appl. \textbf{381}, 1–-9 (2011).
	
	
		
	
	
	
\bibitem{MLWHSZ}	M. Ling, W. Han, and S. Zeng, \emph{A pressure projection stabilized mixed finite element method for a Stokes hemivariational inequality}, J. Sci. Comput., {\bf 92}(1) (2022), Paper No. 13, 21 pp.
		
		
	\bibitem{MLFW}  M. Ling, F. Wang and W. Han, \emph{The nonconforming virtual element method for a stationary Stokes hemivariational inequality with slip boundary condition}, J. Sci. Comput., {\bf 85}(2020), Ar. No. 56. 
		
	\bibitem{JLL}  J.-L. Lions. \emph{Quelques m\'ethodes de r\'esolution des problèmes aux limites non lin\'eaires},  Dunod, Paris, Gauthier-Villars, Paris, 1969.
		
		
	\bibitem{MTT} P. A. Markowich, E. S. Titi and S. Trabelsi, \emph{Continuous data assimilation for the three dimensional Brinkman--Forchheimer-extended Darcy model}, Nonlinearity, \textbf{29}(4) (2016), 1292--1328. 
		
		
		
	
	
\bibitem{SMSDYZ}	S. Mig\'orski, S. Dudek and Y. Zhang, Convergence of a novel penalty method for quasi variational-hemivariational inequality with application, NoDEA Nonlinear Differential Equations Appl. {\bf 33} (2026), no.~2, Paper No. 39, 25 pp.
		
	\bibitem{SMAOMS}   S. Mig\'orski, A. Ochal, M. Sofonea, \emph{Nonlinear inclusions and hemivariational inequalities. models and analysis of contact problems}, in: Advances in Mechanics and Mathematics, vol. 26, Springer, New York, 2013.
		
		
		
		
		
	\bibitem{MT2}  M. T. Mohan, \emph{Well-posedness and asymptotic behavior of stochastic convective Brinkman–Forchheimer equations perturbed by pure jump noise}, Stoch PDE: Anal. Comp., {\bf 10}(2) (2022), 614--690.
		
		
		
		
	
	\bibitem{MTM9} M.~T. Mohan, \emph{Large deviations for the two-time-scale stochastic convective Brinkman--Forchheimer equations}, J. Differential Equations {\bf 376} (2023), 469--537; MR4642545
	
	\bibitem{MTM7} M. T. Mohan, \emph{Existence and uniqueness of time-periodic solutions of the 2D and 3D convective Brinkman--Forchheimer extended Darcy equations}, Z. Angew. Math. Phys. {\bf 77} (2026), no.~4, Paper No. 110;
	
	\bibitem{MTM8} M.~T. Mohan, \emph{Well-posedness of stationary 2D and 3D Convective Brinkman--Forchheimer extended Darcy hemivariational inequalities,} J. Optim. Theory Appl. {\bf 208} (2026), no.~3, Paper No. 118; 
		
		
	\bibitem{ZNPDP} Z.  Naniewicz, and  P. D. Panagiotopoulos,  \emph{Mathematical theory of hemivariational inequalities and applications}, Dekker, New York, 1995. 
		
		
		
		
			
			
			
	\bibitem{PDP}  	P. D. Panagiotopoulos, \emph{Nonconvex energy functions, hemivariational inequalities and substationary principles,} \emph{Acta Mech.}, {\bf 42} (1983), 160--183. 
		
	\bibitem{PDP1} P. D.  Panagiotopoulos, \emph{Hemivariational inequalities: applications in mechanics and engineering}, Springer, Berlin, 1993. 
		
		
		
		
	
	\bibitem{PErnDG} D.~A. Di~Pietro and A. Ern, {\it Mathematical aspects of discontinuous Galerkin methods}, Math\'ematiques \& Applications (Berlin), 69, Springer, Heidelberg, 2012;
		
		
		
	
	\bibitem{RavTho75}  P.-A. Raviart and J.-M. Thomas, \emph{A mixed finite element method for 2nd order elliptic problems}, in {\it Mathematical aspects of finite element methods} (Proc. Conf., Consiglio Naz. delle Ricerche (C.N.R.), Rome, 1975), pp. 292--315, Lecture Notes in Math., Vol. 606, Springer, Berlin-New York, ; MR0483555
		
		
		
		
		
		
	\bibitem{BMR08} B.~M. Rivi\`ere, \emph{ Discontinuous Galerkin methods for solving elliptic and parabolic equations}, Frontiers in Applied Mathematics, 35, SIAM, Philadelphia, PA, 2008; 
	
	\bibitem{BMR06} B.~M. Rivi\`ere and V. Girault, \emph{Discontinuous finite element methods for incompressible flows on subdomains with non-matching interfaces}, Comput. Methods Appl. Mech. Engrg. {\bf 195} (2006), no.~25-28, 3274--3292;
	
	\bibitem{JCBob3D} J. C. Robinson, J. L. Rodrigo, and W. Sadowski. \emph{The three-dimensional Navier--Stokes equations: Classical theory}, volume 157. Cambridge university press, 2016.
	
		
		
	\bibitem{NSa} N. Saito,  \emph{On the Stokes equation with the leak or slip boundary conditions of friction type: regularity ofsolutions}, Publ. Res. Inst. Math. Sci., {\bf 40} (2004), 345--383. 
		
		
	
		
		
		
		
		
		
				
	
		
		
	
		
	\bibitem{Te} R. Temam, \emph{Navier--Stokes equations: theory and numerical analysis,} North-Holland, Amsterdam, 1984.
		
		
	\bibitem{Te1}  R. Temam, \emph{Infinite-dimensional dynamical systems in mechanics and physics}, Second edition, Vol. 68, Applied Mathematical Sciences, Springer, 1997.
		
		
		
		
		\bibitem{FWMLWH} F. Wang,  M. Ling,  W. Han and F. Jing, \emph{Adaptive discontinuous Galerkin methods for solving an incompressible Stokes flow problem with slip boundary condition of frictional type}, J. Comput. Appl. Math., {\bf 371} (2020), 112700. 
		
		
		\bibitem{WWXLCWH} W. Wang, X.-L. Cheng and W. Han, \emph{Analysis and finite element solution of a Navier--Stokes hemivariational inequality for incompressible fluid flows with damping}, Nonlinear Anal. Real World Appl. {\bf 87} (2026), Paper No. 104439, 16 pp.
		
		
		
		
		
		
		
		
	\end{thebibliography}
\end{document}